\documentclass[twoside]{article}
\usepackage{atmp04}
\usepackage{latexsym}
\copyrightnotice{2002}{6}{457}{505} 
\renewcommand{\title}[1] {%
  \begingroup
    \begin{center}
      \vspace{0.4in}
      \bf\huge
      \addtolength{\baselineskip}{5mm}
      #1
    \end{center}
  \endgroup
}

\renewcommand{\author}[1] {%
  \begingroup
    \begin{center}
      \vspace{0.4in}
      \bf
      #1
      \vspace{0.2in}
    \end{center}
  \endgroup
}

\catcode`\@=11
\def\newsymbol#1#2#3#4#5{\let\next@\relax%
 \ifnum#2=\@ne\else%
 \ifnum#2=\tw@\let\next@\msyfam@\fi\fi%
 \mathchardef#1="#3\next@#4#5}
\def\mathhexbox@#1#2#3{\relax%
 \ifmmode\mathpalette{}{\m@th\mathchar"#1#2#3}r
 \else\leavevmode\hbox{$\m@th\mathchar"#1#2#3$}\fi}
\def\hexnumber@#1{\ifcase#1 0\or 1\or 2\or 3\or 4\or 5\or 6\or 7\or 8%
\or 9\or A\or B\or C\or D\or E\or F\fi}
\font\tenmsy=msbm10
\font\sevenmsy=msbm7
\font\fivemsy=msbm5
\newfam\msyfam
\textfont\msyfam=\tenmsy
\scriptfont\msyfam=\sevenmsy
\scriptscriptfont\msyfam=\fivemsy
\edef\msyfam@{\hexnumber@\msyfam}
\def\Bbb#1{\fam\msyfam\relax#1}
\catcode`\@=\active

\catcode`\@=11\relax
\newwrite\@unused
\def\typeout#1{{\let\protect\string\immediate\write\@unused{#1}}}

\def\psglobal#1{
\immediate\special{ps: plotfile #1 }}
\def\psfiginit{\typeout{psfiginit}
\immediate\psglobal{figtex.pro}%
\special{ps:: /TeXMagnification {\the\mag} def}
}

\def\@nnil{\@nil}
\def\@empty{}
\def\@psdonoop#1\@@#2#3{}
\def\@psdo#1:=#2\do#3{\edef\@psdotmp{#2}\ifx\@psdotmp\@empty \else
    \expandafter\@psdoloop#2,\@nil,\@nil\@@#1{#3}\fi}
\def\@psdoloop#1,#2,#3\@@#4#5{\def#4{#1}\ifx #4\@nnil \else
       #5\def#4{#2}\ifx #4\@nnil \else#5\@ipsdoloop #3\@@#4{#5}\fi\fi}
\def\@ipsdoloop#1,#2\@@#3#4{\def#3{#1}\ifx #3\@nnil
       \let\@nextwhile=\@psdonoop \else
      #4\relax\let\@nextwhile=\@ipsdoloop\fi\@nextwhile#2\@@#3{#4}}
\def\@tpsdo#1:=#2\do#3{\xdef\@psdotmp{#2}\ifx\@psdotmp\@empty \else
    \@tpsdoloop#2\@nil\@nil\@@#1{#3}\fi}
\def\@tpsdoloop#1#2\@@#3#4{\def#3{#1}\ifx #3\@nnil
       \let\@nextwhile=\@psdonoop \else
      #4\relax\let\@nextwhile=\@tpsdoloop\fi\@nextwhile#2\@@#3{#4}}
\def\psdraft{
	\def\@psdraft{0}
	\def\@psdraftspecial{100}
}
\def\psdraftspecial{
	\def\@psdraft{0}
	\def\@psdraftspecial{0}
}
\def\psfull{
	\def\@psdraft{100}
}
\psfull

\newif\if@prologfile
\newif\if@postlogfile
\newif\if@bbllx
\newif\if@bblly
\newif\if@bburx
\newif\if@bbury
\newif\if@height
\newif\if@width
\newif\if@rheight
\newif\if@rwidth
\newif\if@clip
\newif\if@right
\newif\if@left
\newif\if@toplines
\newif\if@box
\newif\if@caption
\newif\if@surround
\newif\if@captionwidth
\newif\if@captionwrite
\newif\if@captionopen
\def\@p@@sclip#1{\@cliptrue}
\def\@p@@sfile#1{
		\def\@p@sfile{#1}
}
\def\@p@@sfigure#1{
		\def\@p@sfile{#1}
}
\def\@p@sfake{\hbox to 0pt{\hss Whatever\hss}}
\def\@p@@sbbllx#1{
		\@bbllxtrue
		\@d@mscratch=#1
		\edef\@p@sbbllx{\number\@d@mscratch}
}
\def\@p@@sbblly#1{
		\@bbllytrue
		\@d@mscratch=#1
		\edef\@p@sbblly{\number\@d@mscratch}
}
\def\@p@@sbburx#1{
		\@bburxtrue
		\@d@mscratch=#1
		\edef\@p@sbburx{\number\@d@mscratch}
}
\def\@p@@sbbury#1{
		\@bburytrue
		\@d@mscratch=#1
		\edef\@p@sbbury{\number\@d@mscratch}
}
\def\@p@@sheight#1{
		\@heighttrue
		\@d@mscratch=#1
   		\edef\@p@sheight{\number\@d@mscratch}
}
\def\@p@@swidth#1{
		\@widthtrue
		\@d@mscratch=#1
		\edef\@p@swidth{\number\@d@mscratch}
}
\def\@p@@srheight#1{
		\@rheighttrue
		\@d@mscratch=#1
		\edef\@p@srheight{\number\@d@mscratch}
}
\def\@p@@srwidth#1{
		\@rwidthtrue
		\@d@mscratch=#1
		\edef\@p@srwidth{\number\@d@mscratch}
}
\def\@p@@sright#1{\@righttrue \@surroundtrue}
\def\@p@@sleft#1{\@lefttrue \@surroundtrue}
\def\@p@@sextraheight#1{\@d@mextraheight=#1}
\def\@p@@sbox#1{\@boxtrue}
\def\@p@@scaption#1{\@captiontrue}
\def\@p@@stoplines#1{
		\@toplinestrue
		\@c@ttoplines=#1
}
\def\@p@@scaptionwidth#1{
		\@captionwidthtrue
	  	\@d@mcaptionwidth=#1
}
\def\@p@@scaptionwrite#1{
		\global\@captionwritetrue
		\global\@w@rname=\expandafter{\jobname_captions.tex}
		\typeout{Captions are written to \the\@w@rname}
}

\def\@p@@sprolog#1{\@prologfiletrue\def\@prologfileval{#1}}
\def\@p@@spostlog#1{\@postlogfiletrue\def\@postlogfileval{#1}}
\def\@cs@name#1{\csname #1\endcsname}
\def\@setparms#1=#2,{\@cs@name{@p@@s#1}{#2}}
\def\ps@init@parms{
		\@bbllxfalse \@bbllyfalse
		\@bburxfalse \@bburyfalse
		\@heightfalse \@widthfalse
		\@rheightfalse \@rwidthfalse
		\def\@p@sbbllx{}\def\@p@sbblly{}
		\def\@p@sbburx{}\def\@p@sbbury{}
		\def\@p@sheight{}\def\@p@swidth{}
		\def\@p@srheight{}\def\@p@srwidth{}
		\def\@p@sfile{}
		\def\@p@scost{10}
		\def\@sc{}
		\@prologfilefalse
		\@postlogfilefalse
		\@clipfalse
		\@rightfalse \@leftfalse
		\@boxfalse \@captionfalse
		\@toplinesfalse \@surroundfalse
		\@d@mextraheight=0pt
 		\@c@ttoplines=0
		\@pshape={} \def\@p@srheight@total{}
		\@captionwidthfalse \@d@mcaptionwidth=0pt
}
\def\parse@ps@parms#1{
	 	\@psdo\@psfiga:=#1\do
		   {\expandafter\@setparms\@psfiga,}}
\newif\ifno@bb
\newif\ifnot@eof
\newread\ps@stream
\newtoks\@linetok
\def\bb@missing{
	\typeout{psfig: searching \@p@sfile \space  for bounding box}
	\openin\ps@stream=\@p@sfile
	\no@bbtrue
	\not@eoftrue
	\catcode`\%=12
	\loop
		\read\ps@stream to \line@in
		\global\@linetok=\expandafter{\line@in}
		\ifeof\ps@stream \not@eoffalse \fi
		\@bbtest{\@linetok}
		\if@bbmatch\not@eoffalse\expandafter\bb@cull\the\@linetok\fi
	\ifnot@eof \repeat
	\catcode`\%=14
}	
\catcode`\%=12
\newif\if@bbmatch
\def\@bbtest#1{\expandafter\@a@\the#1
\long\def\@a@#1
     \ifx\@bbtest#2\@bbmatchfalse\else\@bbmatchtrue\fi}
\long\def\bb@cull#1 #2 #3 #4 #5 {
	\@d@mscratch=#2 bp\edef\@p@sbbllx{\number\@d@mscratch}
	\@d@mscratch=#3 bp\edef\@p@sbblly{\number\@d@mscratch}
	\@d@mscratch=#4 bp\edef\@p@sbburx{\number\@d@mscratch}
	\@d@mscratch=#5 bp\edef\@p@sbbury{\number\@d@mscratch}
	\no@bbfalse
}
\catcode`\%=14
\def\compute@bb{
		\no@bbfalse
		\if@bbllx \else \no@bbtrue \fi
		\if@bblly \else \no@bbtrue \fi
		\if@bburx \else \no@bbtrue \fi
		\if@bbury \else \no@bbtrue \fi
		\ifno@bb \bb@missing \fi
		\ifno@bb \typeout{FATAL ERROR: no bb supplied or found}
			\no-bb-error
		\fi
		\count203=\@p@sbburx
		\count204=\@p@sbbury
		\advance\count203 by -\@p@sbbllx
		\advance\count204 by -\@p@sbblly
		\edef\@bbw{\number\count203}
		\edef\@bbh{\number\count204}
}
\def\in@hundreds#1#2#3{\count240=#2 \count241=#3
		     \count100=\count240	
		     \divide\count100 by \count241
		     \count101=\count100
		     \multiply\count101 by \count241
		     \advance\count240 by -\count101
		     \multiply\count240 by 10
		     \count101=\count240	
		     \divide\count101 by \count241
		     \count102=\count101
		     \multiply\count102 by \count241
		     \advance\count240 by -\count102
		     \multiply\count240 by 10
		     \count102=\count240	
		     \divide\count102 by \count241
		     \count200=#1\count205=0
		     \count201=\count200
			\multiply\count201 by \count100
		     	\advance\count205 by \count201
		     \count201=\count200
			\divide\count201 by 10
		     	\multiply\count201 by \count101
			\advance\count205 by \count201
		     \count201=\count200
			\divide\count201 by 100
			\multiply\count201 by \count102
			\advance\count205 by \count201
		     \edef\@result{\number\count205}
}
\def\compute@wfromh{
		\in@hundreds{\@p@sheight}{\@bbw}{\@bbh}
		\edef\@p@swidth{\@result}
}
\def\compute@hfromw{
		\in@hundreds{\@p@swidth}{\@bbh}{\@bbw}
		\edef\@p@sheight{\@result}
}
\def\compute@handw{
		\if@height
			\if@width
			\else
				\compute@wfromh
			\fi
		\else
			\if@width
				\compute@hfromw
			\else
				\edef\@p@sheight{\@bbh}
				\edef\@p@swidth{\@bbw}
			\fi
		\fi
}
\def\compute@resv{
		\if@rheight \else \edef\@p@srheight{\@p@sheight} \fi
		\if@rwidth \else \edef\@p@srwidth{\@p@swidth} \fi
		\edef\@p@srheight@total{\@p@srheight}
}
\newtoks\@pshape
\def\@c@ttoplines{\count120}
\def\@c@theightcount{\count121}
\def\@c@tshapecount{\count122}
\newdimen\@d@mwidthshape
\newdimen\@d@mextraheight
\newdimen\@d@mscratch
\def\compute@parshape{
	\if@right
		\if@left
	   		\typeout{error: Can't have both left and right set}
			\@leftfalse
		\fi
	\fi
	\@d@mscratch=\@p@swidth truesp
	\divide \@d@mscratch by 19
	\multiply \@d@mscratch by 20
	\edef\@p@swidthdimen{\the\@d@mscratch}
	\@c@tshapecount=\@c@ttoplines
 	\@d@mscratch=\baselineskip
	\multiply \@d@mscratch by \@c@ttoplines
	\advance \@d@mscratch by .4\baselineskip
    	\edef\@p@stopdistance{\the\@d@mscratch }
	\@d@mscratch=\@p@sheight truesp
	\divide \@d@mscratch by 2
	\edef\@p@shalfboxheight{\the\@d@mscratch}
	\if@toplines
		\loop \@pshape=\expandafter{\the\@pshape 0pt \hsize}
		\advance\@c@ttoplines by -1
		\ifnum\@c@ttoplines>0 \repeat
	\fi
   	\@c@theightcount=\@p@srheight@total
	\advance \@c@theightcount by \@d@mextraheight
	\divide  \@c@theightcount by \baselineskip
	\advance \@c@theightcount by 1
    	\advance \@c@tshapecount by \@c@theightcount
	\advance \@c@theightcount by -1
	\@d@mwidthshape=\hsize
     	\advance \@d@mwidthshape by -\@p@swidthdimen
	\if@left
		\def\@moveshape{0pt}
		\ifnum\@c@theightcount>0
		      	\loop
			\@pshape=%
			\expandafter{\the\@pshape %
					\@p@swidthdimen \@d@mwidthshape}
			\advance \@c@theightcount by -1
			\ifnum\@c@theightcount>0 \repeat
		\else
			\advance \@c@tshapecount by 1
		\fi
	\fi
	\if@right
		\@d@mscratch=\hsize
		\advance \@d@mscratch by -\@p@swidth truesp
		\edef\@moveshape{\@d@mscratch}
		\ifnum\@c@theightcount>0
			\loop
			\@pshape=\expandafter{\the\@pshape 0pt \@d@mwidthshape}
			\advance \@c@theightcount by -1
			\ifnum\@c@theightcount>0 \repeat
		\else
			\advance \@c@tshapecount by 1
		\fi
	\fi
	\ifnum \@p@srheight=0
		\@pshape={}
		\@c@tshapecount = 0
	\else
	 	\@pshape=\expandafter{\the\@pshape 0pt \hsize}
	\fi
}
\def\@p@ssetsurroundboxes{
	\global\parshape=\count122 \the\@pshape		
 	\moveright\@moveshape
	\vbox to 0pt\bgroup\hskip0pt\vskip\@p@stopdistance
}
\newtoks\@captiontok
\newbox\@b@xcaption
\newdimen\@d@mcaptionwidth
\newdimen\@d@mcaptionheight
\newwrite\@w@rcaption
\newtoks\@w@rname
\def\setcaption#1{\@captiontok={#1}}
\def\@set@caption{
	\setbox\@b@xcaption\vbox{\hsize\@d@mcaptionwidth
	\tolerance=9000 \baselineskip=11.4pt
	\noindent\relax\the\@captiontok}
	\if@captionwrite
		\if@captionopen
		\else
			\global\@captionopentrue
			\immediate\openout\@w@rcaption=\the\@w@rname
		\fi
		\immediate\write\@w@rcaption{\the\@captiontok}
		\immediate\write\@w@rcaption{}
	\fi
}
\def\compute@caption{
	\if@captionwidth
	\else
		\@d@mcaptionwidth = \@p@swidth truesp
		\divide \@d@mcaptionwidth by 20
		\multiply \@d@mcaptionwidth by 17
	\fi
	\@set@caption
	\@d@mcaptionheight=\ht\@b@xcaption
	\if@rheight
	\else
		\count100=\@p@srheight
	   	\advance \count100 by \@d@mcaptionheight
	   	\advance \count100 by \bigskipamount
	   	\advance \count100 by \medskipamount
	   	\edef\@p@srheight@total{\number\count100}
	\fi
}
\newif\if@alreadyjtem \@alreadyjtemfalse
\def\newpar{\hfil\vadjust{\vskip\parskip}%
	{\count100=\parskip
	\count101=\baselineskip
	\divide\count101 by 10  \multiply\count101 by 3
	\advance \count100 by \count101
	\divide\count100 by \baselineskip
	\advance\count100 by \prevgraf
	\global\prevgraf=\count100}%
	\break\if@alreadyjtem\else\indent\fi%
}
\let\sav@par=\par
\def\jtem#1{%
    	\if@alreadyjtem\else\bgroup\fi
	\def\par{\sav@par\egroup\sav@par}
	\if@alreadyjtem\else\leftskip\parindent\fi
	\@alreadyjtemtrue
	\noindent\hskip0pt
	\llap{#1\ }\ignorespaces
}
\def\compute@sizes{%
	\compute@bb
	\compute@handw
  	\compute@resv
	\if@caption
		\compute@caption
	\fi
	\if@surround
		\compute@parshape
	\fi
}
\def\@p@sdospecials{
	\ifnum\@p@scost<\@psdraft
	       	\typeout{psfig: including \@p@sfile \space }
	\fi
	\special{ps::[begin] 	\@p@swidth \space \@p@sheight \space
			\@p@sbbllx \space \@p@sbblly \space
			\@p@sbburx \space \@p@sbbury \space
			startTexFig \space }
	\ifnum\@p@scost<\@psdraft
		\if@clip
			\typeout{(clip)}
			\special{ps:: \@p@sbbllx \space \@p@sbblly \space
				\@p@sbburx \space \@p@sbbury \space
			    	doclip \space }
		\fi
	\fi
	\if@box
		\typeout{(box)}
  		\special{ps:: \@p@sbbllx \space \@p@sbblly \space
			\@p@sbburx \space \@p@sbbury \space
		    	dobox \space }
	\fi
	\ifnum\@p@scost<\@psdraft
		\if@prologfile
	    		\special{ps: plotfile \@prologfileval \space }
		\fi
		\special{ps: plotfile \@p@sfile \space }
    		\if@postlogfile
			\special{ps: plotfile \@postlogfileval \space }
		\fi
	\fi
	\special{ps::[end] endTexFig \space }
}
\newif\if@putinvbox

\def\psfig#1{{%
	\ifhmode%
		\vbox\bgroup
		\@putinvboxtrue
	\else
		\@putinvboxfalse
	\fi
       	\ps@init@parms
	\parse@ps@parms{#1}
       	\compute@sizes
	\if@surround
		\psfig@for@surround
	\else
		\psfig@for@regular
	\fi
	\if@putinvbox
       		\egroup
	\fi
}}
\def\psfig@for@surround{%
	\@p@ssetsurroundboxes
	\ifnum\@p@scost<\@psdraft
		\@p@sdospecials
		\vbox to \@p@srheight true sp{\vss}
       	\else
		\if@box
			\@p@sdospecials
		\fi
		\vbox to \@p@srheight true sp{
			\vskip\@p@shalfboxheight
			\hbox to \@p@srwidth true sp{
				\hss
				\ifnum\@p@scost<\@psdraftspecial
					\@p@sfile
				\else
					\@p@sfake
				\fi
      				\hss
			}
		\vss
		}
	\fi
	\if@caption
		\medskip
		\hbox to \@p@srwidth true sp{
			\hss
			\box\@b@xcaption
			\hss
		}
 		\medskip
	\fi
	\vss\egroup
	\vskip-\parskip
}

\def\psfig@for@regular{%
	\if@putinvbox
	\else
		\vskip\parskip
	\fi
	\ifnum\@p@scost<\@psdraft
		\@p@sdospecials
		\vbox to \@p@srheight true sp{%
			\hbox to \@p@srwidth true sp{
			\hfil
			}
		\vfil
		}
       	\else
		\if@box
			\@p@sdospecials
		\fi
	    	\vbox to \@p@srheight true sp{
			\vss
			\hbox to \@p@srwidth true sp{
				\hss
				\ifnum\@p@scost<\@psdraftspecial
					\@p@sfile
				\else
					\@p@sfake
				\fi
				\hss
			}
		    	\vss
		}
	\fi
	\if@caption
		\medskip
		\hbox to \@p@srwidth true sp{
			\hss
			\box\@b@xcaption
			\hss
		}
		\bigskip
	\fi
	\if@putinvbox
	\else
		\vskip-\parskip
	\fi
}
\catcode`\@=12\relax

\psfiginit

\font\tinybbfont=msbm6
\font\scriptsizebbfont=msbm7 scaled \magstep 1
 1
\font\footnotesizebbfont=msbm9 scaled \magstep 0
 2
\font\bbfont=msbm9 scaled \magstep1  
 1
 2
 3
 4

\def\tinyBbb#1{\hbox{\tinybbfont #1}}
\def\scriptsizeBbb#1{\hbox{\scriptsizebbfont #1}}

\def\footnotesizeBbb#1{\hbox{\footnotesizebbfont #1}}

\def\Bbb#1{\hbox{\bbfont #1}}

\newcommand{\Br}{\mbox{\it Br}\,}

\newcommand{\CDiv}{\mbox{\it CDiv}\,}
\newcommand{\CP}{{\Bbb C}{\rm P}}

\newcommand{\ELB}{\mbox{\it ELB}\,}
\newcommand{\Id}{\mbox{\it Id}}

\newcommand{\Hom}{\mbox{\it Hom}\,}
\newcommand{\Ind}{\mbox{\it Ind}\,}
\newcommand{\Int}{\mbox{\it Int}\,}

\newcommand{\Pic}{\mbox{\rm Pic}\,}

\newcommand{\PDiv}{\mbox{\it PDiv}\,}

\newcommand{\SF}{\mbox{\it SF}\,}
\newcommand{\Sing}{\mbox{\it Sing}\,}
\newcommand{\SL}{\mbox{\it SL}}

\newcommand{\Span}{\mbox{\it Span}\,}
\newcommand{\Spec}{\mbox{\it Spec}\,}

\newcommand{\Star}{\mbox{\it Star}}
\newcommand{\WCP}{\mbox{\rm W{\Bbb C}P}}
\newcommand{\footnotesizeWCP}{\mbox{\rm W{\footnotesizeBbb C}P}}
\newcommand{\scriptsizeWCP}{\mbox{\rm W{\scriptsizeBbb C}P}}

\newcommand{\codim}{\mbox{\it codim}\,}
\newcommand{\coker}{\mbox{\it coker}\:}

\newcommand{\dimm}{\mbox{\it dim}\,}

\newcommand{\im}{\mbox{\rm Im}\,}

\newcommand{\kker}{\mbox{\it ker}\:}

\newcommand{\mult}{\mbox{\it mult}\,}
\newcommand{\pr}{\mbox{\it pr}}

\newcommand{\rank}{\mbox{\it rank}\,}

\newcommand{\res}{\mbox{\it res}\,}

\newcommand{\simeqrightarrow}{\stackrel{\sim}{\rightarrow}}

\begin{document}


 \title{Toric morphisms and
        fibrations of toric Calabi-Yau hypersurfaces}
 \arxurl{math.AG/0010082}     
 \author{Yi Hu$^1$, Chien-Hao Liu$^2$, and Shing-Tung Yau$^3$} 
\makeatletter
 \address{\vspace{-.2in}$^1$Department of Mathematics, 
              University of Arizona,
              Tucson, Arizona 85721 }
\addressemail{yhu@math.arizona.edu}
\address{$^2$Department of Mathematics,
             Harvard University,
             Cambridge, MA 02138}   
\addressemail{chienliu@math.harvard.edu}
\address{$^3$Department of Mathematics,
             Harvard University,
             Cambridge, MA 02138}   
\addressemail{yau@math.harvard.edu}
\makeatletter
 \markboth{\it TORIC MORPHISMS AND FIBRATIONS OF CALABI-YAU \ldots}{\it
               Y.\ HU, C.-H.\ LIU, and S.-T.\ YAU} 
 \begin{abstract}
  Special fibrations of toric varieties have been used by physicists,
   e.g.\ the school of Candelas, to construct dual pairs in the study
   of Het/F-theory duality. Motivated by this, we investigate in this
   paper the details of toric morphisms between toric varieties.
   In particular, a complete toric description of fibers - both generic
   and non-generic -, image, and the flattening stratification of
   a toric morphism are given.
   Two examples are provided to illustrate the discussions.
  We then turn to the study of the restriction of a toric morphism
   to a toric hypersurface. The details of this can be understood by
   the various restrictions of a line bundle with a section 
   that defines the hypersurface.
  These general toric geometry discussions give rise to a computational
   scheme for the details of a toric morphism and the induced fibration
   of toric hypersurfaces therein. We apply this scheme to study
   the family of complex $4$-dimensional elliptic Calabi-Yau
   toric hypersurfaces that appear in a recent work of
   Braun-Candelas-dlOssa-Grassi.
  Some directions for future work are listed in the end.
 \end{abstract}



\bigskip

\noindent {\small
MSC number 2000$\,$: 14M25, 81T30, 14J32, 14D06, 14J35.
} 

\newpage
{\small
\noindent{\bf Acknowledgements.}
This work is a further study of some original ideas from
Philip Candelas and his school. We would like to thank him for
educations and the source of inspirations.
We would like to thank also
 Alex Avram,
 Volker Braun, 
 Hung-Wen Chang,
 Jiun-Cheng Chen, 
 Ti-Ming Chiang,
 Xenia de la Ossa,
 Jacques Distler,
 Daniel Freed, 
 Robert Gompf, 
 Antonella Grassi,
 Brian Greene, 
 Joe Harris,
 Mark Haskins,
 Shinobu Hosono,
 Sheldon Katz, 
 Sean Keel,
 Albrecht Klemm,
 Bong H.\ Lian, 
 Chiu-Chu Liu,
 David Morrison, 
 Eugene Perevalov,
 Govindan Rajesh,
 Lorenzo Sadun, 
 David Saltman,
 Herald Skarke,
 Jason Starr,
 Margaret Symington,
 Richard Thomas,
 Karen Uhlenbeck, 
 Cumrun Vafa, 
 Xiao-Wei Wang,
 and Eric Zaslow
for valuable inspirations, educations, conversations, and
correspondences at various stages of the work. 
C.H.L.\ would like to thank in addition
 Orlando Alvarez and William Thurston for indispensable educations,
 P.C.\ for the two years' group meetings in his school,
 the Geometry and String Theory Group at University of Texas at Austin
 for crucial influences at the brewing stage of the work,
 Nathan Dunfield for making the Adobe figures latex compatible,
 and
 Ling-Miao Chou for tremendous moral support.
The work is supported by
 DOE grant DE-FG02-88ER25065 and NSF grant DMS-9803347.
} 

$ $

\vspace{-4em}  

%
%

i
\setcounter{section}{-1}
\section{Introduction and outline.}

\begin{flushleft}
{\bf Introduction.}
\end{flushleft}
Special fibrations of toric varieties have been used by physicists,
 e.g.\ the school of Candelas
 (cf.\ [A-K-M-S], [B-C-dlO-G], [B-M], [C-F], [C-P-R], [K-S1], [K-S2],
       [P-S], and [Ra]),
 to construct dual pairs in the study of Heterotic-string/F-theory
 duality. This motivates us to investigate the details of toric
 morphisms between toric varieties and the induced morphisms from
 the toric hypersurfaces in the domain toric variety of a toric
 morphism to the target toric variety.
The goal of this paper is to describe these details in terms of
 toric geometry and then to provide a computational scheme for
 actual toric computations.

The toric description of the generic fiber of a toric morphism has
 been discussed extensively e.g.\ in the works of the school of
 Candelas.
Thus, the first theme of this article is to extend their results
 to a more general setting that takes care of all the fibers
 - generic and non-generic alike - in terms of toric geometry.
 The concept of primitive cones and relative stars with respect
 to a toric morphism are introduced along the way.

Once the details of a toric morphism are described in terms of toric
 geometry, our second theme is to understand the details of the
 induced morphism from a toric hypersurface in the domain toric
 variety to the target toric variety. This is an important step in
 understanding fibred toric Calabi-Yau hypersurfaces.
 It turns out that, by turning the problem to a study of the various
 restrictions of the line bundle with a section that defines the
 hypersurface, one can extract certain details of this induced
 morphism in terms of toric geometry as well.

As long as applications to Heterotic-string/F-theory duality is
 concerned, the discussions in this article will be of no use if
 no methods of concrete toric computations for the various details
 outlined above can be provided. Thus our third theme in this article
 is to give a computational scheme to match with the toric geometry
 discussion. This is demonstrated by working out more details of
 an example$\,$: the toric morphism from a toric $5$-variety to
 a toric $3$-variety and the induced elliptic fibration of
 a $4$-dimensional toric Calabi-Yau hypersurfaces from a recent
 collaborative work of Braun, Candelas, de la Ossa, and Grassi
 ([B-C-dlO-G]).
 We expect that the computational scheme presented here together
 with other works from the school of Candelas, notably the
 classification by Kreuzer and Skarke of higher dimensional reflexive
 polytopes that admit fibrations (e.g.\ [K-S1] and [K-S2]), will
 provide us with some more concrete understanding of fibred Calabi-Yau
 manifolds beyond current theory.
 This will await another work in future.

This article is organized as follows. 
After giving the references of the necessary backgrounds for
 physicists in Sec.\ 1, we study in Sec.\ 2 the various details of
 a toric morphism.
 In particular, a complete toric description of fibers
  - both generic and non-generic -, image, and the flattening
  stratification of a toric morphism are worked out in Sec.\ 2.1.
 Two examples are provided in Sec.\ 2.2 to illustrate the discussions.
We then turn to the study of the restriction of a toric morphism
 to a toric hypersurface. After providing some necessary facts 
 about line bundles on a toric variety in Sec.\ 3.1, we relate in
 Sec.\ 3.2 the problem of induced morphism to the equivalent problem
 of the various restrictions of a line bundle and the section in it
 that defines the hypersurface and study the latter.
These general toric geometry discussions give rise to a computational
 scheme for the details of a toric morphism and the induced morphism
 of toric hypersurfaces therein. After recalling some necessary facts
 from the work of Batyrev in Sec.\ 4.1, we lay down this scheme in
 Sec.\ 4.2 and apply it to study the family of
 $4$-dimensional elliptic Calabi-Yau toric hypersurfaces that was
 discussed in a recent work of Braun-Candelas-dlOssa-Grassi.
The Maple codes that are employed for the computation can be found
 in the preprint version of this article$\,$: {\tt math.AG/0010082}.
Some directions for future work are listed in Sec.\ 5.

\bigskip

\noindent
{\bf Convention 0.1 [real vs.\ complex manifolds].} {
A {\it real} $n$-dimensional manifold will be called an
{\it $n$-manifold} while a {\it complex} $n$-dimensional
manifold an {\it $n$-fold}. Also, a real $n$-dimensional orbifold
will be called an {\it $n$-orbifold} and a complex $n$-dimensioanl
variety will be called an {\it $n$-variety}.
All the varieties in this article are over ${\Bbb C}$ and
by {\it a point in a variety} we mean a {\it closed point} of
that variety only.
} 

\bigskip

\noindent
{\bf Convention 0.2 [morphism vs.\ fibrations].} {
 In this article, a {\it morphism} $f:X\rightarrow Y$ means a general
 map between varieties $X$ and $Y$ as in algebraic geometry;
 a {\it fiber} of $f$ means the point-set pre-image $f^{-1}(y)$ of $f$
 at a certain (closed) point $y\in Y$. 
 A morphism $f:X\rightarrow Y$ is called a {\it fibration} if $f$ is
 surjective and all its fibers have the same dimension
 $\dimm X-\dimm Y$.
} 

\bigskip

\begin{flushleft}
{\bf Outline.}
\end{flushleft}
{\small
\baselineskip 11pt  

\begin{quote}
 1. Essential mathematical background for physicists.

 2. The toric geometry of toric morphisms.
    \vspace{-1ex}
    \begin{quote}
     \hspace{-1.3em}
     2.1 The toric geometry of toric morphisms.

     \hspace{-1.3em}
     2.2 Examples of toric morphisms and their fibers.
    \end{quote}

 \vspace{-.8ex}
 3. Induced morphism and fibers for hypersurfaces.
    \vspace{-1ex}
    \begin{quote}
     \hspace{-1.3em}
     3.1 Preparations.

     \hspace{-1.3em}
     3.2 The induced morphism$\,$:
         $Y\subset X_{\Sigma^{\prime}}\rightarrow X_{\Sigma}$.
    \end{quote}

 \vspace{-.8ex}
 4. Fibration of Calabi-Yau hypersurfaces via toric morphisms.
    \vspace{-1ex}
    \begin{quote}
     \hspace{-1.3em}
     4.1 General remarks.

     \hspace{-1.3em}
     4.2 The computational scheme and a detailed study of an example.
    \end{quote}

 \vspace{-.8ex}
 5. Remarks on further study.
\end{quote}
} 

\bigskip

\section{Essential mathematical background for physicists.}

Some references of the necessary mathematical backgrounds are
 provided here for the convenience of physicists.
The notations of the basic ingredients of toric geometry are
 also introduced in this section.
Essential necessary facts will be recalled along the way of
 discussions.

\bigskip

\noindent $\bullet$
{\bf Toric geometry.}
Excellent expositories are given in [C-K], [Fu], and [Gr].
 See also [Da], [Do], [Ew], [G-K-Z3], [Od2], and [Zi].
In this article, we fix the following notations$\,$:

\bigskip

{\it Notation}$\,$:
\begin{quote}
  $N\cong{\Bbb Z}^n$: a lattice;

  $M=\Hom(N,{\Bbb Z})\,$: the dual lattice of $N$;

  $T_N=\Hom(M,{\Bbb C}^{\ast})\,$:
      the (complex) torus associated to $N$;

  $\Sigma\,$: a fan in $N_{\scriptsizeBbb R}$;

  $\tau\prec\sigma\,$: the relation between cones $\tau$ and $\sigma$
      that $\tau$ is in the face of $\sigma$;     

  $X_{\Sigma}\,$: the toric variety associated to $\Sigma$;

  $\Sigma(i)\,$: the $i$-skeleton of $\Sigma$;

  $U_{\sigma}\,$: the local affine chart of $X_{\Sigma}$ associated
     to $\sigma$ in $\Sigma$;

  $x_{\sigma}\in U_{\sigma}\,$: the distinguished point associated
     to $\sigma$;

  $O_{\sigma}\,$: the $T_N$-orbit of $x_{\sigma}$ under the
     $T_N$-action on $X_{\Sigma}$;

  $V(\sigma)\,$: the orbit closure of $O_{\sigma}$;

  $N_{\sigma}\,$: the sublattice of $N$ generated as a subgroup by
                  $\sigma\cap N$;

  $M(\sigma)=\sigma^{\perp}\cap M\,$;

  $\Span_{\scriptsizeBbb R}(\sigma)$:
    the real vector subspace spanned by $\sigma$;

  {\it interior} of $\sigma\,$: the topological interior of $\sigma$
      in $\Span_{\scriptsizeBbb R}(\sigma)$;

  $\lambda_v\,$: the one-parameter subgroup of $T_N$ associated to
      $v\in N$;

  $\chi^m\,$: the character of $T_N$ associated to $m\in M$.
\end{quote}

\bigskip

\noindent $\bullet$
{\bf Cox homogeneous coordinate ring.}
See [C-K] and [Cox].

\bigskip

\noindent $\bullet$
{\bf Toric Calabi-Yau hypersurfaces.}
See [Ba2].
Also [Bo], [B-B], [C-K], and [H\"{u}].

\bigskip

\noindent $\bullet$
{\bf Fibred Calabi-Yau spaces.}
See [Og1], [Og2], and [Og3].
Also [A-M], [H-L-Y], [K-L-M], [K-L-R-Y] and the various works
 from the school of Candelas, e.g.\
 [A-K-M-S], [B-C-dlO-G], [B-M], [C-F], [C-P-R], [K-S1], [K-S2],
 [P-S], and [Ra].

\bigskip

\section{The toric geometry of toric morphisms.}

The toric description of details of a toric morphism is given in
 Sec.\ 2.1.
Two examples are provided in Sec.\ 2.2 to illustrate the ideas.

\bigskip

\subsection{The toric geometry of toric morphisms.}

Recall first the definition:

\bigskip

\noindent
{\bf Definition 2.1.1 [toric morphism].} ([Ew] and [Fu].) {
 Let $\Sigma^{\prime}$, $\Sigma$ be fans in $N^{\prime}$, $N$
 respectively. A {\it map between fans}, in notation
 $\varphi:\Sigma^{\prime}\rightarrow \Sigma$, is a homomorphism
 $\varphi:N^{\prime}\rightarrow N$ of lattices that satisfies the
 condition: For each $\sigma^{\prime}\in\Sigma^{\prime}$, there
 exists a $\sigma\in\Sigma$ such that
 $\varphi(\sigma^{\prime})\subset\sigma$.
 Such $\varphi$ determines a morphism
 $\widetilde{\varphi}:X_{\Sigma^{\prime}}\rightarrow X_{\Sigma}$.
 A morphism between toric varieties that arises in this way is called
 a {\it toric morphism}. (Cf.\ {\sc Figure 2-1-1}.)
} 
 \begin{figure}[htbp]
  \setcaption{{\sc Figure 2-1-1.}
   \baselineskip 14pt
   Both $\pi_1$ and $\pi_2$ are the natural projection map.
   In (a) $\pi_1$ induces a map between fans 
   while in (b) $\pi_2$ does not. The cone of the fan in (b) that
   is not compatible with $\pi_2$ is indicated by a darker shading.
  } 
 \end{figure}

\bigskip

The following two facts, adapted from [Ew] and [Ei], are fundamental
to understanding toric morphisms.

\bigskip

\noindent
{\bf Fact 2.1.2 [equivariance].}
(Cf.\ Theorem 6.4 in [Ew].) {\it
 The toric morphism
 $\widetilde{\varphi}:X_{\Sigma^{\prime}}\rightarrow X_{\Sigma}$
 is equivariant with respect to the homomorphism
 $T_{N^{\prime}}\rightarrow T_N$ induced by
 $\varphi:N^{\prime}\rightarrow N$.
} 

\bigskip

\noindent
{\bf Fact 2.1.3 [semicontinuity of fiber dimension].}
(Cf.\ Theorem 14.8 in [Ei].) {\it
 Let $\varphi:X\rightarrow Y$ be a morphism between projective
 varieties and $s=\max\{\,\dimm X-\dimm Y,\, 0\,\}$.
 Then there is a finite filtration of $Y$ by subvarieties
 (which need not be all distinct),
 $$
  Y\;\supset\; \varphi(X)\,=\, Y_s\;
   \supset Y_{s+1}\; \supset\; \;\cdots\; \supset Y_{\it dim\, X}\,,
 $$
 such that
 $\dimm\varphi^{-1}(y)\ge i$ if and only if $y\in Y_i$.
} 

\bigskip

The main goal of this subsection is to prove Proposition 2.1.4 below,
which extends both a result in [A-K-M-S] related to reflexive
polytopes and the usual description of toric bundles as in
e.g.\ [Ew] and [Fu] to full generality.
It contains the toric analogue of Fact 2.1.3.

\bigskip

\begin{flushleft}
{\bf Toric morphism$\,$: Image, fibers, and a flattening stratification.}
\end{flushleft}

\noindent
{\bf Proposition 2.1.4 [toric morphism].} {\it
 Let $\widetilde{\varphi}:X_{\Sigma^{\prime}}\rightarrow X_{\Sigma}$
 be a toric morphism induced by a map of fans
 $\varphi:\Sigma^{\prime}\rightarrow \Sigma$. Then$\,$:
 \begin{itemize}
  \item
  The image $\widetilde{\varphi}(X_{\Sigma^{\prime}})$ of
  $\widetilde{\varphi}$ is a subvariety of $X_{\Sigma}$.
  It is realized as the toric variety associated to the fan
  $\Sigma_{\varphi}
     \doteq \Sigma\cap\varphi(N^{\prime}_{\scriptsizeBbb R})$.

  \item
  The fiber of $\widetilde{\varphi}$ over a point
  $y\in X_{\Sigma_{\varphi}}$ depends only on the orbit $O_{\sigma}$,
  $\sigma\in\Sigma_{\varphi}$, that contains $y$. Denote this fiber
  by $F_{\sigma}$, then it can be described as follows.
  \begin{quote}
   Define $\Sigma^{\prime}_{\sigma}$ to be the set of cones
   $\sigma^{\prime}$ in $\Sigma^{\prime}$, whose interior is mapped to
   the interior of $\sigma$. Let $\Ind(\sigma)$ be the index of
   $\widetilde{\varphi}$ over $O_{\sigma}$
   (cf.\ Definition/Lemma 2.1.7).
   Then $\varphi^{-1}(y)=F_{\sigma}$ is a disjoint union of
   $\Ind(\sigma)$ identical copies of connected
   {\sl reducible toric variety} $F_{\sigma}^c$, whose irreducible
   components $F^{\tau^{\prime}}_{\sigma}$ are the toric variety
   associated to the relative star $\Star_{\sigma}(\tau^{\prime})$ of
   the primitive elements $\tau^{\prime}$ in $\Sigma^{\prime}_{\sigma}$
   (cf.\ Definition 2.1.9, Definition 2.1.10, and Lemma 2.1.11).
  \end{quote}
  {\rm ({\sc Figure 2-1-2}.)}

  \item
  For $\sigma\in\Sigma_{\varphi}$,
  $\widetilde{\varphi}^{-1}(O_{\sigma})
                  =\widetilde{O_{\sigma}}\times F_{\sigma}^c$,
  where $\widetilde{O_{\sigma}}$ is a connected covering space
  of $O_{\sigma}$ of order $\Ind(\sigma)$.
  Thus, the stratification of $X_{\Sigma}$ by
  $X_{\Sigma}-\widetilde{\varphi}(X_{\Sigma^{\prime}})$ and the
  toric orbits $O_{\sigma}$, $\sigma\in\Sigma_{\varphi}$, of
  $X_{\Sigma_{\varphi}}$ in $\widetilde{\varphi}(X_{\Sigma^{\prime}})$
  gives a flattening stratification of $\widetilde{\varphi}$.
 \end{itemize}
}
%
 \begin{figure}[htbp]
  \setcaption{{\sc Figure 2-1-2.}
   \baselineskip 14pt
    How the various fibers and their index  of a toric morphism may
    look like is indicated. Though they could be wild in general,
    all the information is coded in the map between fans.
  } 
  \centerline{\psfig{figure=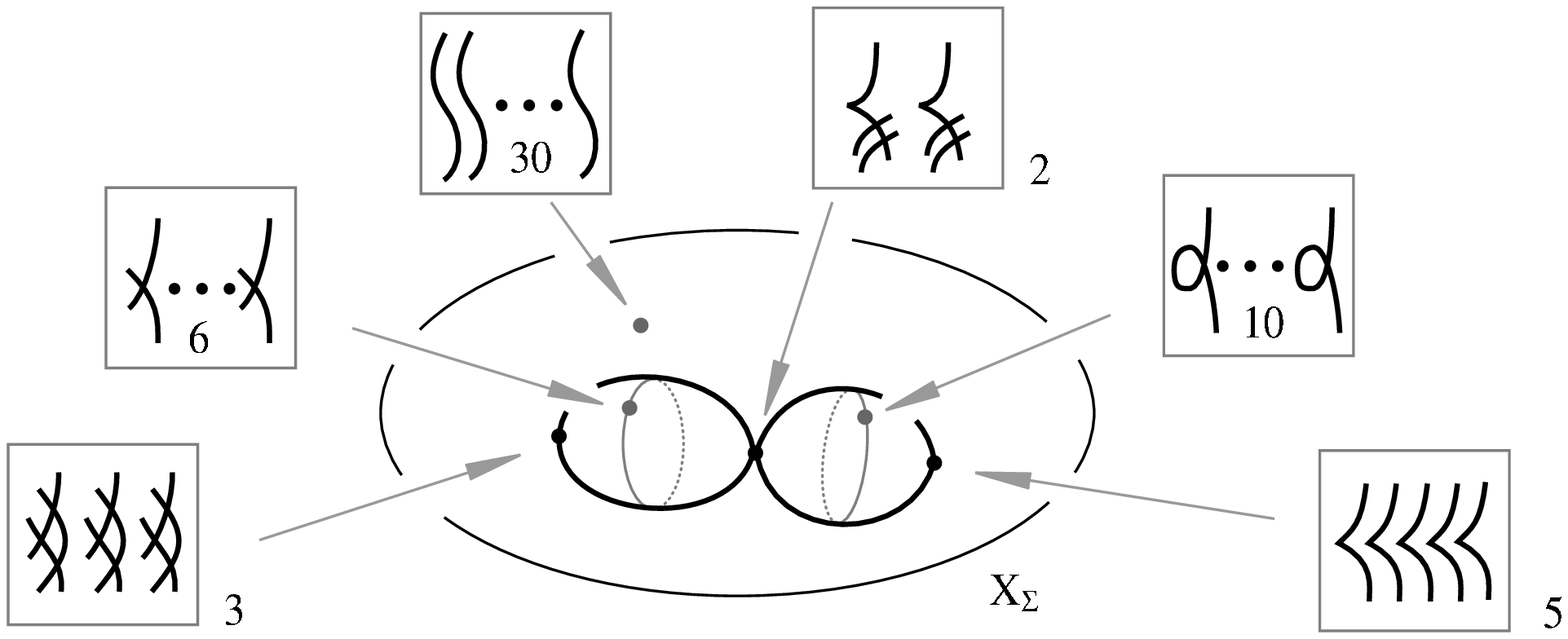,width=11cm,caption=}}
 \end{figure}

\bigskip

\noindent
{\it Remark 2.1.5.} Here the term {\it reducible toric variety} means
a reducible variety obtained by gluing a collection of toric varieties
along some isomorphic toric orbits.

\bigskip

\noindent
{\it Proof of Proposition 2.1.4.}
Let $N_{\varphi}=N\cap\varphi(N^{\prime}_{\scriptsizeBbb R})$.
Then the toric morphism
$\widetilde{\varphi}:X_{\Sigma^{\prime}}\rightarrow X_{\Sigma}$
is the composition of toric morphisms
$X_{\Sigma^{\prime}}
  \stackrel{\widetilde{\varphi}}{\rightarrow}
 X_{\Sigma_{\varphi}}\stackrel{\widetilde{\iota}}{\rightarrow} X_{\Sigma}$
induced from the maps of fans
$\Sigma^{\prime}\stackrel{\varphi}{\rightarrow}
   \Sigma_{\varphi}\stackrel{\iota}{\rightarrow}\Sigma$
from linear maps
$N^{\prime}\stackrel{\varphi}{\rightarrow} N_{\varphi}
                                    \stackrel{\iota}{\rightarrow} N$.
Note that $\varphi(N^{\prime})$ is a subgroup of $N_{\varphi}$
of finite index; hence,
$\widetilde{\varphi}:X_{\Sigma^{\prime}}\rightarrow X_{\Sigma_{\varphi}}$
is surjective.
On the other hand, since
$\Sigma_{\varphi}=\Sigma\cap N_{\varphi\,\scriptsizeBbb R}$,
$\widetilde{\iota}(X_{\Sigma_{\varphi}})$ is the closure of the
$T_{N_{\varphi}}$-orbit $T_{N_{\varphi}}\cdot x_{\{0\}}$ in $X_{\Sigma}$.
Since it contains $T_{N_{\varphi}}$ as an $T_{N_{\varphi}}$-equivariant
open dense subset, from the relation between torus actions and fans
(e.g.\ [Fu]) one concludes that the closure is simply
$X_{\Sigma_{\varphi}}$. This shows that
$\widetilde{\iota}:X_{\Sigma_{\varphi}}\rightarrow X_{\Sigma}$
is an inclusion.
Together, one concludes that $\widetilde{\varphi}(X_{\Sigma^{\prime}})$
is a subvariety of $X_{\Sigma}$ and it is realized as
$X_{\Sigma_{\varphi}}$.

\bigskip

{\it For the remaining part of the proof}, we shall regard
$\widetilde{\varphi}$ as a toric morphism
$X_{\Sigma^{\prime}}\rightarrow X_{\Sigma_{\varphi}}$ induced by
the map of fans $\varphi:\Sigma^{\prime}\rightarrow\Sigma_{\varphi}$.
After relabelling to simplify the notation, without loss of generality
{\it we shall assume in the following discussion that
$\varphi:N^{\prime}_{\scriptsizeBbb R}\rightarrow N_{\scriptsizeBbb R}$
is surjective and, hence $\Sigma_{\varphi}=\Sigma$}.
Note that in this case $\varphi(N^{\prime})$ is a subgroup of finite
index in $N$ and the induced group homomorphism
$T_{N^{\prime}}\rightarrow T_N$ is surjective.

\bigskip

Since $\widetilde{\varphi}$ is toric equivariant, it sends a toric
orbit $O_{\sigma^{\prime}}$ in $X_{\Sigma^{\prime}}$ into a toric
orbit $O_{\sigma}$ in $X_{\Sigma}$. The following lemma characterizes
such pairs $(\sigma^{\prime},\,\sigma)$.

\bigskip

\noindent
{\bf Lemma 2.1.6 [orbit pair].} {\it
 $\widetilde{\varphi}$ maps $O_{\sigma^{\prime}}$ to $O_{\sigma}$
 if and only if $\varphi$ maps the interior of $\sigma^{\prime}$
 to the interior of $\sigma$. When this happens, 
 $\widetilde{\varphi}$ maps $O_{\sigma^{\prime}}$ surjectively
 to $O_{\sigma}$.
} 

\bigskip

\noindent
{\it Proof.}
The if-part follows from the natural presentaion of orbits as affine
subvarieties in a toric variety ([Fu]):
$O_{\sigma^{\prime}}
  =\Spec({\Bbb C}[{\sigma^{\prime}}^{\perp}\cap M^{\prime}])$
and
$O_{\sigma}=\Spec({\Bbb C}[\sigma^{\perp}\cap M])$.
Note that, for each $\sigma^{\prime}\in\Sigma^{\prime}$, there is a
unique $\sigma\in\Sigma$ such that the interior of $\sigma^{\prime}$
is mapped to the interior of $\sigma$. Thus, $\widetilde{\varphi}$
maps $O_{\sigma^{\prime}}$ to a particular $O_{\sigma}$ determined
by $\varphi$. Since $X_{\Sigma}$ is stratified by toric orbits and
the image of an orbit in $X_{\Sigma^{\prime}}$ can be contained only
in a unique orbit in $X_{\Sigma}$, the only-if-part also follows.
For such $(\sigma^{\prime},\,\sigma)$, the surjectivity of
$\widetilde{\varphi}:O_{\sigma^{\prime}}\rightarrow O_{\sigma}$
follows from the surjectivity of $T_{N^{\prime}}\rightarrow T_N$.
This concludes the proof.

\noindent\hspace{12cm} $\Box$

\bigskip

To understand the fiber over each orbit, consider a homomorphism
between lattices $\psi:N_1\rightarrow N_2$ with $\psi(N_1)$
of finite index in $N_2$. Let $[N_2:\psi(N_1)]$ be the index of
$\psi(N_1)$ in $N_2$ and let
$\widetilde{\psi}:T_{N_1}\rightarrow T_{N_2}$ be the induces
morphism, where $T_{N_1}$, $T_{N_2}$ are regarded as the toric
variety associated to the fans that consist only of the origin.
One has then the following exact sequences
$$
 \begin{array}{cccccccccc}
  0 & \longrightarrow & \psi^{-1}(0)  & \longrightarrow   & N_1
    & \longrightarrow & \psi(N_1)     & \longrightarrow   & 0
                                                      &, \\[.6ex]
  0 & \longrightarrow & \psi(N_1)     & \longrightarrow   & N_2
    & \longrightarrow & N_2/\mbox{\raisebox{-.4ex}{$\psi(N_1)$}}
    & \longrightarrow & 0  &. 
 \end{array}
$$
Note that $\psi^{-1}(0)$ is a lattice of rank $\rank N_1-\rank N_2$
and $N_2/\mbox{\raisebox{-.4ex}{$\psi(N_1)$}}$ is a finite abelian
group. From the first exact sequence, the homomorphism
$T_{N_1}\rightarrow T_{\psi(N_1)}$ is a bundle map with fiber
$T_{\psi^{-1}(0)}$. From the second exact sequence, the homomorphism
$T_{\psi(N_1)}\rightarrow T_{N_2}$ is a finite covering map
with fiber a finite set of $[N_2:\psi(N_1)]$-many elements.
Since the composition of these two morphisms gives the morphism
$\widetilde{\psi}:T_{N_1}\rightarrow T_{N_2}$, one concludes
that the fiber of $\widetilde{\psi}$ is a disjoint union of
$[N_2:\psi(N_1)]$ copies of $T_{\psi^{-1}(0)}$.
Notice that the same conclusion can also be phrased in terms of
lattices $M_1$, $M_2$ dual to $N_1$, $N_2$ respectively.

Apply this to our problem. Let $\Sigma^{\prime}_{\sigma}$ be the
set of cones in $\Sigma^{\prime}$ whose interior is mapped to the
interior of $\sigma\in\Sigma$.
Let $\sigma^{\prime}\in\Sigma^{\prime}_{\sigma}$. Since
$\widetilde{\varphi}:O_{\sigma^{\prime}}\rightarrow O_{\sigma}$
is surjective, the image of the homomorphism behind,
$\overline{\varphi}^{\sigma^{\prime}}_{\sigma}:
  N^{\prime}/\mbox{\raisebox{-.4ex}{$N^{\prime}_{\sigma^{\prime}}$}}
                \rightarrow N/\mbox{\raisebox{-.4ex}{$N_{\sigma}$}}$,
must be of finite index in $N/\mbox{\raisebox{-.4ex}{$N_{\sigma}$}}$.
(Equivalently, the induced homomorphism
 $\varphi^{\dagger}_{\sigma^{\prime}}:\sigma^{\perp}\cap M
        \rightarrow {\sigma^{\prime}}^{\perp}\cap M^{\prime}$
 must be injective.)
Since
$$
 {\overline{\varphi}^{\sigma^{\prime}}_{\sigma}}^{-1}(0)\;
  =\; \varphi^{-1}(N_{\sigma})/                 
         \mbox{\raisebox{-.4ex}{$N^{\prime}_{\sigma^{\prime}}$}}\,,
$$
the fiber $F_{\sigma}^{\sigma^{\prime}}$ of $\widetilde{\varphi}$
over a point in $O_{\sigma}$ is thus a disjoint union of
$[\,N/\mbox{\raisebox{-.4ex}{$N_{\sigma}$}}\,:\,
   \overline{\varphi}^{\sigma^{\prime}}_{\sigma}(N^{\prime}/
    \mbox{\raisebox{-.4ex}{$N^{\prime}_{\sigma^{\prime}}$}})\,]$
copies of
$T_{\varphi^{-1}(N_{\sigma})/
  \mbox{\scriptsize\raisebox{-.4ex}{$N^{\prime}_{\sigma^{\prime}}$}}}$.

The following two lemmas play an essential role in understanding
fibers of $\widetilde{\varphi}$ and their relations.

\bigskip

\noindent
{\bf Definition/Lemma 2.1.7
[index of $\widetilde{\varphi}$ over $O_{\sigma}$].}
{\it
 The image
 $\overline{\varphi}^{\sigma^{\prime}}_{\sigma}(
   N^{\prime}/
      \mbox{\raisebox{-.4ex}{$N^{\prime}_{\sigma^{\prime}}$}}
 )$
 in $N/\mbox{\raisebox{-.4ex}{$N_{\sigma}$}}$ is independent of
 the choice of $\sigma^{\prime}$ in $\Sigma^{\prime}_{\sigma}$.
 Thus, the index
 $[\,N/\mbox{\raisebox{-.4ex}{$N_{\sigma}$}}\,:\,
    \overline{\varphi}^{\sigma^{\prime}}_{\sigma}(N^{\prime}/
     \mbox{\raisebox{-.4ex}{$N^{\prime}_{\sigma^{\prime}}$}})\,]$
 depends only on $\sigma$ and will be denoted by $\Ind(\sigma)$
 and called the
 {\sl index of $\widetilde{\varphi}$ over $O_{\sigma}$}.
} 

\bigskip

\noindent
{\it Proof.} Let
$\sigma^{\prime}_1,\,\sigma^{\prime}_2\in\Sigma^{\prime}_{\sigma}$.
Assume first that $\sigma^{\prime}_1$ is in the face of
$\sigma^{\prime}_2$. Then $N^{\prime}_{\sigma^{\prime}_1}$ is a
subspace of $N^{\prime}_{\sigma^{\prime}_2}$ and one has the exact
sequence$\,$:
$$
  0\;\longrightarrow\; 
    N^{\prime}_{\sigma^{\prime}_2}/
       \mbox{\raisebox{-.4ex}{$N^{\prime}_{\sigma^{\prime}_1}$}}\;
      \longrightarrow\;
    N^{\prime}/
       \mbox{\raisebox{-.4ex}{$N^{\prime}_{\sigma^{\prime}_1}$}}\;
      \longrightarrow\;
    N^{\prime}/
       \mbox{\raisebox{-.4ex}{$N^{\prime}_{\sigma^{\prime}_2}$}}\;
      \longrightarrow\; 0\,.
$$
Note that since $N^{\prime}_{\sigma^{\prime}_2}$ is mapped by
$\overline{\varphi}^{\sigma^{\prime}_2}_{\sigma}$ to the $0$
in $N/\mbox{\raisebox{-.4ex}{$N_{\sigma}$}}$,
$N^{\prime}_{\sigma^{\prime}_2}/
  \mbox{\raisebox{-.4ex}{$N^{\prime}_{\sigma^{\prime}_1}$}}$
as a subspace in
$N^{\prime}/
  \mbox{\raisebox{-.4ex}{$N^{\prime}_{\sigma^{\prime}_1}$}}$
is mapped by $\overline{\varphi}^{\sigma^{\prime}_1}_{\sigma}$
also to the $0$ in $N/\mbox{\raisebox{-.4ex}{$N_{\sigma}$}}$.
Since
$N^{\prime}/
  \mbox{\raisebox{-.4ex}{$N^{\prime}_{\sigma^{\prime}_2}$}}\,
 =\, (N^{\prime}/
         \mbox{\raisebox{-.4ex}{$N^{\prime}_{\sigma^{\prime}_1}$}})/
     (N^{\prime}_{\sigma^{\prime}_2}/
         \mbox{\raisebox{-.4ex}{$N^{\prime}_{\sigma^{\prime}_1}$}})$
canonically, $\overline{\varphi}^{\sigma^{\prime}_1}_{\sigma}$
descends to $\overline{\varphi}^{\sigma^{\prime}_2}_{\sigma}$
also canonically after taking the quotient by
$N^{\prime}_{\sigma^{\prime}_2}/
   \mbox{\raisebox{-.4ex}{$N^{\prime}_{\sigma^{\prime}_1}$}}$.
This shows that
$\overline{\varphi}^{\sigma^{\prime}_1}_{\sigma}
  (N^{\prime}/\mbox{\raisebox{-.4ex}{$N^{\prime}_{\sigma^{\prime}_1}$}})
     =\overline{\varphi}^{\sigma^{\prime}_2}_{\sigma}
  (N^{\prime}/\mbox{\raisebox{-.4ex}{$N^{\prime}_{\sigma^{\prime}_2}$}})$
if $\sigma^{\prime}_1$ is in the face of $\sigma^{\prime}_2$.

For general
$\sigma^{\prime}_1,\,\sigma^{\prime}_2\in\Sigma^{\prime}_{\sigma}$,
let $\stackrel{\circ}{\sigma}$ be the interior of $\sigma$ and
regard $\varphi$ as a map
$N^{\prime}_{\scriptsizeBbb R}\rightarrow N_{\scriptsizeBbb R}$.
Then, since $\varphi^{-1}(\stackrel{\circ}{\sigma})$ is connected,
one can find a finite sequence of cones
$\sigma^{\prime}_1,\,\tau^{\prime}_1,\,\cdots,\,
                     \tau^{\prime}_{s},\,\sigma^{\prime}_2$
in $\Sigma^{\prime}_{\sigma}$ such that each element in the sequence
either contains as a face or is in the face of the next element.
By transitivity, this shows that
$\overline{\varphi}^{\sigma^{\prime}_1}_{\sigma}
 (N^{\prime}/\mbox{\raisebox{-.4ex}{$N^{\prime}_{\sigma^{\prime}_1}$}})
   =\overline{\varphi}^{\sigma^{\prime}_2}_{\sigma}
 (N^{\prime}/\mbox{\raisebox{-.4ex}{$N^{\prime}_{\sigma^{\prime}_2}$}})$
and we conclude the proof.

\noindent\hspace{12cm} $\Box$

\bigskip

\noindent
{\bf Lemma 2.1.8.} {\it
 If $\tau$ is in the face of $\sigma$, then
 $\Ind(\tau)$ is an integral multiple of $\Ind(\sigma)$.
} 

\bigskip

\noindent
{\it Proof.}
Since $\tau$ is in the face of $\sigma$, $\varphi^{-1}(\tau)$ must
lie in the boundary of $\varphi^{-1}(\sigma)$. Since both
$\varphi^{-1}(\tau)$ and $\varphi^{-1}(\sigma)$ are unions of cones,
each cone in $\varphi^{-1}(\tau)$, in particular each
$\tau^{\prime}\in\Sigma^{\prime}_{\tau}$,
must lie in the face of some cone $\sigma^{\prime}$ in
$\Sigma^{\prime}_{\sigma}$.
Consider now the following diagram of abelian groups with exact rows:
$$
 \begin{array}{cccccccccc}
  0 & \rightarrow
    & N^{\prime}_{\sigma^{\prime}}/
           \mbox{\raisebox{-.4ex}{$N^{\prime}_{\tau^{\prime}}$}}
    & \longrightarrow
    & N^{\prime}/
           \mbox{\raisebox{-.4ex}{$N^{\prime}_{\tau^{\prime}}$}}
    & \longrightarrow
    & N^{\prime}/
           \mbox{\raisebox{-.4ex}{$N^{\prime}_{\sigma^{\prime}}$}}
    & \longrightarrow  & 0 &                               \\[1ex]
  & & \alpha\downarrow \hspace{1em} &
    & \overline{\varphi}^{\tau^{\prime}}_{\tau}
            \downarrow \hspace{2em} &
    & \overline{\varphi}^{\sigma^{\prime}}_{\sigma}
            \downarrow \hspace{2em} &&& \\[1ex]
  0 & \rightarrow
    & N_{\sigma}/
           \mbox{\raisebox{-.4ex}{$N_{\tau}$}}
    & \longrightarrow 
    & N/\mbox{\raisebox{-.4ex}{$N_{\tau}$}}
    & \longrightarrow
    & N/\mbox{\raisebox{-.4ex}{$N_{\sigma}$}}
    & \longrightarrow  & 0  &.
 \end{array}
$$
Note that the diagram commutes since all the homomorphisms involved
are induced from the morphism $\varphi:N^{\prime}\rightarrow N$.
By the Snake Lemma in homological algebra ([Ei]), one has the exact
sequence
$$
 0\; \rightarrow\; \kker\alpha\; \rightarrow\;
     \kker\overline{\varphi}^{\sigma^{\prime}}_{\sigma}\;
     \rightarrow\;
     \kker\overline{\varphi}^{\sigma^{\prime}}_{\sigma}\;
     \rightarrow\; \coker\alpha\; \rightarrow\;
     \coker\overline{\varphi}^{\tau^{\prime}}_{\tau}\;
     \rightarrow\;
     \coker\overline{\varphi}^{\sigma^{\prime}}_{\sigma}\;
     \rightarrow\; 0\,,
$$
which implies that $\Ind(\sigma)$ divides $\Ind(\tau)$ by considering
the cardinality of the last three terms. This concludes the proof.

\noindent\hspace{12cm} $\Box$

\bigskip

We are now ready to understand how the toroidal pieces from the
orbit-to-orbit fibrations are joined together to form a complete
fiber.

Let $\sigma\in\Sigma$. By assumption, $\Sigma^{\prime}_{\sigma}$ is
non-empty. We shall regard $\Sigma^{\prime}_{\sigma}$ as a subfan in
$\Sigma^{\prime}$. Recall the construction of the star of a cone in
toric geometry. Here we need a relative version defined as follows.

Let $\tau^{\prime}\in\Sigma^{\prime}_{\sigma}$ and
$\{\,\sigma^{\prime}_1,\,\sigma^{\prime}_2,\,\cdots\,\}$
be the set of cones in $\Sigma^{\prime}_{\sigma}$ that contains
$\tau^{\prime}$ as a face. Then each $\sigma^{\prime}_i$ determines
a cone $\overline{\sigma^{\prime}_i}$ in
$\varphi^{-1}((N_{\sigma})_{\scriptsizeBbb R})/
   \mbox{\raisebox{-.4ex}{$(N^{\prime}_{\tau^{\prime}}
                                    )_{\scriptsizeBbb R}$}}$,
defined by
$$
 \overline{\sigma^{\prime}_i}\;
 =\; (\, \sigma^{\prime}_i+
      (N^{\prime}_{\tau^{\prime}})_{\scriptsizeBbb R}\,)/
           \mbox{\raisebox{-.4ex}{$(N^{\prime}_{\tau^{\prime}}
                                           )_{\scriptsizeBbb R}$}}\,.
$$
Note that
$\sigma^{\prime}_i+(N^{\prime}_{\tau^{\prime}})_{\scriptsizeBbb R}$
is contained in $\varphi^{-1}((N_{\sigma})_{\scriptsizeBbb R})$ since
$\tau^{\prime},\, \sigma^{\prime}_i\in \Sigma^{\prime}_{\sigma}$.
Thus,
$\{\,\overline{\sigma^{\prime}_1},\,
                   \overline{\sigma^{\prime}_2},\,\cdots\,\}$
defines a fan in
$\varphi^{-1}((N_{\sigma})_{\scriptsizeBbb R})/
   \mbox{\raisebox{-.4ex}{$(N^{\prime}_{\tau^{\prime}}
                                       )_{\scriptsizeBbb R}$}}$.

\bigskip

\noindent
{\bf Definition 2.1.9 [relative star].} {
 The fan in
 $\varphi^{-1}((N_{\sigma})_{\scriptsizeBbb R})/
   \mbox{\raisebox{-.4ex}{$(N^{\prime}_{\tau^{\prime}}
                                     )_{\scriptsizeBbb R}$}}$
 constructed above will be called the
 {\it relative star of $\tau^{\prime}$ over $\sigma$}
 and will be denoted by $\Star_{\sigma}(\tau^{\prime})$.
} 

\bigskip

\noindent
{\bf Definition 2.1.10 [primitive cone].} {
 A cone $\tau^{\prime}\in \Sigma^{\prime}_{\sigma}$ is called
 {\it primitive} with respect to $\varphi$ if none of the faces
 of $\tau^{\prime}$ are in $\Sigma^{\prime}_{\sigma}$.
 The set of primitive cones in $\Sigma^{\prime}_{\sigma}$ will be
 denoted by $\Sigma^{\prime\,\circ}_{\sigma}$.
} 

\bigskip
\noindent
(Cf.\ {\sc Figure 2-1-3}.)
 \begin{figure}[htbp]
  \setcaption{{\sc Figure 2-1-3.}
  \baselineskip 14pt
  A comparison of the relative star $\Star_{\sigma}(\tau^{\prime})$
  with respect to a map of fans $\varphi$ and the usual star
  $\Star(\tau^{\prime})$.
  } 
  \centerline{\psfig{figure=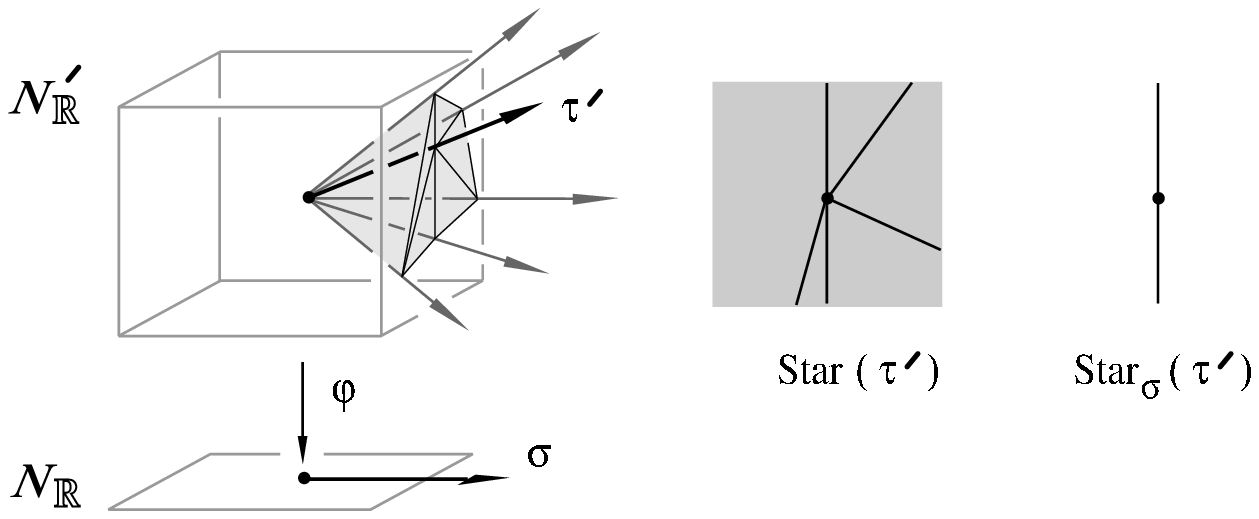,width=11cm,caption=}}
 \end{figure}

The following lemma justifies the two definitions above.

\bigskip

\noindent
{\bf Lemma 2.1.11 [connected/irreducible components of fiber].} {\it
 Given $y\in O_{\sigma}$, let $F_{\sigma}^c$ be a connected component
 of $\widetilde{\varphi}^{-1}(y)$. Then
 \begin{quote}
  \hspace{-1.9em}(a)\hspace{1ex}
  There is a $1-1$ correspondence between
  the irreducible components of the fiber $F_{\sigma}^c$ over
  $O_{\sigma}$ and the elements in $\Sigma^{\prime\,\circ}_{\sigma}$.

  \hspace{-1.9em}(b)\hspace{1ex}
  Let $\tau^{\prime}\in\Sigma^{\prime\,\circ}_{\sigma}$, then the
  component $F_{\sigma}^{\tau^{\prime}}$ of $F_{\sigma}^c$ associated
  to $\tau^{\prime}$ is the toric variety associated to the fan
  $\Star_{\sigma}(\tau^{\prime})$ in
  $\varphi^{-1}((N_{\sigma})_{\scriptsizeBbb R})/
     \mbox{\raisebox{-.4ex}{$(N^{\prime}_{\tau^{\prime}}
                                          )_{\scriptsizeBbb R}$}}$.
  It has (complex) dimension$\,$: $\codim\tau^{\prime}-\codim\sigma$.
  
  \hspace{-1.9em}(c)\hspace{1ex}
  If
  $\tau^{\prime}_1,\,\cdots,\, \tau^{\prime}_s
                       \in\Sigma^{\prime\,\circ}_{\sigma}$
  generate a cone $\sigma^{\prime}$ in $\Sigma^{\prime}$, then
  $\sigma^{\prime}\in\Sigma^{\prime}_{\sigma}$ and  \newline 
  $\;\bigcap\,_{i=1}^s\,F_{\sigma}^{\tau^{\prime}_i}\,
                            =\,\Star_{\sigma}(\sigma^{\prime})$.

  \hspace{-1.9em}(d)\hspace{1ex}
  The whole fiber $F_{\sigma}$ over a point in $O_{\sigma}$ is the
  disjoint union of $\Ind(\sigma)$-many copies of $F_{\sigma}^c$.
 \end{quote}
} 

\bigskip

\noindent
{\it Proof.}
{\it (a) The $1-1$ correspondence$\,$}:
By definition, any non-primitive element $\sigma^{\prime}$ in
$\Sigma^{\prime}_{\sigma}$ must contain some primitive element
$\tau^{\prime}$ in its face. This implies that
$F_{\sigma}^{\,\sigma^{\prime}}$ is in the closure of
$F_{\sigma}^{\,\tau^{\prime}}$. On the other hand, for a primitive
$\tau^{\prime}$ in $\Sigma^{\prime}_{\sigma}$ the fiber
$F_{\sigma}^{\tau^{\prime}}$ cannot be contained in the closure
of $F_{\sigma}^{\sigma^{\prime}}$ for any
$\sigma^{\prime}\ne\tau^{\prime}$. This shows that
$F_{\sigma}^{\tau^{\prime}}$ is a maximal stratum in $F_{\sigma}^c$
and, hence, its closure must be a component of $F_{\sigma}^c$.
Together, these prove the $1-1$ corrrespondence. 

\medskip

\noindent
{\it (b) The toric geometry of a component$\,$}:
Let $\tau^{\prime}$ be a primitive element in
$\Sigma^{\prime}_{\sigma}$. Then the closure
$\overline{F_{\sigma}^{\tau^{\prime}}}$ consists of
$F_{\sigma}^{\tau^{\prime}}$ and all $F_{\sigma}^{\sigma^{\prime}}$
with $\sigma^{\prime}$ containing $\tau^{\prime}$ in its face.
By construction, this is exactly the fan
$\Star_{\sigma}(\tau^{\prime})$ in
$\varphi^{-1}((N_{\sigma})_{\scriptsizeBbb R})/
   \mbox{\raisebox{-.4ex}{$(N^{\prime}_{\tau^{\prime}}
                                          )_{\scriptsizeBbb R}$}}$,
which has dimension $\codim\tau^{\prime}-\codim\sigma$.
This proves Part (b).

\medskip

\noindent
{\it (c) The intersection of irreducible components$\,$}:
This follows from the definition of relative stars and the
order-reversing correspondence between cones $\sigma^{\prime}$ in
$\Sigma^{\prime}_{\sigma}$ and toric varieties associated to
$\Star_{\sigma}(\sigma^{\prime})$.

\medskip

\noindent
{\it (d) The number of connected components$\,$}:
Note that the fiber $F_{\sigma}$ over a point in $O_{\sigma}$ is
naturally stratified by complex tori of the form $({\Bbb C}^{\ast})^r$.
Lemma 2.1.6 implies that each $\sigma^{\prime}$ contributes to
$F_{\sigma}$ exactly $\Ind(\sigma)$-many strata $({\Bbb C}^{\ast})^r$,
where $r=\codim\sigma^{\prime}-\codim\sigma_{\sigma}$.
Each connected components of $F_{\sigma}$ will contain at least one of
them. On the other hand, the restriction of the map
$\widetilde{\varphi}:
   \widetilde{\varphi}^{-1}(O_{\sigma})\rightarrow O_{\sigma}$
factors through a covering morphism
$({\Bbb C}^{\ast})\rightarrow O_{\sigma}$ of degree $\Ind(\sigma)$.
Together this shows that $F_{\sigma}$ must have exactly
$\Ind(\sigma)$-many connected components.

\medskip

This concludes the proof.

\noindent\hspace{12cm} $\Box$

\bigskip

This justifies the description of fibers of a toric morphism in
Proposition 2.1.4.

\bigskip

Finally, let us consider the induced stratification of
$\widetilde{\varphi}$ by toric orbits $O_{\sigma}$ of $X_{\Sigma}$. 
Observe first the following commutative diagram of maps of fans
$$
 \begin{array}{cccccl}
  \Sigma^{\prime}_{\{0\}}
      & \hookrightarrow  & \Sigma^{\prime}
      & =                & \Sigma^{\prime} & \\[.6ex]
  \hspace{1ex}\downarrow \mbox{\scriptsize $\varphi$}
      & & \hspace{2ex}\downarrow \mbox{\scriptsize $\varphi$}
      & & \hspace{2ex}\downarrow \mbox{\scriptsize $\varphi$}  \\[.6ex]
  \hspace{2ex}\{0\}_{\varphi(N^{\prime})}
      & \hookrightarrow  & \hspace{2ex}\Sigma_{\varphi(N^{\prime})}
      & \rightarrow      & \Sigma          &,
 \end{array}
$$
where $\{0\}_{\varphi(N^{\prime})}$ is the fan in
$\varphi(N^{\prime})_{\scriptsizeBbb R}$ that consists only of $\{0\}$,
$\Sigma^{\prime}_{\{0\}}$ is regarded as a fan in
$N^{\prime}_{\scriptsizeBbb R}$, and $\Sigma_{\varphi(N^{\prime})}$ is
the naturally induced fan in $\varphi(N^{\prime})_{\scriptsizeBbb R}$
from $\Sigma$ in $N_{\scriptsizeBbb R}$.
By taking the toric variety associated to the fans in the diagram,
one obtains the torically equivariant commutative diagram
$$
 \begin{array}{cccl}
  \widetilde{O_{\{0\}}}\times F_{\{0\}}^c
   & \hookrightarrow  & X_{\Sigma^{\prime}}         & \\[.6ex]
  \hspace{1ex}\downarrow \mbox{\scriptsize $\widetilde{\varphi}$}
   & & \hspace{2ex}\downarrow \mbox{\scriptsize $\widetilde{\varphi}$}
                                                    & \\[.6ex]
  O_{\{0\}}  & \hookrightarrow  & X_{\Sigma}        & ,
 \end{array}
$$
where $O_{\{0\}}$ is the open dence $T_N$-orbit in $X_{\Sigma}$ and 
$\widetilde{O_{\{0\}}}
  \stackrel{\widetilde{\varphi}}{\rightarrow} O_{\{0\}}$
is a covering map of order $[N:\varphi(N^{\prime})]=\Ind(\{0\})$
This shows that
$\widetilde{\varphi}^{-1}(O_{\{0\}})
  =\widetilde{O_{\{0\}}}\times F_{\{0\}}^c$,
as claimed.
For the stratum over general orbit $O_{\sigma}$, $\sigma\in\Sigma$,
let $\tau^{\prime}$ be a primitive element in
$\Sigma^{\prime}_{\sigma}$. Then the discussion above applied to
the following commutative diagram of maps of fans
$$
 \begin{array}{ccccc}
  \Star_{\sigma}(\tau^{\prime})
      & \hookrightarrow  & \Star(\tau^{\prime})
      & =                & \Star(\tau^{\prime})    \\[.6ex]
  \downarrow & & \downarrow     & & \downarrow     \\[.6ex]
  \{0\}      & \hookrightarrow  & \Star(\sigma)
             & \rightarrow      & \Star(\sigma)          
 \end{array}
$$
from the maps of lattices induced by $\varphi$ (with the fans in the
above diagram lying in the realization of the lattices in the diagram
below at the same relative positions)
$$
 \begin{array}{ccccc}
  N^{\prime}/\mbox{\raisebox{-.4ex}{$N^{\prime}_{\tau^{\prime}}$}}
   & =
   & N^{\prime}/\mbox{\raisebox{-.4ex}{$N^{\prime}_{\tau^{\prime}}$}}
   & =
   & N^{\prime}/\mbox{\raisebox{-.4ex}{$N^{\prime}_{\tau^{\prime}}$}}
                     \\[.6ex]
  \hspace{1ex}\downarrow
   \mbox{\scriptsize $\overline{\varphi}^{\tau^{\prime}}_{\sigma}$}
   & & \hspace{2ex}\downarrow
       \mbox{\scriptsize $\overline{\varphi}^{\tau^{\prime}}_{\sigma}$}
   & & \hspace{2ex}\downarrow
       \mbox{\scriptsize $\overline{\varphi}^{\tau^{\prime}}_{\sigma}$}
                                                                \\[.6ex]
  \overline{\varphi}^{\tau^{\prime}}_{\sigma}(
   N^{\prime}/\mbox{\raisebox{-.4ex}{$N^{\prime}_{\tau^{\prime}}$}})
   & =
   & \overline{\varphi}^{\tau^{\prime}}_{\sigma}(
      N^{\prime}/\mbox{\raisebox{-.4ex}{$N^{\prime}_{\tau^{\prime}}$}})
   & \hookrightarrow
      & N/\mbox{\raisebox{-.4ex}{$N_{\sigma}$}}
 \end{array}
$$
shows that $\widetilde{O_{\sigma}}\times F^{\tau^{\prime}}_{\sigma}$
is contained in $\widetilde{\varphi}^{-1}(O_{\sigma})$, where
$\widetilde{O_{\sigma}}\rightarrow O_{\sigma}$ is a covering map
of order $\Ind(\sigma)$. 
From the pasting of the irreducible components
$F^{\tau^{\prime}}_{\sigma}$ to form $F_{\sigma}^c$ discussed earlier,
one concludes that 
$\widetilde{\varphi}^{-1}(O_{\sigma})
             =\widetilde{O_{\sigma}}\times F_{\sigma}^c$.

\bigskip

This concludes the proof of Proposition 2.1.4.

\noindent\hspace{12cm} $\Box$

\bigskip

\noindent
{\it Remark 2.1.12.} If $[N:\varphi(N^{\prime})]<\infty$, then
$$
 \Ind(\{0\})\;
  =\; [N:\varphi(N^{\prime})]\;
  =\; \Ind(\sigma)\,
     \cdot\, [N_{\sigma}:N_{\sigma}\cap\varphi(N^{\prime})]\,.
$$
In particular, if $\varphi$ is surjective, then $\Ind(\sigma)=1$
for all $\sigma\in\Sigma$.

\bigskip

\begin{flushleft}
{\bf Characterization of toric fibrations.}
\end{flushleft}
To characterize a fibration, let us assume that
$\widetilde{\varphi}:X_{\Sigma^{\prime}}\rightarrow X_{\Sigma}$
is surjective. Let $F^{\tau^{\prime}}_{\sigma}$ be an irreducible
component of the connected fiber $F^c_{\sigma}$ over
$O_{\sigma}\subset X_{\Sigma}$. Then, from Lemma 2.1.11,  
\begin{eqnarray*}
 \lefteqn{
   \dimm F^{\tau^{\prime}}_{\sigma}\;
    =\; \dimm \left(\varphi^{-1}((N_{\sigma})_{\scriptsizeBbb R})/
           \mbox{\raisebox{-.4ex}{$(N^{\prime}_{\tau^{\prime}}
                                  )_{\scriptsizeBbb R}$}}\right)\;
    =\; \dimm\varphi^{-1}(0) + \dimm (N_{\sigma})_{\scriptsizeBbb R}
     - \dimm (N^{\prime}_{\tau^{\prime}})_{\scriptsizeBbb R} } \\[.6ex]
  & &
   \ge \dimm\varphi^{-1}(0)\;
    =\; \dimm F_0\; = \dimm X_{\Sigma^{\prime}}-\dimm X_{\Sigma}\,,
        \hspace{12em}
\end{eqnarray*}
where $F^c_0$ is a connected component of a generic fiber and 
we use the fact that, for $\tau^{\prime}$ primitive, the map
$\varphi_{\scriptsizeBbb R}:
     (N^{\prime}_{\tau^{\prime}})_{\scriptsizeBbb R}
                  \rightarrow (N_{\sigma})_{\scriptsizeBbb R}$ 
is injective. From this, one concludes$\,$:

\bigskip

\noindent
{\bf Corollary 2.1.13 [fibration].} {\it
 Let $\widetilde{\varphi}:X_{\Sigma^{\prime}}\rightarrow X_{\Sigma}$
 is a surjective toric morphism. Then $\widetilde{\varphi}$ is a
 fibration, in the sense that all the irreducible components of
 fibers of $\widetilde{\varphi}$ have the same dimension,
 $\dimm X_{\Sigma^{\prime}}-\dimm X_{\Sigma}$, if and only if
 the map of cones
 $\varphi_{\scriptsizeBbb R}:\tau^{\prime}\rightarrow \sigma$
 is a bijection for every
 $\tau^{\prime}\in\Sigma^{\prime\,\circ}_{\sigma}$.
} 

\bigskip

\noindent
In particular, since a ray in $\Sigma^{\prime}(1)-\varphi^{-1}(0)$
must be a primitive cone, if $\widetilde{\varphi}$ is a fibration,
then $\varphi$ maps $\Sigma^{\prime}(1)$ to $\Sigma(1)$ surjectively.

\bigskip

\begin{flushleft}
{\bf Dual picture via polytopes in $M$ lattice.}
\end{flushleft}
As the above proposition indicates, the study of toric morphisms via
fans in the $N$ lattice is very natural. Nevertheless, for projective
toric varieties, it is convenient to have the dual polyhedral picture
in the $M$ lattice in mind. We summarize some essential features here.
They follow directly from the proof of Proposition 2.1.4 by taking dual
and tracing definitions.

Let $\widetilde{\varphi}:X_{\Sigma^{\prime}}\rightarrow X_{\Sigma}$
be a toric morphism between projective toric varieties.
Without loss of generality, we may assume that
$\varphi_{\scriptsizeBbb R}:N^{\prime}_{\scriptsizeBbb R}
                                  \rightarrow N_{\scriptsizeBbb R}$
is surjective. Note that this implies that
$\varphi^{\dagger}:M\rightarrow M^{\prime}$ is injective.
Let $\Delta^{\prime}$ (resp.\ $\Delta$) be a convex polytope in
$M^{\prime}_{\scriptsizeBbb R}$ whose normal fan is $\Sigma^{\prime}$
(resp.\ $\Sigma$). The face of $\Delta$ whose normal cone is $\sigma$
will be denoted by $\Theta_{\sigma}$. Note that, by definition, the
face associated to $\{0\}$ in $N_{\scriptsizeBbb R}$ is $\Delta$.
Similarly for the face $\Theta_{\sigma^{\prime}}$ of $\Delta^{\prime}$. 

The following equivalences play a key role in translating the fan
picture to the dual polytope picture:
\begin{eqnarray*}
  & & \varphi(\,\mbox{interior of $\sigma^{\prime}$}\,)\,
                     \subset\, \mbox{interior of $\sigma$.}\\
  & \Longleftrightarrow
    & \langle\,u,\,\varphi(v^{\prime})\,\rangle >0 \hspace{1em}
       \mbox{for all $u\in\sigma^{\vee}-\sigma^{\perp}$ and
             $v^{\prime}$ in the interior of $\sigma^{\prime}$.} \\
  & \Longleftrightarrow
    & \langle\,\varphi^{\dagger}(u),\,v^{\prime}\,\rangle >0 \hspace{1em}
       \mbox{for all $u\in\sigma^{\vee}-\sigma^{\perp}$ and
             $v^{\prime}$ in the interior of $\sigma^{\prime}$.} \\
  & \Longleftrightarrow
    & \varphi^{\dagger}(\sigma^{\vee}-\sigma^{\perp})\,
        \subset\, {\sigma^{\prime}}^{\vee}-{\sigma^{\prime}}^{\perp}\,.
        \hspace{1em}\mbox{($\,$Note that 
                           $\varphi^{\dagger}(\sigma^{\perp})\,
                            \subset\,{\sigma^{\prime}}^{\perp}$.$\,$)}
\end{eqnarray*}
For $\sigma\in\Sigma$, define the {\it $\sigma$-lighted part}
$\Delta^{\prime}_{\sigma}\subset\Delta^{\prime}$ associated to a face
$\Theta_{\sigma}\subset\Delta$ by
$$
 \Delta^{\prime}_{\sigma}\;
  =\;\{\, \Theta_{\sigma^{\prime}}\in\Delta^{\prime}\,|\,
            \varphi^{\dagger}(\sigma^{\vee}-\sigma^{\perp})\,
             \subset\,
              {\sigma^{\prime}}^{\vee}-{\sigma^{\prime}}^{\perp} \,\}\,.
$$
(The terminology is inspired by {\sc Figure 2-1-4}.)
Note that, for an interior point $p\in\Theta_{\sigma^{\prime}}$,
the {\it unridged tangent cone}
$T_p\Delta^{\prime}-T_p\Theta_{\sigma^{\prime}}$ at $p$ is
isomorphic to ${\sigma^{\prime}}^{\vee}-{\sigma^{\prime}}^{\perp}$.
Thus, $\Delta^{\prime}_{\sigma}$ is the polyhedral sub-complex of
$\Delta^{\prime}$ that consists of faces, the unridged tangent cone
at an interior point of which contains 
$\varphi^{\dagger}(\sigma^{\vee}-\sigma^{\perp})$.
By definition, $\Theta_{\sigma^{\prime}}\in\Delta^{\prime}_{\sigma}$
if and only if the interior of $\sigma^{\prime}$ is mapped by $\varphi$
to the interior of $\sigma$.
If one regards $\Delta^{\prime}$ as a collection of combinatorial
polytopes related by inclusion relations, then $\Delta^{\prime}$ is
decomposed into a disjoint union
$\sqcup_{\sigma\in\Sigma}\,\Delta^{\prime}_{\sigma}$.
({\sc Figure 2-1-4}.)
\begin{figure}[htbp]
 \setcaption{{\sc Figure 2-1-4.}
  \baselineskip 14pt
  The $\sigma$-lighted part
  $\Delta^{\prime}_{\sigma}\subset\Delta^{\prime}$ associated to
  a face $\Theta_{\sigma}\subset\Delta$ is indicated by the shaded
  polyhedral subcomplex. A cross-section to indicate why these
  faces are $\sigma$-lighted is also indicated.
  (Cf.\ Forward light cones, domain of influence, ...,
        in general relativity.)
 } 
 \centerline{\psfig{figure=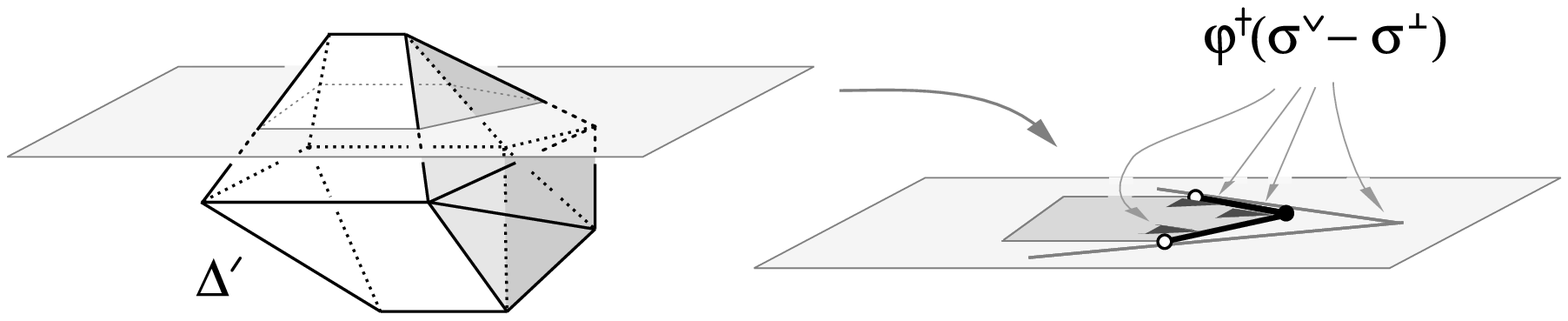,width=11cm,caption=}}
\end{figure}

\noindent
{\bf Lemma 2.1.14 [polytope for relative star].} {\it 
Let $\Theta_{\tau^{\prime}}\in\Delta^{\prime}_{\sigma}$, regarded as a
polytope in ${\tau^{\prime}}^{\perp}$. Let
$\overline{\varphi^{\dagger}}^{\tau^{\prime}}_{\sigma}\,:\,
   {\tau^{\prime}}^{\perp}\, \rightarrow\,
     {\tau^{\prime}}^{\perp}/
        \mbox{\raisebox{-.4ex}{$\varphi^{\dagger}(\sigma^{\perp})$}}$
be the quotient map. Then the normal fan of the polytope 
$\overline{\Theta_{\tau^{\prime}}}
   =\overline{\varphi^{\dagger}}^{\tau^{\prime}}_{\sigma}
                                            (\Theta_{\tau^{\prime}})$
in
${\tau^{\prime}}^{\perp}/
    \mbox{\raisebox{-.4ex}{$\varphi^{\dagger}(\sigma^{\perp})$}}$
is the relative star $\Star_{\sigma}(\tau^{\prime})$.}
({\sc Figure 2-1-5}.)

\begin{figure}[htbp]
\setcaption{{\sc Figure 2-1-5.}
  \baselineskip 14pt
  The polytope $\overline{\Theta_{\tau^{\prime}}}$
  in 
  ${\tau^{\prime}}^{\perp}/
      \mbox{\raisebox{-.4ex}{$\varphi^{\dagger}(\sigma^{\perp})$}}$
  whose normal fan is
  the relative star $\Star_{\sigma}(\tau^{\prime})$ is indicated.
  (In the picture, $\varphi^{\dagger}(\sigma^{\perp})$ is realized
   as a foliation of ${\tau^{\prime}}^{\perp}$.)
 } 
 \centerline{\psfig{figure=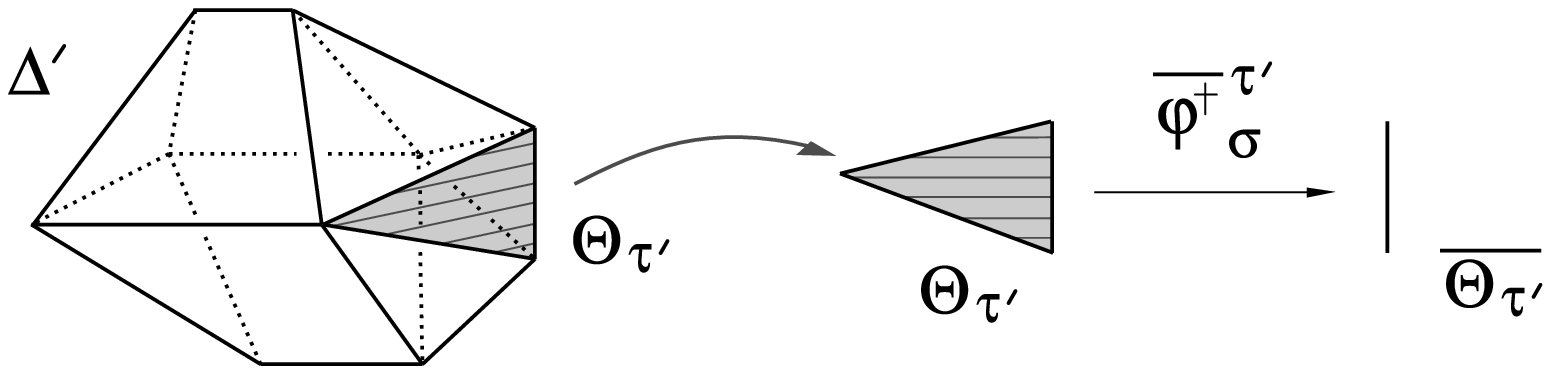,width=11cm,caption=}}
 
\end{figure}

\bigskip

\noindent
{\it Proof.}
Recall Definition 2.1.9 of relative stars.
Let $\tau^{\prime}\in\Sigma^{\prime}_{\sigma}$; then there is
a non-degenerate pairing
$\varphi^{-1}(N_{\sigma})/
    \mbox{\raisebox{-.4ex}{$N^{\prime}_{\tau^{\prime}}$}}\,
 \times\,
 {\tau^{\prime}}^{\perp}/
     \mbox{\raisebox{-.4ex}{$\varphi^{\dagger}(\sigma^{\perp})$}}\,
 \longrightarrow {\Bbb Z}$ 
induced by the restriction of the canonical pairing
$N^{\prime}\times M^{\prime}\longrightarrow {\Bbb Z}$ to
$\varphi^{-1}(N_{\sigma})\times{\tau^{\prime}}^{\perp}
                                         \longrightarrow {\Bbb Z}$.
Consequently, if $\sigma^{\prime}\in\Sigma^{\prime}_{\sigma}$
contains $\tau^{\prime}$ in its face, then the dual cone of 
$\overline{\sigma^{\prime}}\,
 =\, (\, \sigma^{\prime}+
      (N^{\prime}_{\tau^{\prime}})_{\scriptsizeBbb R}\,)/
           \mbox{\raisebox{-.4ex}{$(N^{\prime}_{\tau^{\prime}}
                                           )_{\scriptsizeBbb R}$}}
 \subset
 \varphi^{-1}((N_{\sigma})_{\scriptsizeBbb R})/
   \mbox{\raisebox{-.4ex}{$(N^{\prime}_{\tau^{\prime}}
                                       )_{\scriptsizeBbb R}$}}$
is simply
$(\,({\sigma^{\prime}}^{\vee}\cap{\tau^{\prime}}^{\perp})
                     +\varphi^{\dagger}(\sigma^{\perp})\,)/
        \mbox{\raisebox{-.4ex}{$\varphi^{\dagger}(\sigma^{\perp})$}}
  \subset
    {\tau^{\prime}}^{\perp}/
        \mbox{\raisebox{-.4ex}{$\varphi^{\dagger}(\sigma^{\perp})$}}$,
which is the tangent cone of $\overline{\Theta_{\tau^{\prime}}}$
at an interior point of the face
$\overline{\varphi^{\dagger}}^{\tau^{\prime}}_{\sigma}
                                    (\Theta_{\sigma^{\prime}})
                          \subset \overline{\Theta_{\tau^{\prime}}}$.
This gives a correspondence from $\Star_{\sigma}(\tau^{\prime})$
to the set of faces of $\overline{\Theta_{\tau^{\prime}}}$.
This also implies that the normal fan of
$\overline{\Theta_{\tau^{\prime}}}$ contains
$\Star_{\sigma}(\tau^{\prime})$. Since both fans are complete, they
must be identical. This concludes the proof.

\noindent\hspace{12cm} $\Box$

\bigskip

Similarly, one has the dual characterization for primitive elements
$\sigma^{\prime}$ in $\Sigma^{\prime}_{\sigma}$:

\bigskip

\noindent
{\bf Definition 2.1.15 [primitive face].} {
 A face $\Theta_{\sigma^{\prime}}$ in $\Delta^{\prime}_{\sigma}$ is
 called {\it primitive} if it is not in the boundary of any other face
 in $\Delta^{\prime}_{\sigma}$.
} 

\bigskip

\noindent
One can check that $\sigma^{\prime}$ is a primitive elememt in
$\Sigma^{\prime}_{\sigma}$ if and only if $\Theta_{\sigma^{\prime}}$
is a primitive face of $\Delta^{\prime}_{\sigma}$.
The following corollary follows immediately from Lemma 2.1.11 and
Lemma 2.1.14 (cf.\ the notations in Lemma 2.1.11).

\bigskip

\noindent
{\bf Corollary 2.1.16 [polytope for irreducible component of fiber].}
{\it
 (a) There is a $1-1$ correspondence between the primitive faces in
     $\Delta^{\prime}_{\sigma}$ and the irreducible components of
     the fiber $F_{\sigma}^c$ of $\widetilde{\varphi}$ over a point
     in $O_{\sigma}$.
 (b) Let $\Theta_{\tau^{\prime}}$ be a primitive
     face of $\Delta^{\prime}_{\sigma}$. Then the corresponding
     irreducible component $F^{\tau^{\prime}}_{\sigma}$ in
     $F^c_{\sigma}$ is the toric variety associated to the normal
     fan of the polytope $\overline{\Theta_{\tau^{\prime}}}$
     in
     ${\tau^{\prime}}^{\perp}/
       \mbox{\raisebox{-.4ex}{$\varphi^{\dagger}(\sigma^{\perp})$}}$.
} 

\bigskip

\noindent
In particular, for the generic fiber of the morphism
$\widetilde{\varphi}$, $\{0\}^{\vee}=M$ and $\Delta^{\prime}_{\{0\}}$
consists of $\Delta^{\prime}$ and all the faces
$\Theta_{\sigma^{\prime}}$ of $\Delta^{\prime}$ with
${\sigma^{\prime}}^{\perp}\supset\varphi^{\dagger}(M)$.
The normal fan of the projection of $\Delta^{\prime}_{\{0\}}$ in
$\left(M^{\prime}/\mbox{\raisebox{-.4ex}{$\varphi^{\dagger}(M)$}}
                                          \right)_{\scriptsizeBbb R}$
is exactly $\Sigma^{\prime}_{\{0\}}$ as discussed earlier.
(Cf.\ [A-K-M-S].)

\bigskip

\subsection{Examples of toric morphisms and fibers.}

In this subsection, we give two examples of toric morphisms and
 their fibers to illustrate the previous discussions.
The second example is a toric fibration from [B-C-dlO-G] and
 will be the major example for later discussions of this article.

\bigskip

\noindent
{\bf Example 2.2.1
  [$\widetilde{\varphi}: X_{\Sigma^{\prime}}\rightarrow \CP^1$].}
Let $\Sigma=\{\sigma_+, \{0\}, \sigma_-\}$ be the fan for $\CP^1$.
Then one has the corresponding decomposition
$\Sigma^{\prime}=\Sigma^{\prime}_{\sigma_+}
          \cup\Sigma^{\prime}_{\{0\}}\cup\Sigma^{\prime}_{\sigma_-}$.
Let
$\Sigma^{\prime}(1)
   =\Sigma^{\prime}_+(1)\cup\Sigma^{\prime}_0(1)\cup\Sigma_-(1)$
be the corresponding decomposition. Since $\sigma_+$ and $\sigma_-$
are $1$-dimensional and cones in $\Sigma^{\prime}$ are strongly convex,
the set of primitive elements in $\Sigma^{\prime}_{\sigma_+}$
(resp.\ $\Sigma^{\prime}_{\sigma_-}$) is exactly $\Sigma^{\prime}_+(1)$
(resp.\ $\Sigma^{\prime}_-(1)$). Also, since $\sigma_{\pm}$ is a maximal
cone, for $\sigma^{\prime}\in\Sigma^{\prime}_{\pm}(1)$ the relative
star $\Star_{\sigma_{\pm}}(\sigma^{\prime})$ coincides with the usual
star $\Star(\sigma^{\prime})$. Consequently, all the irreducible
components of fibers of $\widetilde{\varphi}$ are of codimension $1$
in $X_{\Sigma^{\prime}}$. Note that in this case $\CP^1$ has only
three distinct orbits: ${\Bbb C}^{\ast}$, $x_+$, and $x_-$. Thus,
there are at most three different kinds of fibers for
$\widetilde{\varphi}$:
\begin{quote}
 \hspace{-1.9em}(1) {\it generic fiber}$\,$:
  a connected component is given by the toric variety
  $\Sigma^{\prime}\cap\varphi^{-1}(0)$; there are 
  $[N:\varphi(N^{\prime})]$-many connected components.

 \hspace{-1.9em}(2) {\it the fiber over $x_+$}$\,$:
  Since $N/\mbox{\raisebox{-.4ex}{$N_{\sigma_+}$}}=\{0\}$, by
  Lemma 2.1.11 the fiber can have only one connected component.
  The set of irreducible components of the fiber is given then by
  $\{\Star(\sigma^{\prime})\,|\,
              \sigma^{\prime}\in\Sigma^{\prime}_{\sigma_+}\}$.
  
 \hspace{-1.9em}(3) {\it the fiber over $x_-$}$\,$:
  similar to (2) with the label $+$ replaced by $-$.
\end{quote}

For the dual picture, assume that $X_{\Sigma^{\prime}}$ is projective
and $\Sigma^{\prime}$ is the normal fan of a convex polytope
$\Delta^{\prime}$ in $M^{\prime}_{\scriptsizeBbb R}$. Then
$\Delta^{\prime}$ is decomposed into
$\Delta^{\prime}_{\sigma_+}$, $\Delta^{\prime}_{\{0\}}$, and
$\Delta^{\prime}_{\sigma_-}$. The set of primitive faces in
$\Delta^{\prime}_{\sigma_+}$ is exactly the set of maximal faces
of $\Delta^{\prime}_{\sigma_+}$.
Similarly for $\Delta^{\prime}_{\sigma_-}$.
The generic fiber of $\widetilde{\varphi}$ corresponds to the projection
of $\Delta^{\prime}$ in 
$\left(M^{\prime}/\mbox{\raisebox{-.4ex}{$\varphi^{\dagger}(M)$}}
                                          \right)_{\scriptsizeBbb R}$
For the fiber $F_{x^+}$ over $x_+$, since $\sigma_+^{\perp}=\{0\}$,
each maximal face $\Theta_{\sigma^{\prime}}$ in
$\Delta^{\prime}_{\sigma_+}$ as a polytope in
${\sigma^{\prime}}^{\perp}$ corresponds to an irreducible component of
$F_{x_+}$. Similarly for the fiber $F_{x_-}$ over $x_-$.
(Cf.\ {\sc Figure 2-2-1}.)

\begin{figure}[htbp]
 \setcaption{{\sc Figure 2-2-1.}
  \baselineskip 14pt
  A toric morphism
  $\widetilde{\varphi}:X_{\Sigma^{\prime}}\rightarrow\CP^1$
  and its generic and the two special fibers are indicated.
  The primitive faces of $\Delta^{\prime}$ are darkened/shaded.
 } 
\end{figure}

It is worth emphasizing that in this case the fiber over $x_{\pm}$
is a union of codimension-$1$ orbit closures in $X_{\Sigma^{\prime}}$;
in other words, it is a T-Weil divisor in $X_{\Sigma^{\prime}}$.

\noindent\hspace{12cm} $\Box$

\bigskip

\noindent
{\bf Example 2.2.2 [B-C-dlO-G].}
In [B-C-dlO-G], Braun, Candelas, de la Ossa, and Grassi consider
a $5$-dimensional toric variety associated to a reflexive polytope,
in which the Calabi-Yau hypersurfaces admit an elliptic fibration.
In their work, they characterize the existence of the fibration and
the generic fiber at the level of toric varieties and study also
the gauge group of the effective $4$-dimensional supersymmetric
field theory and the divisors of the Calabi-Yau hypersurfaces that
contribute to the superpotential of that effective field theory.
Here let us study their example of toric fibration in more details.
(The restriction of this fibration to Calabi-Yau hypersurfaces will
be studied in detail in Example 4.2.1 in Sec.\ 4.2.) To make the
comparison with their work easier, we shall follow their notations
in describing the fan $\Sigma^{\prime}$ with modifications to fit
into our earlier notations. The notations for the fan $\Sigma$
are modified from [Ra].

\bigskip

\noindent
{\it (a) The toric data of $X_{\Sigma^{\prime}}$ and $X_{\Sigma}$.}
The $1$-skeleton of the fans and the labels of the $1$-rays are given
in the following table$\,$:

\bigskip

\centerline{\footnotesize
\begin{tabular}{|c|c||c|c|} \hline
 \multicolumn{2}{|c||}{ \rule{0ex}{3ex} $\Sigma^{\prime}(1)$}
                   & \multicolumn{2}{c|}{$\Sigma(1)$} \\[.6ex]  \hline
 \rule{0ex}{3ex}
 $v_1^{\prime}$  &  (-1, 0, 0, 2, 3)
                 &  $d_4$   &   (-1, 0, 0) \\[.6ex]
 $v_2^{\prime}$  &  ( 0,-1, 0, 2, 3)
                 &  $d_3$   &   ( 0,-1, 0) \\[.6ex]
 $c_1^{\prime}$  &  ( 0, 0,-1, 2, 3)
                 &  $r_2$   &  ( 0, 0,-1) \\[.6ex]
 $c_2^{\prime}$  &  ( 0, 0,-1, 1, 2)
                 &  $r_1$   &  ( 0, 0, 1) \\[.6ex]
 $v_4^{\prime}$  &  ( 0, 0, 0,-1, 0)
                 &  $d_2$   &  ( 0, 1, 2) \\[.6ex]
 $v_5^{\prime}$  &  ( 0, 0, 0, 0,-1)
                 &  $u$     &  ( 0, 1, 3) \\[.6ex]
 $b^{\prime}$    &  ( 0, 0, 0, 2, 3)
                 &  $d_1$  &   ( 1, 0, 4) \\[.6ex]
 $e_1^{\prime}$  &  ( 0, 0, 1, 2, 3)  & & \\[.6ex]
 $e_2^{\prime}$  &  ( 0, 0, 2, 2, 3)  & & \\[.6ex]
 $e_3^{\prime}$  &  ( 0, 0, 1, 1, 1)  & & \\[.6ex]
 $f^{\prime}$    &  ( 0, 1, 2, 2, 3)  & & \\[.6ex]
 $g^{\prime}$    &  ( 0, 1, 3, 2, 3)  & & \\[.6ex]
 $v_6^{\prime}$  &  ( 1, 0, 4, 2, 3)  & & \\[.6ex]  \hline
\end{tabular}
} 

\bigskip

\noindent
Denote a cone e.g.\
$[\,v_1^{\prime},\,b^{\prime},\, e_1^{\prime},\, v_2^{\prime},\,
                                                v_4^{\prime}\,]$,
generated $v_1^{\prime}$, $b^{\prime}$, $e_1^{\prime}$, $v_2^{\prime}$,
and $v_4^{\prime}$, by
$v_1^{\prime}b^{\prime}e_1^{\prime}v_2^{\prime}v_4^{\prime}$
for brevity of notation.
Then the $5$-dimensional fan $\Sigma^{\prime}$ consists of $54$
maximal (simplicial) cones$\,:$

{\scriptsize
$$
 \left\{\,
  \begin{array}{cccccc}
   v_1^{\prime}b^{\prime}e_1^{\prime}v_2^{\prime}v_4^{\prime}\;,
     & v_1^{\prime}b^{\prime}f^{\prime}c_1^{\prime}v_4^{\prime}\;,
     & v_1^{\prime}b^{\prime}v_2^{\prime}c_1^{\prime}v_4^{\prime}\;,
     & v_1^{\prime}e_3^{\prime}e_1^{\prime}v_2^{\prime}v_5^{\prime}\;,
     & v_1^{\prime}b^{\prime}e_1^{\prime}v_2^{\prime}v_5^{\prime}\;,
     & v_1^{\prime}b^{\prime}f^{\prime}c_1^{\prime}v_5^{\prime}\;,\\[.6ex] 
   v_1^{\prime}b^{\prime}v_2^{\prime}c_1^{\prime}v_5^{\prime}\;,
     & v_1^{\prime}e_3^{\prime}v_2^{\prime}v_4^{\prime}v_5^{\prime}\;,
     & v_1^{\prime}b^{\prime}f^{\prime}v_4^{\prime}g^{\prime}\;,
     & v_1^{\prime}b^{\prime}e_1^{\prime}v_4^{\prime}g^{\prime}\;,
     & v_1^{\prime}b^{\prime}f^{\prime}v_5^{\prime}g^{\prime}\;,
     & v_1^{\prime}e_3^{\prime}e_1^{\prime}v_5^{\prime}g^{\prime}\;,\\[.6ex]
   v_1^{\prime}b^{\prime}e_1^{\prime}v_5^{\prime}g^{\prime}\;,
     & v_1^{\prime}e_3^{\prime}v_4^{\prime}v_5^{\prime}g^{\prime}\;,
     & v_1^{\prime}f^{\prime}v_4^{\prime}v_5^{\prime}g^{\prime}\;,
     & v_1^{\prime}e_3^{\prime}e_1^{\prime}v_2^{\prime}e_2^{\prime}\;,
     & v_1^{\prime}e_3^{\prime}v_2^{\prime}v_4^{\prime}e_2^{\prime}\;,
     & v_1^{\prime}e_1^{\prime}v_2^{\prime}v_4^{\prime}e_2^{\prime}\;,
                                                                 \\[.6ex]
   v_1^{\prime}e_3^{\prime}e_1^{\prime}g^{\prime}e_2^{\prime}\;,
     & v_1^{\prime}e_3^{\prime}v_4^{\prime}g^{\prime}e_2^{\prime}\;,
     & v_1^{\prime}e_1^{\prime}v_4^{\prime}g^{\prime}e_2^{\prime}\;,
     & v_1^{\prime}f^{\prime}c_1^{\prime}v_4^{\prime}c_2^{\prime}\;,
     & v_1^{\prime}v_2^{\prime}c_1^{\prime}v_4^{\prime}c_2^{\prime}\;,
     & v_1^{\prime}f^{\prime}c_1^{\prime}v_5^{\prime}c_2^{\prime}\;,\\[.6ex]
   v_1^{\prime}v_2^{\prime}c_1^{\prime}v_5^{\prime}c_2^{\prime}\;,
     & v_1^{\prime}f^{\prime}v_4^{\prime}v_5^{\prime}c_2^{\prime}\;,
     & v_1^{\prime}v_2^{\prime}v_4^{\prime}v_5^{\prime}c_2^{\prime}\;,
     & v_6^{\prime}b^{\prime}e_1^{\prime}v_2^{\prime}v_4^{\prime}\;,
     & v_6^{\prime}b^{\prime}f^{\prime}c_1^{\prime}v_4^{\prime}\;,
     & v_6^{\prime}b^{\prime}v_2^{\prime}c_1^{\prime}v_4^{\prime}\;,\\[.6ex]
   v_6^{\prime}e_3^{\prime}e_1^{\prime}v_2^{\prime}v_5^{\prime}\;,
     & v_6^{\prime}b^{\prime}e_1^{\prime}v_2^{\prime}v_5^{\prime}\;,
     & v_6^{\prime}b^{\prime}f^{\prime}c_1^{\prime}v_5^{\prime}\;,
     & v_6^{\prime}b^{\prime}v_2^{\prime}c_1^{\prime}v_5^{\prime}\;,
     & v_6^{\prime}e_3^{\prime}v_2^{\prime}v_4^{\prime}v_5^{\prime}\;,
     & v_6^{\prime}b^{\prime}f^{\prime}v_4^{\prime}g^{\prime}\;, \\[.6ex]
   v_6^{\prime}b^{\prime}e_1^{\prime}v_4^{\prime}g^{\prime}\;,
     & v_6^{\prime}b^{\prime}f^{\prime}v_5^{\prime}g^{\prime}\;,
     & v_6^{\prime}e_3^{\prime}e_1^{\prime}v_5^{\prime}g^{\prime}\;,
     & v_6^{\prime}b^{\prime}e_1^{\prime}v_5^{\prime}g^{\prime}\;,
     & v_6^{\prime}e_3^{\prime}v_4^{\prime}v_5^{\prime}g^{\prime}\;,
     & v_6^{\prime}f^{\prime}v_4^{\prime}v_5^{\prime}g^{\prime}\;, \\[.6ex]
   v_6^{\prime}e_3^{\prime}e_1^{\prime}v_2^{\prime}e_2^{\prime}\;,
     & v_6^{\prime}e_3^{\prime}v_2^{\prime}v_4^{\prime}e_2^{\prime}\;,
     & v_6^{\prime}e_1^{\prime}v_2^{\prime}v_4^{\prime}e_2^{\prime}\;,
     & v_6^{\prime}e_3^{\prime}e_1^{\prime}g^{\prime}e_2^{\prime}\;,
     & v_6^{\prime}e_3^{\prime}v_4^{\prime}g^{\prime}e_2^{\prime}\;,
     & v_6^{\prime}e_1^{\prime}v_4^{\prime}g^{\prime}e_2^{\prime}\;,\\[.6ex]
   v_6^{\prime}f^{\prime}c_1^{\prime}v_4^{\prime}c_2^{\prime}\;,
     & v_6^{\prime}v_2^{\prime}c_1^{\prime}v_4^{\prime}c_2^{\prime}\;,
     & v_6^{\prime}f^{\prime}c_1^{\prime}v_5^{\prime}c_2^{\prime}\;,
     & v_6^{\prime}v_2^{\prime}c_1^{\prime}v_5^{\prime}c_2^{\prime}\;,
     & v_6^{\prime}f^{\prime}v_4^{\prime}v_5^{\prime}c_2^{\prime}\;,
     & v_6^{\prime}v_2^{\prime}v_4^{\prime}v_5^{\prime}c_2^{\prime}\;
  \end{array}\, \right\}
$$
} 

\noindent 
while the $3$-dimensional fan $\Sigma$ consists of $10$ maximal
(simplicial) cones$\,$:
$$
 \left\{\;
  d_4d_2r_2\,,\; d_4d_2u\,,\; d_4ur_1\,,\; d_4d_3r_1\,,\; d_4d_3r_2\,,\;
  d_1d_2r_2\,,\; d_1d_2u\,,\; d_1ur_1\,,\; d_1d_3r_1\,,\; d_1d_3r_2 \;
 \right\}\,.
$$
Following [B-C-dlO-G] and denoting a fan as the summation of its
maximal cones, $\Sigma^{\prime}$ can be decomposed into a combination
of suspensions or extensions by extra rays$\,$:
$$
 \Sigma^{\prime}\;=\;(v_1^{\prime}+v_6^{\prime})\Sigma^{\prime\prime}\;
 \hspace{2em}\hbox{and}\hspace{2em}
 \Sigma^{\prime\prime}\;
   =\;(f^{\prime}+v_2^{\prime})\Sigma^{\prime\prime(c)}
      + (g^{\prime}+v_2^{\prime})\Sigma^{\prime\prime(e)}
      + f^{\prime}g^{\prime} \Sigma^{\prime\prime(0)}\,,
$$
where
\begin{eqnarray*}
 \Sigma^{\prime\prime(c)} & =
   & b^{\prime}c_1^{\prime}v_4^{\prime}
     + c_1^{\prime}c_2^{\prime}v_4^{\prime}
     + b^{\prime}c_1^{\prime}v_5^{\prime}
     + c_1^{\prime}c_2^{\prime}v_5^{\prime}
     + c_2^{\prime}v_4^{\prime}v_5^{\prime}\,,   \\[.6ex]
 \Sigma^{\prime\prime(e)} & =
   & e_1^{\prime}e_2^{\prime}e_3^{\prime}
     + b^{\prime}e_1^{\prime}v_4^{\prime}
     + e_1^{\prime}e_2^{\prime}v_4^{\prime}
     + e_2^{\prime}e_3^{\prime}v_4^{\prime}
     + b^{\prime}e_1^{\prime}v_5^{\prime}
     + e_1^{\prime}e_3^{\prime}v_5^{\prime}
     + e_3^{\prime}v_4^{\prime}v_5^{\prime}\,,  \\[.6ex]
 \Sigma^{\prime\prime(0)} & =
   & b^{\prime}v_4^{\prime} + b^{\prime}v_5^{\prime}
     + v_4^{\prime}v_5^{\prime}\,.  
\end{eqnarray*}
Observe that
$\Sigma^{\prime\prime(0)}
  =\Sigma^{\prime\prime(c)}\cap\Sigma^{\prime\prime(e)}$.

Similarly, $\Sigma$ is the suspension by $\{d_4, d_1\}$ of the
$2$-dimensional subfan generated by $\Sigma(1)-\{d_4, d_1\}$ and
can be written as
$$
 \Sigma\;=\; (d_4+d_1)\Sigma^{\sim}\,,
 \hspace{1em}\hbox{where}\hspace{1em}
 \Sigma^{\sim}\;=\;(d_2+d_3)\,r_2 \,+\, (u+d_3)\,r_1 + d_2u\,.
$$
Note that all the cones in $\Sigma$ have multiplicity $1$ and hence
$X_{\Sigma}$ is a smooth toric $3$-fold. It is the blow-up of the
Hirzebruch $3$-fold ${\Bbb F}_{024}$ at a point,
cf.\ Table 4.1 in [Ra].

\bigskip

\noindent 
{\it (b) The toric morphism
  $\widetilde{\varphi}:X_{\Sigma^{\prime}}\rightarrow X_{\Sigma}$.}
The map of fans $\varphi:\Sigma^{\prime}\rightarrow\Sigma$
is induced by the projection map
$\varphi:(y_1, y_2, y_3, y_4, y_5)\mapsto (y_1, y_2, y_3)$.
Under this map,
$$
 \Sigma^{\prime\prime(0)}\;\rightarrow\; 0\,, \hspace{2em}
 \Sigma^{\prime\prime(c)}-\Sigma^{\prime\prime(0)}\;
   \rightarrow\; r_2\,,        \hspace{2em}
 \Sigma^{\prime\prime(e)}-\Sigma^{\prime\prime(0)}\;
   \rightarrow\; r_1\,.                           
$$
Together with the correspondence
$$
 v_1^{\prime}\; \rightarrow\; d_4\,, \hspace{2em}
 v_6^{\prime}\; \rightarrow\; d_1\,, \hspace{2em}
 f^{\prime}  \; \rightarrow\; d_2\,, \hspace{2em}
 g^{\prime}  \; \rightarrow\; u  \,, \hspace{2em}
 v_2^{\prime}\; \rightarrow\; d_3\,,
$$
the set $\Sigma^{\prime}_{\sigma}$ of cones in $\Sigma^{\prime}$
associated to $\sigma\in\Sigma$ under $\varphi$ and the set
$\Sigma^{\prime\circ}_{\sigma}$ of primitive cones therein can be
directly worked out. Furthermore, from Remark 2.1.12,
$\Ind(\sigma)=1$ for all $\sigma\in\Sigma$ since $\varphi$ is
surjective.
All these are listed in Columns 2 - 4 of Table 2-2-2 below.

\bigskip

\begin{minipage}{12cm}
\centerline{\scriptsize 
 \begin{tabular}{|c|l|l|c|l|} \hline
 \rule{0ex}{3ex}
  $\sigma\in\Sigma$
  & \multicolumn{1}{c|}{  $\Sigma^{\prime}_{\sigma}$  }
  & \multicolumn{1}{c|}{ $\tau^{\prime}\in\Sigma^{\prime\circ}_{\sigma}$ }
  & $\Ind(\sigma)$
  & \multicolumn{1}{c|}{
      the irreducible fiber component $F^{\tau^{\prime}}_{\sigma}$ }
         \\[.6ex] \hline\hline
 \rule{0ex}{3ex}
 $0$
  &  $\varphi^{-1}(0)=\Sigma^{\prime\prime(0)}$  
  &  $0$    &  $1$
  &  $\scriptsizeWCP^2(1,2,3)$       \\[.6ex] \hline
 \rule{0ex}{3ex}
  $d_4$
   & $v_1^{\prime}\,\Sigma^{\prime\prime(0)}$
   & $v_1^{\prime}$     &  $1$
   & $\scriptsizeWCP^2(1,2,3)$       \\
  \multicolumn{1}{|c|}{ -------- }
     & \multicolumn{1}{c|}{ ---------------------------------  }
     & \multicolumn{1}{c|}{ --------------------------------  }
     &  \multicolumn{1}{c|}{ --------  }
     & \multicolumn{1}{c|}{ ------------------------------------------  }
      \\
  $d_3$
   & $v_2^{\prime}\,\Sigma^{\prime\prime(0)}$
   & $v_2^{\prime}$     &  $1$
   & $\scriptsizeWCP^2(1,2,3)$       \\[.6ex]
  $r_2$
   & $\Sigma^{\prime\prime(c)}-\Sigma^{\prime\prime(0)}$
   & $c_1^{\prime}$, $\;c_2^{\prime}$      &  $1$
   & $X(4)$, \quad $\;{\scriptsizeBbb C}{\rm P}^2$    \\ [.6ex]
  $r_1$
   & $\Sigma^{\prime\prime(e)}-\Sigma^{\prime\prime(0)}$
   & $e_1^{\prime}$, $\;e_2^{\prime}$, $\;e_3^{\prime}$  & $1$
   & $X(5)$, \quad $\scriptsizeWCP^2(1,1,3)$, \quad
      ${\scriptsizeBbb F}_2$        \\[.6ex]
  $d_2$
   & $f^{\prime}\,\Sigma^{\prime\prime(0)}$
   & $f^{\prime}$       &   $1$
   & $\scriptsizeWCP^2(1,2,3)$      \\[.6ex]
  $u$
   & $g^{\prime}\,\Sigma^{\prime\prime(0)}$
   & $g^{\prime}$       &   $1$
   & $\scriptsizeWCP^2(1,2,3)$      \\
  \multicolumn{1}{|c|}{ -------- }
     & \multicolumn{1}{c|}{ ---------------------------------  }
     & \multicolumn{1}{c|}{ --------------------------------  }
     &  \multicolumn{1}{c|}{ --------  }
     & \multicolumn{1}{c|}{ ------------------------------------------  }
       \\
  $d_1$
   & $v_6^{\prime}\,\Sigma^{\prime\prime(0)}$
   & $v_6^{\prime}$     &     $1$
   & $\scriptsizeWCP^2(1,2,3)$     \\[.6ex] \hline
 \rule{0ex}{3ex}
  $d_2r_2$
   & $f^{\prime}\,(\Sigma^{\prime\prime(c)}-\Sigma^{\prime\prime(0)})$
   & $f^{\prime}c_1^{\prime}$, $\;f^{\prime}c_2^{\prime}$ & $1$
   & $X(4)$, \quad $\;{\scriptsizeBbb C}{\rm P}^2$      \\[.6ex]
  $d_2u$
   &  $f^{\prime}g^{\prime}\Sigma^{\prime\prime(0)}$
   &  $f^{\prime}g^{\prime}$       &   $1$
   &  $\scriptsizeWCP^2(1,2,3)$     \\[.6ex]
  $ur_1$
   &  $g^{\prime}\,(\Sigma^{\prime\prime(e)}-\Sigma^{\prime\prime(0)})$
   &  $g^{\prime}e_1^{\prime}$, $\;g^{\prime}e_2^{\prime}$,
       $\;g^{\prime}e_3^{\prime}$   &     $1$
   & $X(5)$, \quad $\scriptsizeWCP^2(1,1,3)$, \quad
      ${\scriptsizeBbb F}_2$         \\[.6ex]
  $d_3r_1$
   &  $v_2^{\prime}\,(\Sigma^{\prime\prime(e)}-\Sigma^{\prime\prime(0)})$
   &  $v_2^{\prime}e_1^{\prime}$, $\;v_2^{\prime}e_2^{\prime}$,
       $\;v_2^{\prime}e_3^{\prime}$   &   $1$
   & $X(5)$, \quad $\scriptsizeWCP^2(1,1,3)$, \quad
      ${\scriptsizeBbb F}_2$         \\[.6ex]
  $d_3r_2$
   &  $v_2^{\prime}\,(\Sigma^{\prime\prime(c)}-\Sigma^{\prime\prime(0)})$
   &  $v_2^{\prime}c_1^{\prime}$,  $\;v_2^{\prime}c_2^{\prime}$  & $1$
   &  $X(4)$, \quad $\;{\scriptsizeBbb C}{\rm P}^2$     \\
  \multicolumn{1}{|c|}{ -------- }
     & \multicolumn{1}{c|}{ ---------------------------------  }
     & \multicolumn{1}{c|}{ --------------------------------  }
     &  \multicolumn{1}{c|}{ --------  }
     & \multicolumn{1}{c|}{ ------------------------------------------  }
       \\
  $d_4d_3$
   &  $v_1^{\prime}v_2^{\prime}\,\Sigma^{\prime\prime(0)}$
   &  $v_1^{\prime}v_2^{\prime}$     &  $1$
   &  $\scriptsizeWCP^2(1,2,3)$      \\[.6ex]
  $d_4r_2$
   &  $v_1^{\prime}\,(\Sigma^{\prime\prime(c)}-\Sigma^{\prime\prime(0)})$
   &  $v_1^{\prime}c_1^{\prime}$, $\;v_1^{\prime}c_2^{\prime}$
   &  $1$
   & $X(4)$, \quad $\;{\scriptsizeBbb C}{\rm P}^2$     \\[.6ex]
  $d_4r_1$
   &  $v_1^{\prime}(\Sigma^{\prime\prime(e)}-\Sigma^{\prime\prime(0)})$
   &  $v_1^{\prime}e_1^{\prime}$, $\;v_1^{\prime}e_2^{\prime}$,
       $\;v_1^{\prime}e_3^{\prime}$     &  $1$
   & $X(5)$, \quad $\scriptsizeWCP^2(1,1,3)$, \quad
      ${\scriptsizeBbb F}_2$          \\[.6ex]
  $d_4d_2$  
   &  $v_1^{\prime}f^{\prime}\,\Sigma^{\prime\prime(0)}$
   &  $v_1^{\prime}f^{\prime}$         &   $1$
   &  $\scriptsizeWCP^2(1,2,3)$        \\[.6ex]
  $d_4u$   
   &  $v_1^{\prime}g^{\prime}\,\Sigma^{\prime\prime(0)}$
   &  $v_1^{\prime}g^{\prime}$      &   $1$
   &  $\scriptsizeWCP^2(1,2,3)$       \\
  \multicolumn{1}{|c|}{ -------- }
     & \multicolumn{1}{c|}{ ---------------------------------  }
     & \multicolumn{1}{c|}{ --------------------------------  }
     &  \multicolumn{1}{c|}{ --------  }
     & \multicolumn{1}{c|}{ ------------------------------------------  }
       \\
  $d_1d_3$  
   &  $v_6^{\prime}v_2^{\prime}\,\Sigma^{\prime\prime(0)}$
   &  $v_6^{\prime}v_2^{\prime}$    &  $1$
   &  $\scriptsizeWCP^2(1,2,3)$     \\[.6ex]
  $d_1r_2$  
   &  $v_6^{\prime}\,(\Sigma^{\prime\prime(c)}-\Sigma^{\prime\prime(0)})$
   &  $v_6^{\prime}c_1^{\prime}$, $v_6^{\prime}\;c_2^{\prime}$
   &   $1$
   & $X(4)$, \quad $\;{\scriptsizeBbb C}{\rm P}^2$      \\[.6ex]
  $d_1r_1$  
   &  $v_6^{\prime}\,(\Sigma^{\prime\prime(e)}-\Sigma^{\prime\prime(0)})$
   &  $v_6^{\prime}e_1^{\prime}$, $\;v_6^{\prime}e_2^{\prime}$,
       $\;v_6^{\prime}e_3^{\prime}$   &  $1$
   & $X(5)$, \quad $\scriptsizeWCP^2(1,1,3)$, \quad
      ${\scriptsizeBbb F}_2$         \\[.6ex]
  $d_1d_2$ 
   &  $v_6^{\prime}f^{\prime}\,\Sigma^{\prime\prime(0)}$
   &  $v_6^{\prime}f^{\prime}$       &  $1$
   &  $\scriptsizeWCP^2(1,2,3)$      \\[.6ex]
  $d_1u$   
   &  $v_6^{\prime}g^{\prime}\,\Sigma^{\prime\prime(0)}$
   &  $v_6^{\prime}g^{\prime}$     &   $1$
   &  $\scriptsizeWCP^2(1,2,3)$      \\[.6ex]  \hline
 \rule{0ex}{3ex}
  $d_4d_2r_2$  
   &  $v_1^{\prime}f^{\prime}\,
                  (\Sigma^{\prime\prime(c)}-\Sigma^{\prime\prime(0)})$
   &  $v_1^{\prime}f^{\prime}c_1^{\prime}$,
       $\;v_1^{\prime}f^{\prime}c_2^{\prime}$   &  $1$
   & $X(4)$, \quad $\;{\scriptsizeBbb C}{\rm P}^2$     \\[.6ex]
  $d_4d_2u$    
   &  $v_1^{\prime}f^{\prime}g^{\prime}\Sigma^{\prime\prime(0)}$
   &  $v_1^{\prime}f^{\prime}g^{\prime}$    &   $1$
   &  $\scriptsizeWCP^2(1,2,3)$     \\[.6ex]
  $d_4ur_1$   
   &  $v_1^{\prime}g^{\prime}\,
             (\Sigma^{\prime\prime(e)}-\Sigma^{\prime\prime(0)})$
   &  $v_1^{\prime}g^{\prime}e_1^{\prime}$,
       $\;v_1^{\prime}g^{\prime}e_2^{\prime}$,
       $\;v_1^{\prime}g^{\prime}e_3^{\prime}$    &   $1$
   & $X(5)$, \quad $\scriptsizeWCP^2(1,1,3)$, \quad
      ${\scriptsizeBbb F}_2$         \\[.6ex]
  $d_4d_3r_1$ 
   &  $v_1^{\prime}v_2^{\prime}\,
              (\Sigma^{\prime\prime(e)}-\Sigma^{\prime\prime(0)})$
   &  $v_1^{\prime}v_2^{\prime}e_1^{\prime}$,
       $\;v_1^{\prime}v_2^{\prime}e_2^{\prime}$,
       $\;v_1^{\prime}v_2^{\prime}e_3^{\prime}$    &  $1$
   & $X(5)$, \quad $\scriptsizeWCP^2(1,1,3)$, \quad
      ${\scriptsizeBbb F}_2$        \\[.6ex]
  $d_4d_3r_2$  
   &  $v_1^{\prime}v_2^{\prime}\,
               (\Sigma^{\prime\prime(c)}-\Sigma^{\prime\prime(0)})$
   &  $v_1^{\prime}v_2^{\prime}c_1^{\prime}$,
       $\;v_1^{\prime}v_2^{\prime}c_2^{\prime}$    &   $1$
   &  $X(4)$, \quad $\;{\scriptsizeBbb C}{\rm P}^2$       \\
  \multicolumn{1}{|c|}{ -------- }
     & \multicolumn{1}{c|}{ ---------------------------------  }
     & \multicolumn{1}{c|}{ --------------------------------  }
     &  \multicolumn{1}{c|}{ --------  }
     & \multicolumn{1}{c|}{ ------------------------------------------  }
       \\
  $d_1d_2r_2$ 
   &  $v_6^{\prime}f^{\prime}\,
                 (\Sigma^{\prime\prime(c)}-\Sigma^{\prime\prime(0)})$
   &  $v_6^{\prime}f^{\prime}c_1^{\prime}$,
       $\;v_6^{\prime}f^{\prime}c_2^{\prime}$      &  $1$
   &  $X(4)$, \quad $\;{\scriptsizeBbb C}{\rm P}^2$     \\[.6ex]
  $d_1d_2u$   
   &  $v_6^{\prime}f^{\prime}g^{\prime}\Sigma^{\prime\prime(0)}$
   &  $v_6^{\prime}f^{\prime}g^{\prime}$    &    $1$
   &  $\scriptsizeWCP^2(1,2,3)$       \\[.6ex]
  $d_1ur_1$   
   &  $v_6^{\prime}g^{\prime}\,
            (\Sigma^{\prime\prime(e)}-\Sigma^{\prime\prime(0)})$
   &  $v_6^{\prime}g^{\prime}e_1^{\prime}$,
       $\;v_6^{\prime}g^{\prime}e_2^{\prime}$,
       $\;v_6^{\prime}g^{\prime}e_3^{\prime}$      &  $1$
   & $X(5)$, \quad $\scriptsizeWCP^2(1,1,3)$, \quad
      ${\scriptsizeBbb F}_2$         \\[.6ex]
  $d_1d_3r_1$ 
   &  $v_6^{\prime}v_2^{\prime}\,
           (\Sigma^{\prime\prime(e)}-\Sigma^{\prime\prime(0)})$
   &  $v_6^{\prime}v_2^{\prime}e_1^{\prime}$,
       $\;v_6^{\prime}v_2^{\prime}e_2^{\prime}$,
       $\;v_6^{\prime}v_2^{\prime}e_3^{\prime}$    &  $1$
   & $X(5)$, \quad $\scriptsizeWCP^2(1,1,3)$, \quad
      ${\scriptsizeBbb F}_2$         \\[.6ex]
  $d_1d_3r_2$  
   &  $v_6^{\prime}v_2^{\prime}\,
            (\Sigma^{\prime\prime(c)}-\Sigma^{\prime\prime(0)})$
   &  $v_6^{\prime}v_2^{\prime}c_1^{\prime}$,
       $\;v_6^{\prime}v_2^{\prime}c_2^{\prime}$    &  $1$
   &  $X(4)$, \quad $\;{\scriptsizeBbb C}{\rm P}^2$   \\[.6ex] \hline
\end{tabular}
} 

 \bigskip

 \centerline{
   \parbox{9cm}{ {\sc Table 2-2-2.}
   A toric description of the toric morphism $\widetilde{\varphi}$.
   The irreducible components $F^{\tau^{\prime}}_{\sigma}$ over the
   toric orbit $O_{\sigma}$ is listed in the same order as 
   $\tau^{\prime}$.
   The first four columns are explained in Part (a) and Part (b);
   the last column is explained in Part (c), in which the notations
   $X(4)=X_{\{(2,3),(-1,0),(-1,-1),(0,-1)\}}$ and
   $X(5)=X_{\{(2,3),(-1,0),(-2,-3),(-1,-2),(0,-1)\}}$ are used. }
 } 
\end{minipage}

\bigskip

\bigskip

\noindent
{\it (c) The fibers of $\widetilde{\varphi}$.}
From the table above, all the fibers of $\widetilde{\varphi}$ are
connected. Among the $33$ torus orbits of $X_{\Sigma}$, the fiber
is irreducible only over $15$ of them.
There are $59$ primitive cones $\tau^{\prime}$ in $\Sigma^{\prime}$
with respect to the map of fans $\varphi$. To see their relative star
$\Star_{\sigma}(\tau^{\prime})$, observe that
\begin{eqnarray*}
 \Sigma^{\prime\prime(c)}
  & = & c_1^{\prime}\,(b^{\prime}v_4^{\prime}+c_2^{\prime}v_4^{\prime}
                     + b^{\prime}v_5^{\prime}+c_2^{\prime}v_5^{\prime})\,
        +\, (\,\mbox{cones without $c_1^{\prime}$}\,)       \\[.6ex]
  & = & c_2^{\prime}\,(c_1^{\prime}v_4^{\prime}+c_1^{\prime}v_5^{\prime}
                         + v_4^{\prime}v_5^{\prime})\,
        +\,(\,\mbox{cones without $c_2^{\prime}$}\,) \\[.6ex]
  & = & c_1^{\prime}c_2^{\prime}\,(v_4^{\prime}+v_5^{\prime})\,
        +\,(\,\mbox{cones without $c_1^{\prime}c_2^{\prime}$}\,)\,.
\end{eqnarray*}
These identities together with that fact that all the cones explicitly
shown are in $\Sigma^{\prime}_{r_2}$ give the combinatorial type of
the fans for fibers $F_{r_2}$, $F_{d_2r_2}$, $F_{d_3r_2}$, $F_{d_4r_2}$,
$F_{d_1r_2}$, $F_{d_4d_2r_2}$, $F_{d_4d_3r_2}$, $F_{d_1d_2r_2}$,
$F_{d_1d_3r_2}$. They imply that each of these fibers has two
irreducible components, which intersect at a toric-invariant $\CP^1$.
Observe also that
\begin{eqnarray*}
 \Sigma^{\prime\prime(e)}
  & = & e_1^{\prime}\,(e_2^{\prime}e_3^{\prime} + b^{\prime}v_4^{\prime}
                       + e_2^{\prime}v_4^{\prime} + b^{\prime}v_5^{\prime}
                       + e_3^{\prime}v_5^{\prime} )\,
        +\, (\,\mbox{cones without $e_1^{\prime}$ }\,)   \\[.6ex]
  & = & e_2^{\prime}\,(e_1^{\prime}e_3^{\prime}
               + e_1^{\prime}v_4^{\prime} + e_3^{\prime}v_4^{\prime})\,
        +\, (\,\mbox{cones without $e_2^{\prime}$}\,)   \\[.6ex]
  & = & e_3^{\prime}\,(e_1^{\prime}e_2^{\prime} + e_2^{\prime}v_4^{\prime}
               + e_1^{\prime}v_5^{\prime} + v_4^{\prime}v_5^{\prime})\,
        +\, (\,\mbox{cones without $e_3^{\prime}$}\,)   \\[.6ex]
  & = & e_1^{\prime}e_2^{\prime}\,(e_3^{\prime} + v_4^{\prime})\,
        +\, (\,\mbox{cones without $e_1^{\prime}e_2^{\prime}$}\,) \\[.6ex]
  & = & e_2^{\prime}e_3^{\prime}\,(e_1^{\prime} + v_4^{\prime})\,
        +\, (\,\mbox{cones without $e_2^{\prime}e_3^{\prime}$}\,) \\[.6ex]
  & = & e_1^{\prime}e_3^{\prime}\,(e_2^{\prime} + v_5^{\prime})\,
        +\, (\,\mbox{cones without $e_1^{\prime}e_3^{\prime}$}\,) \\[.6ex]
  & = & e_1^{\prime}e_2^{\prime}e_3^{\prime}\,
        +\, (\,\mbox{cones without
                             $e_1^{\prime}e_2^{\prime}e_3^{\prime}$}\,)
\end{eqnarray*}
These identities together with the fact that all the cones explicitly
shown are in $\Sigma^{\prime}_{r_1}$ give the combinatorial type of
the fans for fibers $F_{r_1}$, $F_{ur_1}$, $F_{d_3r_1}$, $F_{d_4r_1}$,
$F_{d_1r_1}$, $F_{d_4ur_1}$, $F_{d_4d_3r_1}$, $F_{d_1ur_1}$,
$F_{d_1d_3r_1}$. They imply also that each of these fibers has
three irreducible components, any pair of which intersect at
a toric-invariant $\CP^1$ and the three components intersect
at a single point.

To determine $F^{\tau^{\prime}}_{\sigma}$, one needs to know 
the lattice structure in
$\varphi^{-1}((N_{\sigma})_{\scriptsizeBbb R})/
  \mbox{\raisebox{-.4ex}{$(N^{\prime}_{\tau^{\prime}}
                                     )_{\scriptsizeBbb R}$}}$
as well. It turns out that, in the current example, the lattice in 
$\varphi^{-1}((N_{\sigma})_{\scriptsizeBbb R})/
  \mbox{\raisebox{-.4ex}{$(N^{\prime}_{\tau^{\prime}}
                                     )_{\scriptsizeBbb R}$}}$
can be canonically identified with the lattice in
$(\varphi^{-1}(0))_{\scriptsizeBbb R}$ except for $\tau^{\prime}$
that contain $e_2^{\prime}$. In the latter case, the lattice
structure involved can be identified with the lattice
${\Bbb Z}\oplus\frac{1}{2}{\Bbb Z}$ in
$(\varphi^{-1}(0))_{\scriptsizeBbb R}$ when an appropriate basis
for the lattice $\varphi^{-1}(0)$ is chosen. One can also check
that $F^{\tau^{\prime}}_{\sigma}$ depends only on which or none of
$c_1^{\prime}$, $c_2^{\prime}$, $e_1^{\prime}$, $e_2^{\prime}$,
$e_3^{\prime}$ is contained in $\tau^{\prime}$.

All these observations together with some explicit calculations
give us a complete description of all the irreducible components
$F^{\tau^{\prime}}_{\sigma}$ of $F_{\sigma}$ and their intersections
in $F_{\sigma}$, as listed in the last column of the Table 2-2-2
and indicated in {\sc Figure 2-2-3} via fans.
In particular, there are $6$ different kinds of toric varieties
that an irreducible fiber components $F^{\tau^{\prime}}_{\sigma}$
can assume.

\bigskip

\centerline{
\begin{tabular}{|l|cccccc|}  \hline
 \rule{0ex}{3ex}
  ray involved  & $\bullet$         & $c_1^{\prime}$   & $c_2^{\prime}$
                & $e_1^{\prime}$    & $e_2^{\prime}$   & $e_3^{\prime}$
                   \\[.6ex]  \hline
 \rule{0ex}{3ex}
  toric variety & $\WCP^2(1,2,3)$   & $X(4)$           & $\CP^2$
                & $X(5)$            & $\WCP^2(1,1,3)$  & ${\Bbb F}_2$
                   \\[.6ex]  \hline
 \rule{0ex}{3ex}
  reflexivity   & yes               & yes              & yes
                & no                & no               & yes 
                   \\[.6ex]  \hline
\end{tabular}
}  

\bigskip

\noindent
where $\WCP^2(1,2,3)$ and $\WCP^2(1,1,3)$ are weighted projective
spaces, $X(4)$ (resp.\ $X(5)$, ${\Bbb F}_2$) is the toric variety
associated to the $2$-dimensional fan determined by the set of
$1$-cones $\{(2,3),(-1,0),(-1,-1),(0,-1)\}$  \newline
(resp.\ $\{(2,3),(-1,0),(-2,-3),(-1,-2),(0,-1)\}$,
 $\{(1,0),(0,1),(-1,2),(0,-1)\}$).
\begin{figure}[htbp]
\vspace{2in}
  \caption{{\sc Figure 2-2-3.}
   \baselineskip 14pt
   The fans for the generic fiber and the reducible nongeneric fibers
    are indicated. An reducible toric variety can be coded in
    a multi-fan; the intersection relations of irreducible components
    can also be read off from it. 
   The $1$-cones in $\Sigma^{\prime}(1)$ that give rise to the
    corresponding cones for the fiber are used as labels.
   The T-Weil diviser $D_i$ associated to each ray is also indicated
    for later use.
   (For clarity of picture, all the "$^{\prime}$"'s are omitted.)
  } 
 \end{figure}

\noindent\hspace{12cm} $\Box$

\bigskip

\noindent
{\it Remark 2.2.3.}
Note that in Example 2.2.2, each non-degenerate fiber contains
an irreducible component that is associated to a reflexive polytope
(while it may have other components that are not). It is interesting
to know if this is always so for toric morphisms from fans related
to reflexive polytopes. If not, then how to characterize those
fibrations that have this property.

\bigskip

\section{Induced morphism and fibers for hypersurfaces.}

Having understood the fibers of a toric morphism
 $\widetilde{\varphi}:X_{\Sigma^{\prime}}\rightarrow X_{\Sigma}$,
 our next theme is the study of the induced morphism
 $\widetilde{\varphi}|_{Y^{\prime}}:Y^{\prime}\rightarrow X_{\Sigma}$
 for a hypersurface $Y^{\prime}$ in $X_{\Sigma^{\prime}}$.
 Since hypersurfaces in $X_{\Sigma^{\prime}}$ that we will be
 interested in are realized as the zero locus of sections of some
 line bundle over $X_{\Sigma^{\prime}}$, the problem can be converted
 to the study of the various restrictions of a line bundle over
 $X_{\Sigma^{\prime}}$ with a section.
The general study of this is given in this section.
The results presented here will then be applied to the case of
 Calabi-Yau hypersurfaces in Sec.\ 4.

\bigskip

\subsection{Preparations.}

We give first some basic lemmas relating $\widetilde{\varphi}$ and
sections of a line bundle ${\cal L}$ over $X_{\Sigma^{\prime}}$ and
then use the flattening stratification of $\widetilde{\varphi}$
described in Proposition 2.1.4 in Sec.\ 2.1 to study
$\widetilde{\varphi}|_{Y^{\prime}}$ for $Y^{\prime}$ realized as
the zero-locus of a section of ${\cal L}$.

Recall first the following facts$\,$:
(Proposition 2.4 and Corollary 2.5 in [Od2])

\bigskip

\noindent
{\bf Fact 3.1.1 [equivariant line bundle].} {\it
 Let
 $X^{\prime}=X_{\Sigma^{\prime}}$,
 $\SF(\Sigma^{\prime})$ be the group of continuous functions from
  $N^{\prime}_{\scriptsizeBbb R}$ to ${\Bbb R}$ that are linear on
  each $\sigma^{\prime}\in\Sigma^{\prime}$, and
 $\CDiv_{T_{N^{\prime}}}(X^{\prime})$ be the group of
  $T_{N^{\prime}}$-invariant Cartier divisors on $X$.
 For $h^{\prime}\in \SF(\Sigma^{\prime})$, let $D_{h^{\prime}}$ be
 the Cartier divisor on $X^{\prime}$ associated to $h^{\prime}$.
 Then
 \begin{quote}
  \hspace{-1.9em}(1)\hspace{1ex}
  Given a $T_{N^{\prime}}$-invariant Cartier divisor $D^{\prime}$ on
  $X^{\prime}$, then the restriction $D^{\prime}$ to each affine open
  set $U_{\sigma^{\prime}}$, $\sigma^{\prime}\in\Sigma^{\prime}$, 
  coincides with a principal divisor on $U_{\sigma^{\prime}}$
  so that $D^{\prime}=D_{h^{\prime}}$ for some
  $h^{\prime}\in\SF(\Sigma^{\prime})$.
  The map $h^{\prime}\mapsto D_{h^{\prime}}$ gives an isomorphism
  $\SF(\Sigma^{\prime})
                \simeqrightarrow\CDiv_{T_{N^{\prime}}}(X^{\prime})$.

  \hspace{-1.9em}(2)\hspace{1ex}
  For any Cartier divisor $D^{\prime}$ on $X^{\prime}$, one has an
  ${\cal O}_{X^{\prime}}$-module isomorphism
  ${\cal O}_{X^{\prime}}(D^{\prime})
     \simeq{\cal O}_{X^{\prime}}(D_{h^{\prime}})$
  for some $h^{\prime}\in\SF(\Sigma^{\prime})$.
  Thus, the composition of natural homomorphisms
  $\SF(\Sigma^{\prime})\rightarrow\CDiv(X^{\prime})
                                         \rightarrow\Pic(X^{\prime})$
  is surjective.

  \hspace{-1.9em}(3)\hspace{1ex}
  The following are equivalent for $h^{\prime}\in\SF(\Sigma^{\prime})\,$:
  \vspace{-1ex}
  \begin{quote}
  \hspace{-1.9em}(a)\hspace{1ex}
  $h^{\prime}\in M^{\prime}$.

  \hspace{-1.9em}(b)\hspace{1ex}
  $D_{h^{\prime}}$ is a principal divisor.

  \hspace{-1.9em}(c)\hspace{1ex}
  ${\cal O}_{X^{\prime}}(D_{h^{\prime}})\simeq{\cal O}_{X^{\prime}}$
  as ${\cal O}_{X^{\prime}}$-modules.
  \end{quote}

  \hspace{-1.9em}(4)\hspace{1ex}
  As a corollary of (1), (2), and (3), one has the canonical
  isomorphism
  {\footnotesize
   \begin{eqnarray*}
   \SF(\Sigma^{\prime})/
           \mbox{\raisebox{-.4ex}{$M^{\prime}$}}
     & \simeqrightarrow
     & \CDiv_{T_{N^{\prime}}}(X^{\prime})/
           \mbox{\raisebox{-.4ex}{$\CDiv_{T_{N^{\prime}}}(X^{\prime})
                                             \cap\PDiv(X^{\prime})$}} \\
     & \simeqrightarrow
     & \ELB(X^{\prime})/
        \mbox{\raisebox{-.4ex}{$\{{\cal O}_{X^{\prime}}(m^{\prime})\,|\,
                                         m^{\prime}\in M^{\prime}\}$}}
     \simeqrightarrow
   \Pic(X^{\prime})\,,
   \end{eqnarray*}
  {\normalsize where}
  } 
  $\PDiv(X^{\prime})$ is the group of principal divisors on $X^{\prime}$,
  $\ELB(X^{\prime})$ is the group of equivariant line bundles over
   $X^{\prime}$, and
  ${\cal O}_{X^{\prime}}(m^{\prime})$ is the equivariant trivial line
   bundle ${\cal O}_{X^{\prime}}$ with the linearization given by the
   character $\chi^{m^{\prime}}$.
 \end{quote}
} 

\bigskip

\noindent
{\it Remark 3.1.2} ([Fu] and [Od2], also [Ka]). {
An element $h^{\prime}\in\SF(\Sigma^{\prime})$ determines
(non-uniquely) a system of weights
$W_{h^{\prime}}
 =\{\,m^{\prime}_{\sigma^{\prime}}\in M^{\prime}\,|\,
                          \sigma^{\prime}\in\Sigma^{\prime}\,\}$,
from which the following data are specified$\,$:
\begin{quote}
 \hspace{-1.9em}(1)\hspace{1ex}
  The {\it linearization for the trivialization of ${\cal L}_{h^{\prime}}$
  when restricted to the affine charts $U_{\sigma^{\prime}}$} for
  $X_{\Sigma^{\prime}}$. Explicitly, this is given by
  $t^{\prime}\cdot(x^{\prime},a)
     =(\,t^{\prime}\cdot x^{\prime},\,
                       \chi^{m_{\sigma^{\prime}}}(t^{\prime})a\,)$,
  where $t^{\prime}\in T_{N^{\prime}}$ and
  $(x^{\prime}, a)\in U_{\sigma^{\prime}}\times {\Bbb C}$.

 \hspace{-1.9em}(2)\hspace{1ex}
 The {\it transition function} between local trivializations:
 $f_{\sigma^{\prime}_1\sigma^{\prime}_2}:
   U_{\sigma^{\prime}_1}\cap U_{\sigma^{\prime}_2}\rightarrow
                                              {\Bbb C}^{\times}$
 from
 ${\cal L}_{h^{\prime}}|_{
     U_{\sigma^{\prime}_1}\cap U_{\sigma^{\prime}_2}
                                   \subset U_{\sigma^{\prime}_1} }$ 
 to 
 ${\cal L}_{h^{\prime}}|_{
     U_{\sigma^{\prime}_1}\cap U_{\sigma^{\prime}_2}
                                   \subset U_{\sigma^{\prime}_2} }$
 by
 $\chi^{m^{\prime}_{\sigma^{\prime}_2}-m^{\prime}_{\sigma^{\prime}_1}}$.

 \hspace{-1.9em}(3)\hspace{1ex}
 A {\it $T_{N^{\prime}}$-equivariant meromorphic section}, namely
 $\{\chi^{m^{\prime}_{\sigma^{\prime}}}
                           \}_{\sigma^{\prime}\in\Sigma^{\prime}}$,
 in ${\cal L}_{h^{\prime}}$.
\end{quote}
The weights
$\{m^{\prime}_{\sigma^{\prime}}\in M^{\prime}\,|\,
              \sigma^{\prime}\in\Sigma^{\prime}(n^{\prime})\}$
determines an integral polytope $\Delta_{h^{\prime}}$ in
$M^{\prime}_{\scriptsizeBbb R}$ by taking the convex hull and one has
$$
 \begin{array}{cccl}
  \oplus_{m^{\prime}\in\Delta_{h^{\prime}}\cap M^{\prime}}
   {\Bbb C}\cdot\chi^{m^{\prime}}
   & \simeq  & H^0(X_{\Sigma^{\prime}},{\cal L}_{h^{\prime}})
                                                         & \\[.6ex]
   m^{\prime} & \mapsto & s_{m^{\prime}}  &, 
 \end{array}
$$
where 
$s_{m^{\prime}}=\chi^{m^{\prime}-m^{\prime}_{\sigma^{\prime}}}$
over each affine chart $U_{\sigma^{\prime}}$,
$\sigma^{\prime}\in\Sigma^{\prime}$.
The compatibility of elements in $W_{h^{\prime}}$ implies that
$\Sigma^{\prime}$ is a refinement of the normal fan of
$\Delta_{h^{\prime}}$.
} 

\bigskip

\noindent
{\it Remark 3.1.3} [{\it $\,T_{N^{\prime}}$-action on
 $H^0(X_{\Sigma^{\prime}},{\cal O}_{X^{\prime}}(D_{h^{\prime}}))\,$}].{
 The $T_{N^{\prime}}$-action on $X_{\Sigma^{\prime}}$ together with
 the linearization of ${\cal O}_{X^{\prime}}(D_{h^{\prime}})$
 determined by the weight system $W_{h^{\prime}}$ in Remark 3.1.2
 determines a $T_{N^{\prime}}$-action on
 $H^0(X_{\Sigma^{\prime}},{\cal O}_{X^{\prime}}(D_{h^{\prime}}))$.
 Let $t^{\prime}\in T_{N^{\prime}}$; then Remark 3.1.2 implies that 
 $t^{\prime}\cdot s_{m^{\prime}}
              =\chi^{-m^{\prime}}(t^{\prime})\cdot s_{m^{\prime}}$
 over $U_{\sigma^{\prime}}$ for                                 
 $m^{\prime}\in\Delta^{\prime}\cap M^{\prime}$.
} 

\bigskip

\noindent
{\bf Fact 3.1.4 [restriction of line bundle to orbit closure].}
 (Cf.\ [Fu].)
{\it
 Given a line bundle ${\cal L}_{\Delta^{\prime}}$ over
 $X_{\Sigma^{\prime}}$. Let $h^{\prime}\in\SF(\Sigma^{\prime})$
 be the associated continuous piecewise linear function on
 $N^{\prime}_{\scriptsizeBbb R}$ with respect to $\Sigma^{\prime}$.
 Then the restriction of ${\cal L}_{\Delta^{\prime}}$
 to an orbit closure $V(\tau^{\prime})$ can be described as follows$\,$:
 \begin{quote}
  \hspace{-1.9em}(1)\hspace{1ex}
  Let $m^{\prime}_{\tau^{\prime}}\in M^{\prime}$ be a linear function
  on $N^{\prime}_{\scriptsizeBbb R}$ such that
  $m^{\prime}_{\tau^{\prime}}|_{\tau^{\prime}}
                                     =h^{\prime}|_{\tau^{\prime}}$.
  Then
  $h^{\prime\prime}
      =h^{\prime}-m^{\prime}_{\tau^{\prime}}\in\SF(\Sigma^{\prime})$
  defines a line bundle ${\cal L}_{h^{\prime\prime}}$ isomorphic to
  ${\cal L}_{h^{\prime}}$ (with different linearization).
  ${\cal L}_{h^{\prime\prime}}$ corresponds to the polytope
  $\Delta^{\prime}-m^{\prime}_{\tau^{\prime}}$.

  \hspace{-1.9em}(2)\hspace{1ex}
  After this shift, the restriction of $h^{\prime\prime}$ to the
  set of cones $\sigma^{\prime}$ with
  $\tau^{\prime}\prec\sigma^{\prime}$ descends to a continuous
  piecewise linear function on
  $(N^{\prime}/\mbox{\raisebox{-.4ex}{$
               N^{\prime}_{\sigma^{\prime}}$}})_{\scriptsizeBbb R}$
  with respect to $\Star(\tau^{\prime})$.
  This defines a line bundle over $V(\tau^{\prime})$ isomorphic to
  ${\cal L}_{\Delta^{\prime}}|_{V(\tau^{\prime})}$
  (with perhaps different linearization).

  \hspace{-1.9em}(3)\hspace{1ex}
  In terms of polytopes, let $\sigma^{\prime}_i\,$, $i=1,\ldots, k$,
  be the maximal cones in $\Sigma^{\prime}$ with
  $\tau^{\prime}\prec\sigma^{\prime}_i$. Let $m^{\prime}_i$ be the
  vertices of $\Delta^{\prime}$ corresponding to the restriction
  $h^{\prime}|_{\sigma^{\prime}_i}$. Since all $m^{\prime}_i$ give
  the same restriction to $\tau^{\prime}$, all the differences
  $m^{\prime}_i-m^{\prime}_j$ lie in sublattice
  $\tau^{\prime\perp}\cap M^{\prime}$.
  The convex hull $\Delta^{\prime}_{\tau^{\prime}}$ of
  $\{\,m^{\prime}_1,\,\ldots,\,m^{\prime}_k\,\}$ in
  $M^{\prime}_{\scriptsizeBbb R}$ is then contained in, say,
  $m^{\prime}_1+\tau^{\prime\perp}$ and hence determines an integral
  polytope in $\tau^{\prime\perp}$, unique up to a translation by
  a lattice point. This is the polytope in $\tau^{\prime\perp}$
  that gives ${\cal L}_{\Delta^{\prime}}|_{V(\tau^{\prime})}$,
  up to a linearization.
 \end{quote}
} 

\bigskip

The following lemma characterizes $\Delta^{\prime}_{\tau^{\prime}}$
up to a translation.

\bigskip

\noindent
{\bf Lemma 3.1.5.} {\it
 $\Delta^{\prime}_{\tau^{\prime}}
  =(\Delta^{\prime}-m^{\prime}_{\tau^{\prime}})\cap \tau^{\prime\perp}$
 up to an integral translation.
} 

\bigskip

\noindent
{\it Proof.}
After the translation, $\Delta^{\prime}_{\tau^{\prime}}$ becomes the
convex hull of
$\{\,m^{\prime}_1-m^{\prime}_{\tau^{\prime}},\,\ldots,\,
                       m^{\prime}_k-m^{\prime}_{\tau^{\prime}}\,\}$
contained in
$(\Delta^{\prime}-m^{\prime}_{\tau^{\prime}})\cap\tau^{\prime\perp}$.
If there is an extra integral point $m^{\prime}_0$  in
$(\Delta^{\prime}-m^{\prime}_{\tau^{\prime}})\cap\tau^{\prime\perp}$
besides those in $\Delta^{\prime}_{\tau^{\prime}}$, then, since
$m^{\prime}_0\in (\Delta^{\prime}-m^{\prime}_{\tau^{\prime}})$,
it gives rise to a holomorphic section in
${\cal L}={\cal L}_{\Delta^{\prime}-m^{\prime}_{\tau^{\prime}}}$
and hence on ${\cal L}|_{V(\tau^{\prime})}$. On the other hand,
since it is outside $\Delta^{\prime}_{\tau^{\prime}}$, it cannot
be holomorphic on $V(\tau^{\prime})$. This leads to a contradiction
and one concludes the lemma.

\noindent\hspace{12cm} $\Box$

\bigskip

Consider now the restriction homomorphism
$$
 \res\;:\; H^0(X_{\Sigma^{\prime}}, {\cal L}_{\Delta^{\prime}})\,
  \longrightarrow\,H^0(V(\tau^{\prime})\,,
                       {\cal L}_{\Delta^{\prime}_{\tau^{\prime}}})
$$
that sends a holomorphic section of ${\cal L}_{\Delta^{\prime}}$
to its restriction over $V(\sigma^{\prime})$.
Note that for
$m^{\prime}\in(\Delta^{\prime}\cap M^{\prime} )
                                -\Delta^{\prime}_{\tau^{\prime}}$,
the associated holomorphic section in ${\cal L}_{\Delta^{\prime}}$
restricted to $V(\tau^{\prime})$ gives the $0$-section in
${\cal L}_{\Delta^{\prime}_{\tau^{\prime}}}$.
The lemma below follows immediately from Remark 3.1.2$\,$:

\bigskip

\noindent
{\bf Lemma 3.1.6 [sections restricted to orbit closure].} {\it
 The restriction homomorphism
 $\res: H^0(X_{\Sigma^{\prime}}, {\cal L}_{\Delta^{\prime}})\,
   \rightarrow\,H^0(V(\sigma^{\prime})\,,
                       {\cal L}_{\Delta^{\prime}_{\tau^{\prime}}})$
 is surjective with kernel generated by the sections in
 ${\cal L}_{\Delta^{\prime}}$ associated to
 $\Delta^{\prime}\cap M^{\prime}-\Delta^{\prime}_{\tau^{\prime}}$.
 The restriction to $V(\tau^{\prime})$ of the section $s_{m^{\prime}}$
 in ${\cal L}_{\Delta^{\prime}}$ associated to
 $m^{\prime}\in\Delta^{\prime}_{\tau^{\prime}}\cap M^{\prime}
                                              \subset\Delta^{\prime}$
 is the section in ${\cal L}_{\Delta^{\prime}_{\tau^{\prime}}}$
 given by the same $m^{\prime}\in\Delta^{\prime}_{\tau^{\prime}}$.
}

\bigskip

Finally, let us turn to the restriction of a line bundle to a generic
fiber of a toric morphism.
The following lemma characterizes pullback line bundles and sections
from a toric morphism in terms of toric data$\,$:

\bigskip

\noindent
{\bf Lemma 3.1.7 [pullback of line bundle and section].} {\it
 Given lattices $N_1$, $N_2$, and fans $\Sigma_i$ in
 ${N_i}_{\scriptsizeBbb R}$, $i=1,\, 2$.
 Let $\widetilde{f}:X_{\Sigma_1}\rightarrow X_{\Sigma_2}$ be a toric
 morphism associated to a map of fans $f:\Sigma_1\rightarrow \Sigma_2$.
 Let $M_i$ be the dual lattice of $N_i$ and
 $f^{\dagger}:M_2\rightarrow M_1$ be the dual of $f$.
 Suppose that $\Delta_2$ is an integral polytope in
 ${M_2}_{\scriptsizeBbb R}$. Then
 $\widetilde{f}^{\ast}{\cal L}_{\Delta_2}
                                 ={\cal L}_{f^{\dagger}(\Delta_2)}$.
 Furthermore, if $s$ is a section of ${\cal L}_{\Delta_2}$ that
 corresponds to a lattice point $m\in\Delta_2$, then $f^{\ast}s$
 is the section of $f^{\ast}{\cal L}_{\Delta_2}$ that corresponds to
 $f^{\dagger}(m)\in f^{\dagger}(\Delta_2)$.
} 

\bigskip

\noindent
{\it Proof.}
Let $h_2\in \SF(\Sigma_2)$ be the piecewise linear function on
${N_2}_{\scriptsizeBbb R}$ such that
${\cal L}_{h_2}={\cal L}_{\Delta_2}$.
Then the first assertion follows from the isomorphisms
$\widetilde{f}^{\ast}{\cal L}_{\Delta_2}
  = \widetilde{f}^{\ast}{\cal L}_{h_2}
  = {\cal L}_{h_2\circ f} = {\cal L}_{f^{\dagger}(\Delta_2)}$
as equivariant line bundles.

For the second assertion, recall from [Od2] that the piecewise
linear function $h_2$ determines a local linearized trivialization
of ${\cal L}_{\Delta_2}$ over the affine chart $U_{\sigma_2}$ for
each $\sigma_2\in\Sigma_2$ together with a meromorphic section up
to a constant multiple. This meromorphic section over $U_{\sigma_2}$
is determined by the restriction $h_2|_{\sigma_2}$, which can be
represented by an $m_{\sigma_2}\in\Delta_2\cap M_2$. The section
$s$ associated to $m_2$ is then realized by the character
$\chi^{m_2-m_{\sigma_2}}$ over $U_{\sigma_2}$.
Note also that, for $\sigma^{\prime}_1\in\Sigma_1$, let
$m_{\sigma^{\prime}_1}=f^{\dagger}(m_{\sigma^{\prime}_2})$,
where $\sigma^{\prime}_2$ is the unique cone in $\Sigma_2$ such
that the interior of $\sigma^{\prime}_1$ is mapped by $f$ to the
interior of $\sigma^{\prime}_2$, then $f^{\dagger}(\Delta_2)$ is
the convex hull of
$\{m_{\sigma^{\prime}_1}\,|\, \sigma^{\prime}_1\in\Sigma_1\}$.

Now suppose that the interior of $\sigma_1\in\Sigma_1$ is mapped
by $f$ to the interior of $\sigma_2$. Then
$\widetilde{f}(U_{\sigma_1})\subset U_{\sigma_2}$
and $\widetilde{f}^{\ast}(s|_{U_{\sigma_2}})$ is a holomorphic section
in $\widetilde{f}^{\ast}({\cal L}_{\Delta_2}|_{U_{\sigma_2}})$
represented by the character
$\chi^{f^{\dagger}(m_2)-f^{\dagger}(m_{\sigma_2})}$.
Since the local meromorphic section of
$\widetilde{f}^{\ast}({\cal L}_{\Delta_2})|_{U_{\sigma_1}}$
determined by $f^{\dagger}(\Delta_2)$ is given by
$\chi^{f^{\dagger}(m_{\sigma_2})}$ and this discussion holds for
every $\sigma_2\in\Sigma_2$ that contains some cone in $f(\Sigma_1)$,
$f^{\dagger}(m_2)$ must correspond to $\widetilde{f}^{\ast}s$.
This concludes the proof.

\noindent\hspace{12cm} $\Box$

\bigskip

Recall the toric morphism
$\widetilde{\varphi}:X_{\Sigma^{\prime}}\rightarrow X_{\Sigma}$
and the notations from Sec.\ 2.

\bigskip

\noindent
{\bf Corollary 3.1.8 [line bundle and sections restricted to fiber].}
{\it
 Let ${\cal L}_{\Delta^{\prime}}$ be a line bundle over
 $X_{\Sigma^{\prime}}$.
 \begin{quote}
  \hspace{-1.9em}(a)\hspace{1ex}
  The restriction of ${\cal L}_{\Delta^{\prime}}$ to the irreducible
  component $F^{\tau^{\prime}}_{\sigma}$ of a fiber is given by the
  polytope $\overline{\Delta^{\prime}_{\tau^{\prime}}}$ in
  ${\tau^{\prime}}^{\perp}/
      \mbox{\raisebox{-.4ex}{$\varphi^{\dagger}(\sigma^{\perp})$}}$.
  In particular, the (non-linearized) isomorphism class of the
  restriction of ${\cal L}_{\Delta^{\prime}}$ to each fiber over
  the orbit $O_{\sigma}$ in $X_{\Sigma}$ depends only on $\sigma$.

  \hspace{-1.9em}(b)\hspace{1ex}
  The restriction to
  ${\cal L}_{\Delta^{\prime}}|_{F^{\tau^{\prime}}_{\sigma}}$
  of a section in ${\cal L}_{\Delta^{\prime}}$ that
  corresponds to $m^{\prime}\in\Delta_{\tau^{\prime}}\cap M^{\prime}$
  is a constant multiple of the section in
  ${\cal L}_{\Delta^{\prime}}|_{F^{\tau^{\prime}}_{\sigma}}$
  that corresponds to the lattice point
  $\overline{\varphi^{\dagger}}^{\tau^{\prime}}_{\sigma}(m^{\prime})
                              \in \overline{\Delta_{\tau^{\prime}}}$.
 \end{quote}
 {\sc Figure 3-1-1}.
} 

\bigskip

\noindent
{\it Proof.}
Since $F^{\tau^{\prime}}_{\sigma}$ is the generic fiber of the
restriction of $\widetilde{\varphi}$ to the toric variety
$V(\tau^{\prime})$, without loss of generality, we only need to
prove both assertions for the generic fiber of $\widetilde{\varphi}$.
In this case, consider the map of fans
$\varphi^{-1}(0)\rightarrow\Sigma^{\prime}$ induced by the
inclusion map
$\widetilde{\varphi}^{-1}(0)\hookrightarrow
                                    N^{\prime}_{\scriptsizeBbb R}$.
Then both assertions follow from Lemma 3.1.7.

\noindent\hspace{12cm} $\Box$

\begin{figure}[htbp]
 \setcaption{{\sc Figure 3-1-1.}
  \baselineskip 14pt
  Restriction of a line bundle ${\cal L}_{\Delta^{\prime}}$ with 
  a section $s_{m^{\prime}}$ to an irreducible fiber
  $F^{\tau^{\prime}}_{\sigma}$ can be understood in two steps$\,$:
  (1) restriction to $V(\tau^{\prime})$, which is described by
      $\Delta^{\prime}_{\tau^{\prime}}$; and
  (2) further restriction to the generic fiber 
      $F^{\tau^{\prime}}_{\sigma}$ of the toric morphism
      $\widetilde{\varphi}\,:\,V(\tau^{\prime})\rightarrow V(\sigma)$
      from the restriction of
      $\widetilde{\varphi}\,:\,X_{\Sigma^{\prime}}\rightarrow X_{\Sigma}$.
      The latter is described by the quotient
      $\overline{\varphi^{\dagger}}^{\tau^{\prime}}_{\sigma}(m^{\prime})
                        \in \overline{\Delta^{\prime}_{\tau^{\prime}}}$.
  ($m^{\prime}$ and its image are indicated by the thickened darkest
   lattice points.)
 } 
 \centerline{\psfig{figure=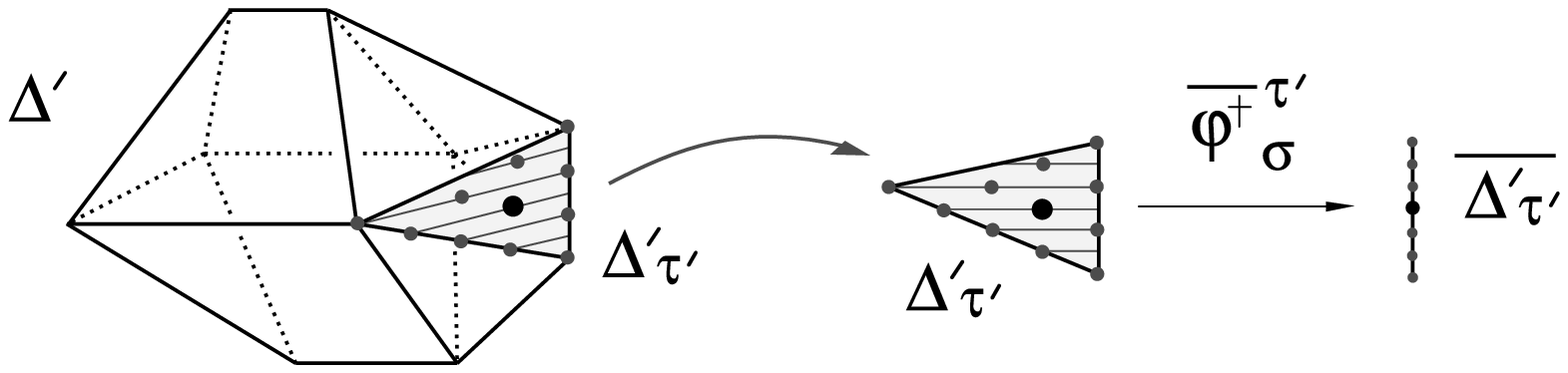,width=11cm,caption=}}
\end{figure}

\bigskip

\noindent
{\it Remark 3.1.9.} {
 When $\Sigma^{\prime}$ is the normal fan of $\Delta^{\prime}$,
 in which case $h^{\prime}$ is strictly convex,
 $\Delta^{\prime}_{\tau^{\prime}}=\Theta_{\tau^{\prime}}$,
 the facet of $\Delta^{\prime}$ dual to $\tau^{\prime}$.
} 

\bigskip

\subsection{The induced morphism$\,$:
    $Y^{\prime}\subset X_{\Sigma^{\prime}}\rightarrow X_{\Sigma}$.}
Let $Y^{\prime}$ be a hypersurface of $X_{\Sigma^{\prime}}$ realized
as the zero-locus of a section of ${\cal L}_{\Delta^{\prime}}$ and 
recall the flattening stratification for the toric morphism
$\widetilde{\varphi}:X_{\Sigma^{\prime}}\rightarrow X_{\Sigma}$,
discussed in Sec.\ 2.1, Proposition 2.1.4.
Our strategy of understanding the induced morphism
$\varphi|_{Y^{\prime}}:Y^{\prime}\rightarrow X_{\Sigma}$
is to understand the restriction of the flattening stratification
of $\widetilde{\varphi}$ to $\widetilde{\varphi}|_{Y^{\prime}}$.
A flattening stratification of $\widetilde{\varphi}|_{Y^{\prime}}$
with each stratum a topological bundle is then obtained by
a refinement of this restriction.
We discuss first a fundamental case and then all other cases,
which can be reduced to the fundamental case.
The notations here follow those in Sec.\ 2.

\bigskip

\begin{flushleft}
{\bf Case (1) [fundamental]$\,$:
     Restriction to generic fibers for $\varphi$ surjective.}
\end{flushleft}
Let $\xi:N\rightarrow N^{\prime}$ be a homomorphism such that
$\varphi\circ\xi=\Id_N$. Then $\xi$ determines a splitting of
$N^{\prime}$ by 
$$
 \begin{array}{ccccc}
  N^{\prime} & =       & \ker\varphi + \im\xi 
             & \simeq  & \ker\varphi\oplus N     \\[.6ex]
   v         & =       & (v-\xi\circ\varphi(v))\,+\,\xi\circ\varphi(v)
             & \mapsto & (v-\xi\circ\varphi(v)\,,\, \varphi(v))\,.
 \end{array}
$$
On the dual side, recall the quotient map
$\overline{\varphi^{\dagger}}:M^{\prime}\rightarrow
        M^{\prime}/\mbox{\raisebox{-.4ex}{$\varphi^{\dagger}(M)$}}$,
then 
$M^{\prime}/\mbox{\raisebox{-.4ex}{$\varphi^{\dagger}(M)$}}$
is canonically isomorphic to $\ker\xi^{\dagger}$.
Thus, there is a map
$i_{\xi}:M^{\prime}/\mbox{\raisebox{-.4ex}{$\varphi^{\dagger}(M)$}}
                                             \rightarrow M^{\prime}$
such that $\xi^{\dagger}\circ i_{\xi}=0$ and that
$\overline{\varphi^{\dagger}}\circ i_{\xi}
  =\Id_{M^{\prime}/
        \mbox{\scriptsize\raisebox{-.4ex}{$\varphi^{\dagger}(M)$}}}$.
From this, one obtains a decomposition:
$$
 \begin{array}{ccccc}
  M^{\prime}  & =   & \ker\xi^{\dagger}+\im\varphi^{\dagger}
              & \simeq
   & M^{\prime}/\mbox{\raisebox{-.4ex}{$\varphi^{\dagger}(M)$}}
                                                   \oplus M \\[.6ex]
  m^{\prime}  & =
   & i_{\xi}\circ\widetilde{\varphi^{\dagger}}(m^{\prime})
     +(m^{\prime}-i_{\xi}\circ\widetilde{\varphi^{\dagger}}(m^{\prime}))
              & \mapsto
              & (\overline{\varphi^{\dagger}}(m^{\prime})\,,\,
                                           \xi^{\dagger}(m^{\prime}))\,.
 \end{array}
$$
Note that different choices of $\xi$ give different projections to
the second compoment of the decomposition
$(\overline{\varphi^{\dagger}},\xi^{\dagger}):
  M^{\prime}\stackrel{\sim}{\rightarrow}
  M^{\prime}/\mbox{\raisebox{-.4ex}{$\varphi^{\dagger}(M)$}}\oplus M$.

Let $\Delta^{\prime}$ be a rational convex polytope in
$M^{\prime}_{\scriptsizeBbb R}$ and ${\cal L}_{\Delta^{\prime}}$
be the associated line bundle over $X_{\Sigma^{\prime}}$.
Note that, since $\xi$ in general does not induce a map of fans
from $\Sigma$ to $\Sigma^{\prime}$, the polytope
$\xi^{\dagger}(\Delta^{\prime})$ in $M_{\scriptsizeBbb R}$ may 
not give a line bundle over $X_{\Sigma}$. However,
$\Sigma_{\xi} \doteq
 \widetilde{\varphi}(\xi(N_{\scriptsizeBbb R})\cap\Sigma^{\prime})$
is a fan in $N_{\scriptsizeBbb R}$ that refines $\Sigma$.
Thus $\xi^{\dagger}(\Delta^{\prime})$ corresponds to the line
bundle $\widetilde{\xi}^{\ast}{\cal L}_{\Delta^{\prime}}$ over
$X_{\Sigma_{\xi}}$. After pushing forward by the morphism
$X_{\Sigma_{\xi}}\rightarrow X_{\Sigma}$, we shall regard 
$\widetilde{\xi}^{\ast}{\cal L}_{\Delta^{\prime}}$ as a rank-one
torsion-free sheaf over $X_{\Sigma}$.

\bigskip

\noindent
{\bf Lemma 3.2.1 [factorization over orbit].} {\it
 Let $O_0=T_N$ be the $T_N$-orbit in $X_{\Sigma}$ corresponding to
 $\{0\}\in\Sigma$ and $F_0$ be the fiber of $\widetilde{\varphi}$
 over $O_0$. Then$\,$:
 \begin{quote}
  \hspace{-1.9em}(a)\hspace{1ex}
  $\xi$ determines a product decomposition
  $\widetilde{\varphi}^{-1}(O_0)=O_0\times F_0$.
  Let $\pr_2$ be the corresponding projection map to the second
  factor, then
  ${\cal L}_{\Delta^{\prime}}|_{\widetilde{\varphi}^{-1}(O_0)}\,
    \simeq\, \pr_2^{\ast}\,{\cal L}_{\overline{\Delta^{\prime}}}$.
 
  \hspace{-1.9em}(b)\hspace{1ex}
  The surjective maps
  $\xi^{\dagger}(\Delta^{\prime})
    \stackrel{\xi^{\dagger}}{\longleftarrow} \Delta^{\prime}
    \stackrel{\overline{\varphi^{\dagger}}}{\longrightarrow}
                                   \overline{\Delta^{\prime}}$
  induce a factorization
  $$
   \begin{array}{cccl}
     H^0(\widetilde{\varphi}^{-1}(O_0),
       {\cal L}_{\Delta^{\prime}}|_{\widetilde{\varphi}^{-1}(O_0)})\;
      & \hookrightarrow
      & H^0(O_0, {\cal O}_{O_0})
            \otimes\,H^0(F_0, {\cal L}_{\overline{\Delta^{\prime}}})
                                                       &   \\[.6ex]
     s_{m^{\prime}} & \mapsto
      & s_{\xi^{\dagger}(m^{\prime})}
         \otimes s_{\overline{\varphi^{\dagger}}(m^{\prime})}  & ,
   \end{array}
  $$
  where $m^{\prime}\in\Delta^{\prime}$.
  (Here we are identifying ${\cal O}_{O_0}$ with  
  ${\cal L}_{\xi^{\dagger}(\Delta^{\prime})}|_{O_0}$ implicitly.)
  
  \hspace{-1.9em}(c)\hspace{1ex}
  Let $\xi_1,\,\xi_2:N\rightarrow N^{\prime}$ be two homomorphisms
  such that $\varphi\circ\xi_i=\Id_N$, $i=1,\,2$. Then, the
  $H^0(O_0, {\cal O}_{O_0})$-factor of the factorizations in Part (b)
  with respect to $\xi_1$ and $\xi_2$ respectively are related by
  $s_{\xi_2^{\dagger}(m^{\prime})}
   =\chi^{(\xi_2^{\dagger}-\xi_1^{\dagger})(m^{\prime})}\,
                                   s_{\xi_1^{\dagger}(m^{\prime})}$
  while the $H^0(F_0,{\cal L}_{\overline{\Delta^{\prime}}})$-factors
  are the same.
 \end{quote}
} 

\bigskip

\noindent
{\it Proof.}
 As a toric variety, $\widetilde{\varphi}^{-1}(O_0)$ is described by
 the fan $\Sigma^{\prime}_{\{0\}}$ in $N^{\prime}_{\scriptsizeBbb R}$
 while $F_0$ is described by the same fan but in the subspace
 $\varphi^{-1}_{\scriptsizeBbb R}(0)$ in
 $N^{\prime}_{\scriptsizeBbb R}$. The map $\xi$ induces a projection
 map $\pi_{\xi}:\widetilde{\varphi}^{-1}(O_0)\rightarrow F_0$ that 
 projects the affine chart $U_{\sigma^{\prime}}$ of
 $\widetilde{\varphi}^{-1}(O_0)$ to the affine chart
 $\overline{U}_{\sigma^{\prime}}$ of $F_0$ for
 $\sigma^{\prime}\in\Sigma^{\prime}_{\{0\}}$.
 It also induces a $T_N$-action on $\widetilde{\varphi}^{-1}(O_0)$
 that descends under $\pi_{\xi}$ to the identity group action on
 $F_0$. Together, this specifies a decomposition
 $\widetilde{\varphi}^{-1}(O_0)=O_0\times F_0$.
 Let $h^{\prime}\in\SF(\Sigma^{\prime})$ be a piecewise linear
 function on $N^{\prime}_{\scriptsizeBbb R}$ that determines
 ${\cal L}_{\Delta^{\prime}}$. Then
 ${\cal L}_{\Delta^{\prime}}|_{\widetilde{\varphi}^{-1}(O_0)}$
 is determined by $h^{\prime}|_{\varphi_{\scriptsizeBbb R}^{-1}(0)}$,
 whose corresponding polytope in $M^{\prime}$ can be chosen to be 
 $i_{\xi}(\overline{\Delta^{\prime}})$. This implies that
 ${\cal L}_{\Delta^{\prime}}|_{\widetilde{\varphi}^{-1}(O_0)}\,
    \simeq\, \pr_2^{\ast}\,{\cal L}_{\overline{\Delta^{\prime}}}$
 and, hence, proves Part (a).

 For Part (b), first observe that the set of affine charts
 ${\cal U}=\{ U_{\sigma^{\prime}} |
                     \sigma^{\prime}\in\Sigma^{\prime}_{\{0\}} \}$
 forms a covering of $\widetilde{\varphi}^{-1}(O_0)$.
 Over the affine chart $U_{\sigma^{\prime}}$,
 $\sigma^{\prime}\in\Sigma^{\prime}_{\{0\}}$, of
 $\widetilde{\varphi}^{-1}(O_0)$, one can decompose the section
 $s_{m^{\prime}}$ in
 ${\cal L}_{\Delta^{\prime}}|_{\widetilde{\varphi}^{-1}(O_0)}$
 associated to $m^{\prime}\in\Delta^{\prime}\cap M^{\prime}$ as
 $$
  s_{m^{\prime}}|_{U_{\sigma^{\prime}}}\;
   =\;\chi^{m^{\prime}-m^{\prime}_{\sigma^{\prime}}}\;
   =\; \chi^{\xi^{\dagger}(m^{\prime})
                -\xi^{\dagger}(m^{\prime}_{\sigma^{\prime}})}\,
     \cdot\,
     \chi^{\overline{\varphi^{\dagger}}(m^{\prime})
      -\overline{\varphi^{\dagger}}(m^{\prime}_{\sigma^{\prime}})}
   =\; \chi^{\xi^{\dagger}(m^{\prime})}\,
     \cdot\,
     \chi^{\overline{\varphi^{\dagger}}(m^{\prime})
      -\overline{\varphi^{\dagger}}(m^{\prime}_{\sigma^{\prime}})}
 $$
 since $\xi^{\dagger}\circ i_{\xi}=0$. One notes that
 $\chi^{\xi^{\dagger}(m^{\prime})}$ is the section
 $s_{\xi^{\dagger}(m^{\prime})}$ of ${\cal O}_{O_0}$ and
 $\chi^{\overline{\varphi^{\dagger}}(m^{\prime})
        -\overline{\varphi^{\dagger}}(m^{\prime}_{\sigma^{\prime}})}$
 is the section $s_{\overline{\varphi^{\dagger}}(m^{\prime})}$ of
 ${\cal L}_{\overline{\Delta^{\prime}}}$ restricted to
 $\overline{U}_{\sigma^{\prime}}$. Since the transition map between
 the affine charts of $\widetilde{\varphi}^{-1}(O_0)$ in ${\cal U}$
 is trivial along the $O_0$-component of the product decomposition
 of $\widetilde{\varphi}^{-1}(O_0)$, this proves Part (b).
 
 Finally, Part (c) follows from the proof of Part (b).
 This concludes the proof.

\noindent\hspace{12cm} $\Box$

\bigskip

Let $\xi_1,\,\xi_2:N\rightarrow N^{\prime}$ be two homomorphisms
such that $\varphi\circ\xi_i=\Id_N$, $i=1,\,2$.
Then the difference $\xi_2-\xi_1$ induces morphisms
$$
 \begin{array}{cccl}
  N           & \stackrel{\xi_{21}}{\longrightarrow}
              & \varphi^{-1}(0)                        &  \\[.6ex]
  M^{\prime}/\mbox{\raisebox{-.4ex}{$\varphi^{\dagger}(M)$}}
              & \stackrel{\xi_{21}^{\dagger}}{\longrightarrow}
              & M                                      &  \\[.6ex]
  T_N         & \stackrel{\xi_{21}^T}{\longrightarrow}
              & T_{\varphi^{-1}(0)}                    &  \\[.6ex]
  \{\,\mbox{character of $T_{\varphi^{-1}(0)}$}\,\}
              & \stackrel{\xi_{21}^{\chi}}{\longrightarrow}
              & \{\,\mbox{character of $T_N$}\,\}      &.
 \end{array}
$$
In terms of these, the transition factor
$\chi^{(\xi_2^{\dagger}-\xi_1^{\dagger})(m^{\prime})}$ in Part (c)
of Lemma 3.2.1, as a character of $T_N$, is simply
$\xi_{21}^{\chi}(\chi^{\overline{\varphi^{\dagger}}(m^{\prime})})
   =\chi^{\overline{\varphi^{\dagger}}(m^{\prime})}\circ\xi_{21}^T$.
Thus, Lemma 3.2.1 and Remark 3.1.3 together imply$\,$:

\bigskip

\noindent
{\bf Corollary 3.2.2 [torically equivalent restriction].} {\it 
 The description of the restriction of a section of
 ${\cal L}_{\Delta^{\prime}}$ to a fiber $F_0$ over $O_0$ by
 different choices of $\xi$ give rise to torically equivalent
 sections in $H^0(F_0, {\cal L}_{\overline{\Delta^{\prime}}})$.
 } 

\bigskip

\noindent
In particular, when one considers the hypersurfaces from the
zero-locus of sections of ${\cal L}_{\Delta^{\prime}}$, though
a $\xi$ has to be chosen to give an explicit description of its
restriction to the generic fiber, the choice of $\xi$ does not change
the toric equivalence class of the hypersurfaces in the fibers.

\bigskip

\noindent
{\it Remark 3.2.3.}
The situation here is pretty much like local trivialization of
a bundle: There is no canonical one; however, different local
trivializations are related by automorphisms of the fiber.

\bigskip

\begin{flushleft}
{\bf Case (2)$\,$:
     Restriction to exceptional fibers for $\varphi$ surjective.}
\end{flushleft}
Recall from Remark 2.1.12 that, in this case,
$\Ind(\sigma)
 =[\,N/\mbox{\raisebox{-.4ex}{$N_{\sigma}$}}\,:\,
    \overline{\varphi}^{\sigma^{\prime}}_{\sigma}(N^{\prime}/
     \mbox{\raisebox{-.4ex}{$N^{\prime}_{\sigma^{\prime}}$}})\,]
 =1
$
for any $\sigma\in\Sigma$ and
$\sigma^{\prime}\in\Sigma^{\prime}_{\sigma}$. Hence,
$\overline{\varphi}^{\sigma^{\prime}}_{\sigma}:
   N^{\prime}/\mbox{\raisebox{-.4ex}{$N^{\prime}_{\sigma^{\prime}}$}}
                  \rightarrow N/\mbox{\raisebox{-.4ex}{$N_{\sigma}$}}$
is surjective and 
$\widetilde{\overline{\varphi}}^{\sigma^{\prime\;-1}}_{\sigma}(O_{\sigma})
 =\widetilde{\varphi}^{-1}(O_{\sigma})=O_{\sigma}\times F^c_{\sigma}$.
Let $\tau^{\prime}$ be a primitive element in $\Sigma^{\prime}_{\sigma}$,
then the orbit closure $V(\tau^{\prime})$ of $\tau^{\prime}$ in
$X_{\Sigma^{\prime}}$ is an irreducible component of the closure of 
$\widetilde{\varphi}^{-1}(O_{\sigma})$ in $X_{\Sigma^{\prime}}$
and $F_{\sigma}^{\tau^{\prime}}$ is a generic fiber of the restriction
$\widetilde{\varphi}|_{V(\tau^{\prime})}:
                       V(\tau^{\prime})\rightarrow V(\sigma)$.
Consequently, for each pair $(\sigma, \tau^{\prime})$, where
$\sigma\in\Sigma$ and $\tau^{\prime}\in\Sigma_{\sigma}$ primitive,
one can apply the discussion in Preparations to obtain the
restriction of ${\cal L}_{\Delta^{\prime}}$ and the sections therein
to $V({\tau^{\prime}})$ and then the discussions in
Case (1) [fundamental] to this bundle over $V(\tau^{\prime})$.
In this way, one can write down how the generic fiber of the morphism
$Y^{\prime}\cap V(\tau^{\prime})\rightarrow V(\sigma)$
varies along $O(\sigma)$.

Since each irreducible component of a fiber of $\widetilde{\varphi}$
must appear as a generic fiber for some $(\sigma,\tau^{\prime})$,
the discussions in Case (1) and above together gives us a description
of every fiber of $Y^{\prime}\rightarrow X_{\Sigma}$ orbit by orbit. 

\bigskip

Together this concludes the discussion for the case
$\varphi:N^{\prime}\rightarrow N$ surjective.

\bigskip

\begin{flushleft}
{\bf Case (3)$\,$:
     Restriction to fibers when the index
                            $[N:\varphi(N^{\prime})]$ is finite.}
\end{flushleft}
The maps $N^{\prime}\rightarrow\varphi(N^{\prime})\hookrightarrow N$
induce a decomposition of
$\widetilde{\varphi}:
  X_{\Sigma^{\prime}}
   \stackrel{\widetilde{\varphi}_1}{\rightarrow} X_{\Sigma}
              \stackrel{\widetilde{\iota}_1}{\rightarrow} X_{\Sigma}$,
where the first $\Sigma$ is a fan in
$\varphi(N^{\prime})_{\scriptsizeBbb R}$ while
the second $\Sigma$ is in $N_{\scriptsizeBbb R}$.
The induced fibration
$\widetilde{\varphi}_1:Y^{\prime}\rightarrow X_{\Sigma}$ is already
discussed in Cases (1) and (2) above.
Thus, it remains to understand the map $\widetilde{\iota}_1$,
which is a branched covering of order $[N:\varphi(N^{\prime})]$
and is a special case of our discussion in Sec.\ 2.1.
The number of sheets over $O_{\sigma}$, $\sigma\in\Sigma$,
is the same as $\Ind(\sigma)$, which must satisfy
$\Ind(\sigma)\le [N:\varphi(N^{\prime})]$.
The branched locus $\Br(\widetilde{\varphi})$ of $\widetilde{\varphi}$
is a union of $V(\sigma)$ with $\Ind(\sigma)<[N:\varphi(N^{\prime})]$.
The fiber of $\widetilde{\varphi}:Y^{\prime}\rightarrow X_{\Sigma}$
over a point $p$ is then the disjoint union 
$\cup_{p_i\in\widetilde{\iota}_1^{-1}(p)}\,
             \widetilde{\varphi}_1|_{Y^{\prime}}^{-1}(p_i)$.
This concludes the discussion.

\bigskip

\begin{flushleft}
{\bf Case (4)$\,$:
     Restriction to fibers when the index
                            $[N:\varphi(N^{\prime})]$ is infinite.}
\end{flushleft}
This is the remaining case after Cases (1), (2), and (3); and 
it can be reduced to the discussions in Cases (1), (2), and (3)
by considering the maps of fans
$\Sigma^{\prime}
   \stackrel{\varphi_2}{\rightarrow} \Sigma_{\varphi}
                          \stackrel{\iota_2}{\rightarrow} \Sigma$
induced by
$\varphi:N^{\prime}\rightarrow
   N\cap\varphi(N^{\prime}_{\scriptsizeBbb R})\hookrightarrow N$.
This gives a decomposition
$\widetilde{\varphi}=\widetilde{\varphi}_2\circ\widetilde{\iota}_2$. 
The surjective map $\widetilde{\varphi}_2|_{Y^{\prime}}$ is already
discussed in Case (3). Since our goal is to understand the induced
fibration of $Y^{\prime}$, we may replace $\widetilde{\varphi}$
by $\widetilde{\varphi}_2$. This concludes our discussion.

\bigskip

\bigskip

In the string literatures, e.g.\ [K-S1] and [K-S2], 
Cox homogeneous coordinates are used to understand the fibrations
of toric Calabi-Yau hypersurfaces. For completeness of discussion,
we recall in the following remark their essential points and fit
them into the setting above.

\bigskip

\noindent 
{\it Remark 3.2.4 [Cox homogeneous coordinates].}
Assume that 
$\widetilde{\varphi}:X_{\Sigma^{\prime}}\rightarrow X_{\Sigma}$
is a fibration. Let 
$\Sigma^{\prime}(1)
  =\{v^{\prime}_1,\,\cdots,\,v^{\prime}_{I^{\prime}}\}$,
$\Sigma(1)=\{v_1,\,\cdots,\, v_I\}$, and
$(z^{\prime}_1,\,\cdots,\,z^{\prime}_{I^{\prime}})$ and
$(z_1,\,\cdots,\, z_I)$ be the corresponding homogeneous coordinates
for $X_{\Sigma^{\prime}}$ and $X_{\Sigma}$ respectively. Since
$\varphi:\Sigma^{\prime}(1)\rightarrow\Sigma(1)$ is now surjective,
$\varphi(v^{\prime}_i)=r_{ij}z_j$ for some $r_{ij}\in {\Bbb Z}_{\ge 0}$
and $\varphi$ determines a polynomial map from ${\Bbb C}^{I^{\prime}}$
to ${\Bbb C}^I$ defined by
$$
 (z^{\prime}_1,\,\cdots,\,z^{\prime}_{I^{\prime}})\;
  \longrightarrow\;
 (z_1,\,\cdots,\, z_I)\,
 =\,(\,\prod_i\,{z^{\prime}_i}^{r_{i1}},\,\cdots,\,
       \prod_i\,{z^{\prime}_i}^{r_{iI}} \,)\,.
$$
This map is equivariant with respect to the group actions involved
and indeed descends to the toric morphism
$\widetilde{\varphi}:X_{\Sigma^{\prime}}\rightarrow X_{\Sigma}$.
Let $Y^{\prime}$ be a $\Sigma^{\prime}$-regular (cf.\ Sec.\ 4.1)
hypersurface in $X_{\Sigma^{\prime}}$ in the linear system
${\cal O}_{X_{\Sigma^{\prime}}}(D^{\prime})$, where
$D^{\prime}=\sum_{i}\,a_iD_{v^{\prime}_i}$. Then, in terms of the
given homogeneous coordinates, $Y^{\prime}$ is described by the zero
locus of the polynomial
$$
 {z^{\prime}}^{D^{\prime}}\,
 \sum_{m^{\prime}\in\Delta^{\prime}\cap M^{\prime}}\, c_{m^{\prime}}\,
   \prod_{i=1}^{|\Sigma^{\prime}(1)|}\,
               {z_{i}^{\prime}}^{m^{\prime}(v_i^{\prime})}\;
 =\; \sum_{m^{\prime}\in\Delta^{\prime}\cap M^{\prime}}\, c_{m^{\prime}}\,
      \prod_{i=1}^{|\Sigma^{\prime}(1)|}\,
               {z_{i}^{\prime}}^{m^{\prime}(v_i^{\prime})+ a_i}\,,
$$
where $\Delta^{\prime}$ is the polytope in $M^{\prime}$ determined
by $D^{\prime}$. To manifest the fibration $\widetilde{\varphi}$, 
rewrite this polynomial as (cf.\ Equation (18) in [K-S2])
$$
 \sum_{\overline{\varphi^{\dagger}}(m^{\prime})
                                 \in\overline{\Delta^{\prime}}}\,
 c_{\overline{\varphi^{\dagger}}(m^{\prime})}\,
  \prod_{v^{\prime}_k\in\varphi^{-1}(0)}\,
       {z^{\prime}_k}^{\langle\,
                          \overline{\varphi^{\dagger}}(m^{\prime})\,,
                                               \,v^{\prime}_k\,\rangle}
$$
with
$$
 c_{\overline{\varphi^{\dagger}}(m^{\prime})}\;
 =\; \sum_{m^{\prime}\in \overline{\varphi^{\dagger}}(m^{\prime})}\,
      c_{m^{\prime}}\,
      \prod_{v^{\prime}_k\notin \varphi^{-1}(0)}\,
       {z^{\prime}}_k^{m^{\prime}(v^{\prime}_k)+a_k}\,.
$$
Formally, this is an equation for the fiber with coefficients
$c_{\overline{\varphi^{\dagger}}(m^{\prime})}$ regarded as polynomial
of coordinates of the base $X_{\Sigma}$.

To realize this explicitly, pick a $\xi:N\rightarrow N^{\prime}$
as in our discussion of induced morphism
$Y^{\prime}\subset X_{\Sigma^{\prime}}\rightarrow X_{\Sigma}$.
The map $\xi$ determines a decomposition into three factors$\,$:
\begin{eqnarray*} 
 \lefteqn{
   \prod_{v^{\prime}_k\notin \varphi^{-1}(0)}\,
               {z^{\prime}}_k^{m^{\prime}(v^{\prime}_k)+a_k} } \\
    & & =\; \prod_{v_j\in\Sigma(1)}\,
       z_j^{\langle\xi^{\dagger}m\,,\,v_j\rangle}\,.
     \prod_{v^{\prime}_k\notin \varphi^{-1}(0)}\,
       {z^{\prime}}_k^{\langle\,
                         \overline{\varphi^{\dagger}}(m^{\prime})\,,\,
                         v^{\prime}_k-\xi\circ\varphi(v_k)\, \rangle}\,.
     \prod_{v^{\prime}_k\notin \varphi^{-1}(0)}\,
     {z^{\prime}}_k^{\,a_k}\,.
\end{eqnarray*} 
Observe that the first factor is the polynomial in the homogeneous
coordinates of $X_{\Sigma}$ associated to $\xi^{\dagger}(m)\in M$
and is the function of the orbit $O_{\{0\}}$ of the base that appear
in the proof of Lemma 3.2.1, the last factor is an overall factor
that is independent of $\overline{\varphi^{\dagger}}(m^{\prime})$,
and the middle factor is projected out in this correspondence.
This gives a correspondence of our discussion and the ones in
homogeneous coordinates for the fiber over $O_{\{0\}}$ as seen
in physics literature. This concludes the remark.

\bigskip

\noindent
{\it Remark.\ 3.2.5.}
When $\widetilde{\varphi}:Y^{\prime}\rightarrow X_{\Sigma}$ is
a fibration, though in principle one can write down the equation
of the discriminant locus of $\widetilde{\varphi}$ but practically
it remains very difficult to understand its details
(cf.\ Example 4.2.1).

\bigskip
             
\section{Fibration of Calabi-Yau hypersurfaces via toric morphisms.}

The theme of this section is a computational scheme given and
 demonstrated in Sec.\ 4.2.
Before getting into it, we give a couple of related remarks
 in Sec.\ 4.1.

\bigskip

\subsection{General remarks.}

\noindent
{\it Remark 4.1.1 [Batyrev].} We recall some definitions and facts
about toric hypersurfaces from [Ba2] that will be needed in Sec.\ 4.2.
\begin{itemize}
 \item
  Let $\Sigma^{\prime\prime}$ be a fan in
  $N^{\prime\prime}_{\scriptsizeBbb R}$.
  Then the toric orbit decomposition of the ambient toric variety
  $X_{\Sigma^{\prime\prime}}$ induces a stratification of
  a hypersurface $Y^{\prime\prime}$ in $X_{\Sigma^{\prime\prime}}\,$:
  $Y^{\prime\prime}
   =\cup_{\sigma^{\prime\prime}\in\Sigma^{\prime\prime}}\,
                      Y^{\prime\prime}_{\sigma^{\prime\prime}}$.
 $Y^{\prime\prime}$ is called {\it $\Sigma^{\prime\prime}$-regular}
  if each stratum $Y^{\prime\prime}_{\sigma^{\prime\prime}}$ is
  either empty or a smooth subvariety of codimension $1$ on
  $O_{\sigma^{\prime\prime}}$.
 In this case, this decomposition is a refinement of
  the stratification of the hypersurface by the analytic isomorphism
  classes of singularities.
 When $\Sigma^{\prime\prime}$ is the normal fan of a polytope
  $\Delta^{\prime}$, a $\Sigma^{\prime\prime}$-regular
  hypersurface is also called a $\Delta^{\prime}$-regular
  hypersurface.

 \item
 {\bf Fact [$\Delta^{\prime}$-regular dense].}
 ([Proposition 3.1.3 in [Ba2].) {\it
 The set of $\Delta^{\prime}$-regular hypersurfaces is a Zariski
 open subset in ${\Bbb P}(L(\Delta^{\prime}))$, where
 $L(\Delta^{\prime})$ is the space of Laurent polynomials with
 Newton polytope $\Delta^{\prime}$.
 For such hypersurfaces, the singularities are induced from the
 ambient toric variety.
 } 

 \item
 The local model for the transverse singularity along the stratum
 $Y^{\prime\prime}_{\sigma^{\prime\prime}}$ of the induced
 stratification is given by the affine toric variety
 $U_{\sigma^{\prime\prime},\,
                     N^{\prime\prime}_{\sigma^{\prime\prime}}}$
 associated to the cone $\sigma^{\prime\prime}$ regarded as a fan
 in ${N^{\prime\prime}_{\sigma^{\prime\prime}}}_{\scriptsizeBbb R}$.
\end{itemize}

\bigskip

\noindent
{\bf Lemma 4.1.2 [regularity].} {\it
 Let $\Sigma^{\prime}$ corresponds to a maximal projective
 triangulation of the dual polytope ${\Delta^{\prime}}^{\ast}$ of
 $\Delta^{\prime}$, then a $\Delta^{\prime}$-regular hypersurface
 in $X_{\Delta^{\prime}}$ becomes a $\Sigma^{\prime}$-regular
 hypersurface in $X_{\Sigma^{\prime}}$ after the associated maximal
 projective crepant partial (MPCP-)desingularization.
} 

\bigskip

\noindent
{\it Proof.}
Let $\Sigma_{\Delta^{\prime}}$ be the normal fan of $\Delta^{\prime}$,
which consists of cones over the faces of ${\Delta^{\prime}}^{\ast}$,
and consider the toric morphism 
$\zeta:X_{\Sigma^{\prime}}\rightarrow X_{\Delta^{\prime}}$
associated to the triangulation of ${\Delta^{\prime}}^{\ast}$.
Recall from Proposition 2.1.4 and its proof that
 the exceptional locus of
 $\zeta:X_{\Sigma^{\prime}}\rightarrow X_{\Delta^{\prime}}$
 is a union of the orbit closures of $X_{\Sigma^{\prime}}$ associated
 to primitive cones in $\Sigma^{\prime}$ with respect to $\zeta$
and that, over an orbit $O_{\sigma^{\prime\prime}}$ of
$X_{\Delta^{\prime}}$,
 $\zeta^{-1}(O_{\sigma^{\prime\prime}})
  =O_{\sigma^{\prime\prime}}\times X(\sigma^{\prime\prime})$
 for some possibly reducible toric variety $X(\sigma^{\prime\prime})$
and is the union of toric orbits
$\cup_{\sigma^{\prime}
    \subset\mbox{\scriptsize\it Int}(\sigma^{\prime\prime})}\,
                                          O_{\sigma^{\prime}}$
in $X_{\Sigma^{\prime}}$.
Consider now the restriction to hypersurfaces
$\zeta: Y^{\prime}\rightarrow Y^{\prime\prime}$.
For $\sigma^{\prime}\subset\Int(\sigma^{\prime\prime})$,
$Y^{\prime}_{\sigma^{\prime}}\,
 = \zeta^{-1}(Y^{\prime\prime})\cap O_{\sigma^{\prime}}\,
 = Y^{\prime\prime}_{\sigma^{\prime\prime}}
   \times O_{\sigma^{\prime},\sigma^{\prime\prime}}$,
where $O_{\sigma^{\prime},\sigma^{\prime\prime}}$ is the toric
orbit in $X(\sigma^{\prime\prime})$ associated to $\sigma^{\prime}$.
From this, $Y^{\prime}_{\sigma^{\prime}}$ must be empty if
$Y^{\prime\prime}_{\sigma^{\prime\prime}}$ is empty or smooth of
codimension $1$ in $O_{\sigma^{\prime}}$ if
$Y^{\prime\prime}_{\sigma^{\prime\prime}}$ is smooth of codimension
$1$ in $O_{\sigma^{\prime\prime}}$. This concludes the proof.

\noindent\hspace{12cm} $\Box$

\bigskip

\subsection{The computational scheme and a detailed study of an example.}

The discussions in the previous subsections together provide us with
a toric computational scheme to study the induced fibration of a
toric Calabi-Yau hypersurface, as outlined below$\,$:

\bigskip

\baselineskip 12pt
{\footnotesize
\begin{itemize}
 \item
 The polytope $\Delta^{\prime}$
 and the toric variety $X_{\Sigma^{\prime}}\,$:
 \begin{itemize}
  \item
  The reflexive polytope $\Delta^{\prime}$.

  \item
  The lattice points in $\Delta^{\prime}\cap M^{\prime}$.

  \item
  The dual reflexive polytope ${\Delta^{\prime}}^{\ast}$ in
   $N^{\prime}_{\tinyBbb R}$.

  \item
  The triangulation of ${\Delta^{\prime}}^{\ast}$ that gives
   $\Sigma^{\prime}$.

  \item
  The singular locus of $X_{\Sigma^{\prime}}$.
 \end{itemize}

 \item
 The line bundle ${\cal L}_{\Delta^{\prime}}$ and its restrictions$\,$:
 \begin{itemize}
  \item
  The restriction of ${\cal L}_{\Delta^{\prime}}$ to
  $V(\tau^{\prime})$ for $\tau^{\prime}$ primitive in
  $\Sigma^{\prime}_{\sigma}$ for some $\sigma\in\Sigma\,$:
  \begin{itemize}
   \item
   the polytope $\Delta^{\prime}_{\tau^{\prime}}$ in
   $\tau^{\prime\perp}_{\tinyBbb R}$, up to translation.

   \item
   the lattice points in
   $\Delta^{\prime}_{\tau^{\prime}}\cap\tau^{\prime\perp}$.
  \end{itemize}

  \item
  The irreducible component $F_{\sigma}^{\tau^{\prime}}$ of
  a generic fiber of
  $\widetilde{\varphi}:V(\tau^{\prime})\rightarrow V(\sigma)$.

  \item
  The restriction ${\cal L}_{F_{\sigma}^{\tau^{\prime}}}$ of
   ${\cal L}$ to $F_{\sigma}^{\tau^{\prime}}$:
  \begin{itemize}
   \item
   the polytope $\overline{\Delta^{\prime}_{\tau^{\prime}}}$ in
    ${\tau^{\prime}}^{\perp}/
      \mbox{\raisebox{-.4ex}{$\varphi^{\dagger}(\sigma^{\perp})$}}$.

   \item
   The lattice points of the polytope
   $\overline{\Delta^{\prime}_{\tau^{\prime}}}$ in
   ${\tau^{\prime}}^{\perp}/
     \mbox{\raisebox{-.4ex}{$\varphi^{\dagger}(\sigma^{\perp})$}}$.
  \end{itemize}
 \end{itemize}

 \item
 The discriminant locus associated to
 $\overline{\Delta^{\prime}_{\tau^{\prime}}}$ in
 ${\tau^{\prime}}^{\perp}/
      \mbox{\raisebox{-.4ex}{$\varphi^{\dagger}(\sigma^{\perp})$}}$.
 (I.e.\ the discriminant locus of the family of sections of the
  line bundle ${\cal L}_{\overline{\Delta^{\prime}_{\tau^{\prime}}}}$
  over $F_{\sigma}^{\tau^{\prime}}$.)

 \item
 The $\Sigma^{\prime}$-regular sections in ${\cal L}_{\Delta^{\prime}}$
 and their restrictions$\,$:
 \begin{itemize}
  \item
  A $\Sigma^{\prime}$-regular section $s$ of
  ${\cal L}_{\Delta^{\prime}}$
  and its zero-locus, the Calabi-Yau hypersurface $Y^{\prime}$.

  \item
  The singular loci of $Y^{\prime}$.

  \item
  The restriction of $s$ to $V(\tau^{\prime})$ for
  $\tau^{\prime}$ primitive in $\Sigma^{\prime}_{\sigma}$
  for some $\sigma\in\Sigma$.

  \item
  The restriction of $s$ to $F_{\sigma}^{\tau^{\prime}}$:
  \begin{itemize}
   \item
   description of the irreducible components of a generic fiber of
    $\widetilde{\varphi}:Y^{\prime}\rightarrow X_{\Sigma}$
    over each of the toric orbits of $X_{\Sigma}$.

   \item
   restriction of a section of ${\cal L}_{\Delta^{\prime}}$ to
    $V(\tau^{\prime})$ and then rewrite it in the fibred form
    with respect to a map
    $\xi_{\sigma}^{\tau^{\prime}}:
      N/\mbox{\raisebox{-.4ex}{$N_{\sigma}$}} \rightarrow
      N^{\prime}/\mbox{\raisebox{-.4ex}{$N^{\prime}_{\tau^{\prime}}$}}$.
  \end{itemize}
 \end{itemize}

 \item
 The induced fibration $Y^{\prime}\rightarrow X_{\Sigma}$
 and the stratification of the discriminant locus of the fibration.

 \item
 Beyond$\,$:
 \begin{itemize}
  \item
  Note that before resolving the singularities of $Y^{\prime}$,
  all the fibers, generic and special alike, are described as
  hypersurfaces of the (perhaps reducible) fiber toric varieties.

  \item
  Fibers after crepant resolution of the singularities of
  $Y^{\prime}$.
  Note that now some of the fibers may not be realizable as
  hypersurfaces of fibers of
  $X_{\Sigma^{\prime}}\rightarrow X_{\Sigma}$.

  \item
  Monodromy of the fibration.

  \item
  Other details.
 \end{itemize}
\end{itemize}
} 

\baselineskip 14pt

\bigskip

\noindent
The sub-items in Item Beyond are not discussed in general in the
current work. We will study some part of it in the specific example
below. The purely toric items of this scheme should be programmable,
(cf.\ Problem 5.1 in Sec.\ 5).

We now apply this to the example of elliptic Calabi-Yau $4$-folds
in [B-C-dlO-G].

\bigskip

\noindent
{\bf Example 4.2.1 [B-C-dlO-G].} (Continuing Example 2.2.2.)
Let us now consider the induced fibration of a $4$-dimensional
Calabi-Yau hypersurface $Y^{\prime}$ in $X_{\Sigma^{\prime}}$ over
$X_{\Sigma}$ from Example 2.2.2.
All the notations follow from that example.

\begin{itemize}
 \item
 {\it The polytope $\Delta^{\prime}$
 and the toric variety $X_{\Sigma^{\prime}}\,$}:
 \begin{itemize}
  \item
  {\it The reflexive polytope $\Delta^{\prime}\,$}:
  \begin{quote}
   \hspace{-1.9em}(1)\hspace{1ex}
   $\Delta^{\prime}$ is the convex hull in
   $M^{\prime}_{\scriptsizeBbb R}$ of the following set of
   $14$ generating vertices in $M^{\prime}\,$:
   {\footnotesize
   $$
    \left\{
     \begin{array}{lll}
      m^{\prime}_1=(-22,-14,4,1,1)\,, & m^{\prime}_2=(-22,6,4,1,1)\,,\\[.6ex]
      m^{\prime}_3=(-10,-6,2,-1,1)\,, & m^{\prime}_4=(-10,2,2,-1,1)\,,
                                                                     \\[.6ex]  
      m^{\prime}_5=(-6, -6,0,1,1)\,, & m^{\prime}_6=(0,0,0,-2,1)\,,\\[.6ex]
      m^{\prime}_7=(0,0,0,1,-1)\,, & m^{\prime}_8 =(2, -6, 2,-1,1)\,,\\[.6ex]
      m^{\prime}_9  =(2,2,2,-1,1)\,, 
                         & m^{\prime}_{10}=(6, -14,4,1,1)\,, \\[.6ex]
      m^{\prime}_{11}=(6,-6,0,1,1)\,,
         & m^{\prime}_{12}  =(6,3,-3,1,1)\,,        \\[.6ex]
      m^{\prime}_{13}=(6,6,-3,1,1)\,,    & m^{\prime}_{14}=(6,6,4,1,1)
     \end{array}
    \right\}\,.
   $$
   } 

   \hspace{-1.9em}(2)\hspace{1ex}
    There are $9$ codimension-1 faces, i.e.\ the maximal polytope
    of the boundary $\partial\Delta^{\prime}\,$:
    (For brevity of notaion, only the indices $i$ of the generating
     $m_i^{\prime}$ are indicated.
     Their dual vertices in $N^{\prime}$ are labelled as subscripts.)
    {\footnotesize
    $$
     \left\{
      \begin{array}{lll}
       [6,  7, 8, 9, 10, 11, 12 , 13, 14]_{v_1^{\prime}}\,,
         & [2, 4, 6, 7, 9, 13, 14 ]_{v_2^{\prime}}\,, \\[.6ex]
       [3, 4, 6, 7, 8, 9, ]_{c_1^{\prime}}\,,
         & [1, 2, 3, 4, 7, 8, 9, 10, 14 ]_{c_2^{\prime}}\,, \\[.6ex]
       [1, 2, 5, 7, 10, 11, 12, 13, 14]_{v_4^{\prime}}\,,
         & [1, 2, 3, 4, 5, 6, 8, 9, 10, 11, 12, 13, 14]_{v_5^{\prime}}\,,
                                               \\[.6ex]
       [1, 3, 5, 6, 7, 8, 10, 11,]_{f^{\prime}}\,,
         & [5, 6, 7, 11, 12]_{g^{\prime}}\,, \\[.6ex]
       [1, 2, 3, 4, 5, 6, 7, 12, 13]_{v_6^{\prime}}
      \end{array}
     \right\}\,.
    $$
    {\normalsize Faces} } 
    of higher codimensions are obtained by the intersection
    of subcollections of the above faces.
  \end{quote}
  \hspace{1ex}

  \item
  {\it The lattice points in $\Delta^{\prime}\cap M^{\prime}\,$}:
  \begin{quote}
   The cardinality $|\Delta^{\prime}\cap M^{\prime}| = 3365$.
   Some sample lattice points in $\Delta^{\prime}\cap M^{\prime}$ are
   {\footnotesize
    $$
     \begin{array}{rllll}
        & (-22, -14, 4, 1, 1)\,, & (-22, -13, 4, 1, 1)\,,
        & (-22, -12, 4, 1, 1)\,, & \cdots        \\[.6ex]
      \cdots\,,
        & (0, 0, 0, -1, 1)\,,    & (0, 0, 0, 0, 0)\,,
        & (0, 0, 0, 0 , 1)\,,    & \cdots        \\[.6ex]
      \cdots\,,
        & (6, 6, 2, 1, 1)\,,
        & (6, 6, 3, 1, 1)\,,     & (6, 6, 4, 1, 1)\,.  &
    \end{array}
   $$
   } 
  \end{quote}
  \vspace{1ex}

  \item
  {\it The dual reflexive polytope ${\Delta^{\prime}}^{\ast}$ in
       $N^{\prime}_{\scriptsizeBbb R}\,$}:
  \begin{quote}
   \hspace{-1.9em}(1)\hspace{1ex}
   ${\Delta^{\prime}}^{\ast}$ is the convex hull in
   $N^{\prime}_{\scriptsizeBbb R}$ of the following set of
   $9$ generating vertices in $N^{\prime}\,$:
   $\{\, v_1^{\prime},\; v_2^{\prime},\; c_1^{\prime},\;
         c_2^{\prime},\; v_4^{\prime},\; v_5^{\prime},\;
         f^{\prime},\; g^{\prime},\; v_6^{\prime}  \,\}$.
   \vspace{1ex}

   \hspace{-1.9em}(2)\hspace{1ex}
   There are $14$ codimension-1 faces$\,$:
   (Their dual vertices in $M^{\prime}$ are labelled as subscripts.)
    {\footnotesize
    $$
     \left\{
      \begin{array}{ll}
       [c_2^{\prime}, v_4^{\prime}, v_5^{\prime}, f^{\prime},
                      v_6^{\prime}]_{m_1^{\prime}}\,,
         & [v_2^{\prime}, c_2^{\prime}, v_4^{\prime}, v_5^{\prime},
                          v_6^{\prime}]_{m_2^{\prime}}\,,   \\[.6ex]
       [c_1^{\prime}, c_2^{\prime}, v_5^{\prime}, f^{\prime},
                          v_6^{\prime}]_{m_3^{\prime}}\,,   
         & [v_2^{\prime}, c_1^{\prime}, c_2^{\prime}, v_5^{\prime},
                      v_6^{\prime}]_{m_4^{\prime}}\,,       \\[.6ex]
       [v_4^{\prime}, v_5^{\prime}, f^{\prime}, g^{\prime},
                          v_6^{\prime}]_{m_5^{\prime}}\,,
         & [v_1^{\prime}, v_2^{\prime}, c_1^{\prime}, v_5^{\prime},
                  f^{\prime}, g^{\prime}, v_6^{\prime}]_{m_6^{\prime}}\,,
                                                             \\[.6ex] 
       [v_1^{\prime}, v_2^{\prime}, c_1^{\prime}, c_2^{\prime},
              v_4^{\prime}, f^{\prime}, g^{\prime},
              v_6^{\prime}]_{m_7^{\prime}}\,,
         & [v_1^{\prime}, c_1^{\prime}, c_2^{\prime}, v_5^{\prime},
                          f^{\prime}]_{m_8^{\prime}}\,,      \\[.6ex]
       [v_1^{\prime}, v_2^{\prime}, c_1^{\prime}, c_2^{\prime},
                          v_5^{\prime}]_{m_9^{\prime}}\,,
         & [v_1^{\prime}, c_2^{\prime}, v_4^{\prime}, v_5^{\prime},
                      f^{\prime}]_{m_{10}^{\prime}}\,,     \\[.6ex]
       [v_1^{\prime}, v_4^{\prime}, v_5^{\prime}, f^{\prime},
                          g^{\prime}]_{m_{11}^{\prime}}\,,
         & [v_1^{\prime}, v_4^{\prime}, v_5^{\prime}, g^{\prime},
                          v_6^{\prime}]_{m_{12}^{\prime}}\,,  \\[.6ex]
       [v_1^{\prime}, v_2^{\prime}, v_4^{\prime}, v_5^{\prime},
                      v_6^{\prime}]_{m_{13}^{\prime}}\,,
         & [v_1^{\prime}, v_2^{\prime}, c_2^{\prime}, v_4^{\prime},
                      v_5^{\prime}, ]_{m_{14}^{\prime}} 
      \end{array}
     \right\}\,.
    $$
    {\normalsize Faces} } 
    of higher codimensions are obtained by the intersection
    of subcollections of the above faces.
  \end{quote}
  \vspace{1ex}

  \item
  {\it The triangulation of ${\Delta^{\prime}}^{\ast}$ that gives
       $\Sigma^{\prime}\,$}:
  \begin{quote}
   \hspace{-1.9em}(0)\hspace{1ex}
   In general, triangulations of a polytope by lattice points can be
   generated, using the packages
    "PUNTOS", developed by Loera ([Lo]) and
    "PORTA", developed by Christof ([Ch1-2] and [Zi]).\footnote{
      We would like to thank Shinobu Hosono for providing us with
      these codes and many helpful guidances of their use.}
   In this example, the triangulation we use is given in Example 2.2.2
   in Sec.\ 2.2, following [B-C-dlO-V].
   \vspace{1ex}

   \hspace{-1.9em}(1)\hspace{1ex}
   The vertices $b^{\prime}$, $e_1^{\prime}$, $e_2^{\prime}$, and
   $e_3^{\prime}$ that lie in the boundary of
   ${\Delta^{\prime}}^{\ast}$ are added to the vertex set to give
   a triangulation of ${\Delta^{\prime}}^{\ast}$, the cones 
   - with apex at the origin - over which form the fan $\Sigma^{\prime}$.
   \vspace{1ex}

   \hspace{-1.9em}(2)\hspace{1ex}
   The relation of $b^{\prime}$, $e_1^{\prime}$, $e_2^{\prime}$, and
   $e_3^{\prime}$ with the generating vertices of
   ${\Delta^{\prime}}^{\ast}$ is given by 
   $b^{\prime}=(v_1^{\prime}+4c_1^{\prime}+v_6^{\prime})/6$,
   $e_1^{\prime}=(v_1^{\prime}+c_1^{\prime}+v_6^{\prime})/3$,
   $e_2^{\prime}=(v_1^{\prime}+v_6^{\prime})/2$, and 
   $e_3^{\prime}=(v_1^{\prime}+2v_5^{\prime}+v_6^{\prime})/4\,$.
  \end{quote}
  \vspace{1ex}

  \item
  {\it The singular locus of $X_{\Sigma^{\prime}}$}:
  \begin{quote}
   \hspace{-1.9em}(1)\hspace{1ex}
   The multiplicity $\mult(\sigma^{\prime})$ of
   $\sigma^{\prime}\in\Sigma^{\prime}$ is in the set $\{1, 2, 3\}$.
   Since $2$ and $3$ are primes and $\Sigma^{\prime}$ is
   a simplicial fan, if $\mult(\tau^{\prime})>1$ and
   $\tau^{\prime}\prec\sigma^{\prime}$, then
   $\mult(\tau^{\prime})=\mult(\sigma^{\prime})$.
   \vspace{1ex}

   \hspace{-1.9em}(2)\hspace{1ex}
   Cones in $\Sigma^{\prime}$ whose multiplicity $>1$ but all of
   whose faces have multiplicity $1$ are given by$\,$:
   $$
    \begin{array}{cl}
      \mbox{multiplicity}\;2\;:\;
       & v_5^{\prime}b^{\prime}\,; \\[.6ex]
      \mbox{multiplicity}\;3\;:\;
       & v_4^{\prime}b^{\prime}\,,\;
                   v_4^{\prime}e_1^{\prime}e_2^{\prime}\,.
    \end{array}
   $$
   From the fan structure of $\Sigma^{\prime}$, the three orbit
   closurs $V(v_5^{\prime}b^{\prime})$, $V(v_4^{\prime}b^{\prime})$,
   and $V(v_4^{\prime}e_1^{\prime}e_2^{\prime})$ are disjoint from
   each other. Their union forms the singular locus of
   $X_{\Sigma^{\prime}}$.
   Any other cone with multiplicity $>1$ contains exactly one of
   cones $v_5^{\prime}b^{\prime}$, $v_4^{\prime}b^{\prime}$,
   $v_4^{\prime}e_1^{\prime}e_2^{\prime}$ in its face and hence
   is contained in the corrresponding orbit closure of the latter.
   ({\sc Figure 4-2-1}.)
   \vspace{1ex}

   \hspace{-1.9em}(3)\hspace{1ex}
   The toric geometry of the singular locus and the transverse
   singularity along it is given below$\,$:
   \begin{quote}
    \hspace{-1.9em}(3.1)\hspace{1ex}
    $V(v_5^{\prime}b^{\prime})\,$:
    A toric $3$-variety whose fan is given by $\Sigma$
    (recall the fan of the base).
    The transverse singularity of $X_{\Sigma^{\prime}}$ along the
    open dense orbit of $V(v_5^{\prime}b^{\prime})$ is given by
    the surface $A_1$-singularity.
    \vspace{1ex}

    \hspace{-1.9em}(3.2)\hspace{1ex}
    $V(v_4^{\prime}b^{\prime})\,$:
    A toric $3$-variety whose fan is also given by $\Sigma$.
    The transverse singularity of $X_{\Sigma^{\prime}}$ along the
    open dense orbit of $V(v_5^{\prime}b^{\prime})$ is given by
    the surface $A_2$-singularity.
    \vspace{1ex}

    \hspace{-1.9em}(3.3)\hspace{1ex}
    $V(v_4^{\prime}e_1^{\prime}e_2^{\prime})\,$:
    A toric $2$-variety whose fan is given by
    $$
      \{\,\,(1,0),\,(0,1),\,(-1,0),\,(0,-1)\,\}\,.
    $$
    This is $\CP^1\times\CP^1$.
    The transverse singularity of $X_{\Sigma^{\prime}}$ along the
    open dense orbit of $V(v_4^{\prime}e_1^{\prime}e_2^{\prime})$
    is given by the quotient
    ${\Bbb C}^3/\mbox{\raisebox{-.4ex}{${\Bbb Z}_3$}}$, where
    ${\Bbb Z}_3$ acts on ${\Bbb C}^3$ by
    $(z_1, z_2, z_3)\mapsto(\omega z_1, \omega z_2, \omega z_3)$,
    $\omega^3=1$.
    \vspace{1ex}

    \hspace{-1.9em}(3.4)\hspace{1ex}
    Since $b^{\prime}$, $v_4^{\prime}$, $v_5^{\prime}$ are all mapped
     to zero under $\varphi$, $V(v_5^{\prime}b^{\prime})$ and
     $V(v_4^{\prime}b^{\prime})$ are indeed two sections of the
     fibration $\widetilde{\varphi}$. They correspond to the two
     singular points of the generic fiber $\WCP^2(1,2,3)$
     - one $A_1$-singularity and the other $A_2$-singularity.
    Similarly, since the cone
     $v_4^{\prime}e_1^{\prime}e_2^{\prime}$ is mapped to $r_1$,
     $V(v_4^{\prime}e_1^{\prime}e_2^{\prime})$ sits only over $V(r_1)$.
   \end{quote}
   \medskip
   The structure of these orbit closures are obtained from the fact
   that these cones come from extensions of the $3$-dimensional fan
   $\Sigma^{\prime\prime}$ in Example 2.2.2.
   The transverse singularities follow from [Ba2] and [Fu] and
   the symmetry of the fans involved with respect to the lattice
   structure. Along their lower dimensional orbits, it is the product
   $$
    \mbox{(the above singularity)}\,
    \times{\Bbb C}^{\,\mbox{\scriptsize codimension of that orbit}}\,.
   $$
  \end{quote}
 \end{itemize}
\end{itemize}
\begin{figure}[htbp]
 \setcaption{{\sc Figure 4-2-1.}
  \baselineskip 14pt
  Multiplicity of $[v_5^{\prime}, b^{\prime}]$,
  $[v_4^{\prime},b^{\prime}]$,
  and $[v_4^{\prime}, e_1^{\prime}, e_2^{\prime}]$.
 } 
 \centerline{\psfig{figure=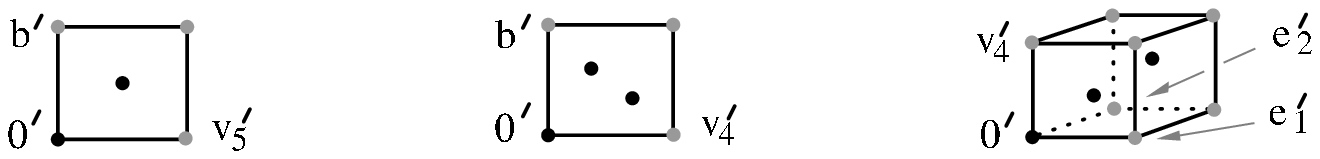,width=11cm}}
\end{figure}

\begin{itemize}
 \item
 {\it The line bundle ${\cal L}_{\Delta^{\prime}}$ and
   its restrictions$\,$}:
 \begin{itemize}
  \item
  ${\cal L}_{\Delta^{\prime}}$ has $3365$ independent holomorphic
  sections.

  \item
  {\it The restriction of ${\cal L}_{\Delta^{\prime}}$ to
  $V(\tau^{\prime})$ for $\tau^{\prime}$ primitive in
  $\Sigma^{\prime}_{\sigma}$ for some $\sigma\in\Sigma\,$}:
  \begin{itemize}
   \item
   {\it the polytope $\Delta^{\prime}_{\tau^{\prime}}$ in
    $\tau^{\prime\perp}_{\scriptsizeBbb R}$, up to translation}$\,$:
    These are the convex hull of a subset of
    $\{\,m_1^{\prime},\,\ldots,\,m_{14}^{\prime}\,\}$ determined by
    $\tau^{\prime}$. They are listed in the $3k+2$ rows of
    {\sc Table 4-2-2}.
    \vspace{.4ex}

   \item
   {\it the lattice points in
    $\Delta^{\prime}_{\tau^{\prime}}\cap\tau^{\prime\perp}$}$\,$:
    Their numbers are listed in the $3k$ rows of {\sc Table 4-2-2}. 

  \end{itemize}

  \item
  {\it The irreducible component $F_{\sigma}^{\tau^{\prime}}$ of
   a generic fiber of
   $\widetilde{\varphi}:V(\tau^{\prime})\rightarrow V(\sigma)$}$\,$:
   These are given in Example 2.2.2, last culumn of {\sc Table 2-2-2}.
   The intersection patterns for these are discussed in that example
   also, cf.\ {\sc Figure 2-2-3}.

  \item
  {\it The restriction ${\cal L}_{F_{\sigma}^{\tau^{\prime}}}$ of
   ${\cal L}$ to $F_{\sigma}^{\tau^{\prime}}\,$}:
  \begin{itemize}
   \item
   {\it the polytope $\overline{\Delta^{\prime}_{\tau^{\prime}}}$ in
    ${\tau^{\prime}}^{\perp}/
      \mbox{\raisebox{-.4ex}{$\varphi^{\dagger}(\sigma^{\perp})$}}\,$}:
    It turns out that, in this example,
    ${\overline{\Delta^{\prime}_{\tau^{\prime}}}}$
    and hence ${\cal L}_{F_{\sigma}^{\tau^{\prime}}}$ depend only
    on the toric geometry of $F_{\sigma}^{\tau^{\prime}}$. They are
    given in {\sc Figure 4-2-3}.
    From the $\overline{\Delta^{\prime}_{\tau^{\prime}}}$ worked out,
    one concludes that except for $F^{\tau^{\prime}}_{\sigma}$ that
    is isomorphic to the generic fiber $\WCP^2(1,2,3)$, the
    restriction ${\cal L}_{F^{\tau^{\prime}}_{\sigma}}$ of
    ${\cal L}_{\Delta^{\prime}}$ to $F^{\tau^{\prime}}_{\sigma}$ is
    not the anticanonical line bundle of $F^{\tau^{\prime}}_{\sigma}$.
    One also observe that the various
    $\overline{\Delta^{\prime}_{\tau^{\prime}}}$ all come from
    some chopping of the polytope
    $\overline{\Delta^{\prime}_{\{0\}}}$ associated to the generic
    fiber $\WCP^2(1,2,3)$.

   \item
   {\it The lattice points of the polytope
    $\overline{\Delta^{\prime}_{\tau^{\prime}}}$ in
    ${\tau^{\prime}}^{\perp}/
     \mbox{\raisebox{-.4ex}{$\varphi^{\dagger}(\sigma^{\perp})$}}\,$}:
    The polytopes are indicated in {\sc Figure 4-2-3}.
    The numbers of lattice points they contain are listed below$\,$:
    \newline

    \vspace{1ex}
    \hspace{1em}
    {\footnotesize
    \hspace{-3em}
    \begin{tabular}{|ccccccccccc|}  \hline 
     \rule{0ex}{3ex}
      $\mbox{\rm W{\footnotesizeBbb C}P}^2(1,2,3)$  && $X(4)$
      && $\mbox{\rm {\footnotesizeBbb C}P}^2$
      && $X(5)$       && $\mbox{\rm W{\footnotesizeBbb C}P}^2(1,1,3)$
      && ${\footnotesizeBbb F}_2$                     \\[.6ex]\hline
     \rule{0ex}{3ex}
      $7$      && $4$     && $6$
      && $2$   && $5$    && $4$ \\[.6ex]\hline
    \end{tabular}  } 
  \end{itemize}
 \end{itemize}
\end{itemize}

\bigskip

\begin{minipage}{11cm}
{
 \tiny
 \centerline{
 \begin{tabular}{|cccc|}  \hline
  \rule{0ex}{3ex}
   $O^{\prime}\;\; (54)$
     & $v_1^{\prime}\;\;(27)$
       & $v_2^{\prime}\;\;(24)$
         & $c_1^{\prime}\;\;(16)$
           \\[.6ex] \hline
  \rule{0ex}{3ex}
   all of $m^{\prime}_i$
     & $[6,7,8,9,10,11,12,13,14]$
       & $[2,4,6,7,9,13,14]$
         & $[3,4,6,7,8,9]$
           \\[.6ex] \hline
  \rule{0ex}{3ex}
   $3365$    & $227$   & $262$   & $154$ 
           \\[.6ex] \hline\hline
  \rule{0ex}{3ex}
    $c_2^{\prime}\;\;(12)$
      & $e_1^{\prime}\;\;(20)$
        & $e_2^{\prime}\;\;(12)$
          & $e_3^{\prime}\;\;(16)$
            \\[.6ex] \hline
  \rule{0ex}{3ex}
   $[1,2,3,4,7,8,9,10,14]$
    & $[6,7]$                
     & $[6, 7, 12, 13]$
      & $[6, 12, 13]$
          \\[.6ex] \hline
  \rule{0ex}{3ex}
    $1242$  & $2$   & $11$    & $10$  
     \\[.6ex] \hline\hline
  \rule{0ex}{3ex}
   $f^{\prime}\;\;(16)$
         & $g^{\prime}\;\;(20)$
           & $v_6^{\prime}\;\;(27)$
             & $f^{\prime}c_1^{\prime}\;\;(8)$
      \\[.6ex] \hline
  \rule{0ex}{3ex}  
   $[1, 3, 5, 6, 7, 8, 10, 11]$
         & $[5, 6, 7, 11, 12]$
           & $[1, 2, 3, 4, 5, 6, 7, 12, 13]$
             & $[3, 6, 7, 8]$ 
      \\[.6ex] \hline
   \rule{0ex}{3ex} 
   $237$   & $64$   & $227$   & $22$
                                        \\[.6ex] \hline\hline
  \rule{0ex}{3ex}
   $f^{\prime}c_2^{\prime}\;\;(6)$
     & $f^{\prime}g^{\prime}\;\;(6)$
       & $g^{\prime}e_1^{\prime}\;\;(10)$
         & $g^{\prime}e_2^{\prime}\;\;(6)$
               \\[.6ex] \hline
  \rule{0ex}{3ex}
   $[3, 7, 8, 10]$
     & $[5, 6, 7, 11]$
       & $[6, 7]$
         & $[6, 7, 12]$
               \\[.6ex] \hline
  \rule{0ex}{3ex}
   $67$    & $39$   & $2$   & $5$
          \\[.6ex] \hline\hline
  \rule{0ex}{3ex}
   $g^{\prime}e_3^{\prime}\;\;(8)$
     & $v_2^{\prime}e_1^{\prime}\;\;(10)$
      &   $v_2^{\prime}e_2^{\prime}\;\;(6)$
       & $v_2^{\prime}e_3^{\prime}\;\;(8)$
         \\[.6ex] \hline
  \rule{0ex}{3ex}
   $[6, 12]$
     & $[6, 7]$ 
      &   $[6, 7, 13]$
        & $[6, 13]$
                       \\[.6ex] \hline
  \rule{0ex}{3ex}
    $4$  & $2$  &  $5$    & $4$  \\[.6ex] \hline\hline
  \rule{0ex}{3ex}
    $v_2^{\prime}c_1^{\prime}\;\;(8)$
         & $v_2^{\prime}c_2^{\prime}\;\;(6)$
           & $v_1^{\prime}v_2^{\prime}\;\;(12)$
             & $v_1^{\prime}c_1^{\prime}\;\;(8)$
                \\[.6ex] \hline 
  \rule{0ex}{3ex}
   $[4, 6, 7, 9]$
         & $[2, 4, 7, 9, 14]$
           & $[6, 7, 9, 13, 14]$
             & $[6, 7, 8, 9]$
                \\[.6ex] \hline
  \rule{0ex}{3ex}
   $22$   & $86$   &  $26$  & $16$
                                        \\[.6ex] \hline\hline
  \rule{0ex}{3ex}
   $v_1^{\prime}c_2^{\prime}\;\;(6)$
     & $v_1^{\prime}e_1^{\prime}\;\;(10)$
       & $v_1^{\prime}e_2^{\prime}\;\;(6)$
         & $v_1^{\prime}e_3^{\prime}\;\;(8)$
                                       \\[.6ex] \hline
  \rule{0ex}{3ex}
   $[7, 8, 9, 10, 14]$
     & $[6, 7]$
       & $[6, 7, 12, 13]$
         & $[6, 12, 13]$
                     \\[.6ex] \hline
  \rule{0ex}{3ex}
   $62$    & $2$   & $11$   & $10$ 
                                        \\[.6ex] \hline\hline
  \rule{0ex}{3ex}
   $v_1^{\prime}f^{\prime}\;\;(8)$
     & $v_1^{\prime}g^{\prime}\;\;(10)$
        & $v_6^{\prime}v_2^{\prime}\;\;(12)$
           & $v_6^{\prime}c_1^{\prime}\;\;(8)$
              \\[.6ex] \hline
  \rule{0ex}{3ex}
    $[6, 7, 8, 10, 11]$
      & $[6, 7, 11, 12]$
        &  $[2, 4, 6, 7, 13]$
           & $[3, 4, 6, 7]$
                    \\[.6ex] \hline
  \rule{0ex}{3ex}
   $19$  & $14$  &  $26$    & $16$  
                   \\[.6ex] \hline\hline
  \rule{0ex}{3ex}
     $v_6^{\prime}c_2^{\prime}\;\;(6)$
         & $v_6^{\prime}e_1^{\prime}\;\;(10)$
           & $v_6^{\prime}e_2^{\prime}\;\;(6)$
             & $v_6^{\prime}e_3^{\prime}\;\;(8)$
                 \\[.6ex] \hline
  \rule{0ex}{3ex}
     $[1, 2, 3, 4, 7]$
         & $[6, 7]$
           & $[6, 7, 12, 13]$
             & $[6, 12, 13]$   \\[.6ex] \hline 
  \rule{0ex}{3ex}
    $62$   & $2$   &  $11$  & $10$  \\[.6ex] \hline\hline
  \rule{0ex}{3ex}
   $v_6^{\prime}f^{\prime}\;\;(8)$
     & $v_6^{\prime}g^{\prime}\;\;(10)$
       & $v_1^{\prime}f^{\prime}c_1^{\prime}\;\;(4)$
         & $v_1^{\prime}f^{\prime}c_2^{\prime}\;\;(3)$
                  \\[.6ex] \hline
  \rule{0ex}{3ex}
   $[1, 3, 5, 6, 7]$
     & $[5, 6, 7, 12]$
       & $[6, 7, 8]$
         & $[7, 8, 10]$
                 \\[.6ex] \hline
  \rule{0ex}{3ex}
   $19$   & $14$   & $4$   &  $6$ 
           \\[.6ex] \hline\hline
  \rule{0ex}{3ex}
   $v_1^{\prime}f^{\prime}g^{\prime}\;\;(3)$
     & $v_1^{\prime}g^{\prime}e_1^{\prime}\;\;(5)$
       & $v_1^{\prime}g^{\prime}e_2^{\prime}\;\;(3)$
           & $v_1^{\prime}g^{\prime}e_3^{\prime}\;\;(4)$
              \\[.6ex] \hline
  \rule{0ex}{3ex}
   $[6, 7, 11]$
     & $[6, 7]$     
       & $[6, 7, 12]$
          & $[6, 12]$
                       \\[.6ex] \hline
  \rule{0ex}{3ex}
    $7$  & $2$  &  $5$    & $4$ 
      \\[.6ex] \hline\hline
  \rule{0ex}{3ex}  
     $v_1^{\prime}v_2^{\prime}e_1^{\prime}\;\;(5)$
         & $v_1^{\prime}v_2^{\prime}e_2^{\prime}\;\;(3)$
           & $v_1^{\prime}v_2^{\prime}e_3^{\prime}\;\;(4)$
             & $v_1^{\prime}v_2^{\prime}c_1^{\prime}\;\;(4)$
       \\[.6ex] \hline  
  \rule{0ex}{3ex}
       $[6, 7]$
         & $[6, 7, 13]$
           & $[6, 13]$
             & $[6, 7, 9]$                   \\[.6ex] \hline
  \rule{0ex}{3ex}  
     $2$   & $5$   &  $4$  & $4$
                                        \\[.6ex] \hline\hline
  \rule{0ex}{3ex}
   $v_1^{\prime}v_2^{\prime}c_2^{\prime}\;\;(3)$
     & $v_6^{\prime}f^{\prime}c_1^{\prime}\;\;(4)$
       & $v_6^{\prime}f^{\prime}c_2^{\prime}\;\;(3)$
         & $v_6^{\prime}f^{\prime}g^{\prime}\;\;(3)$
             \\[.6ex] \hline
  \rule{0ex}{3ex}
   $[7, 9, 14]$
     & $[3, 6, 7]$
       & $[1, 3, 7]$
         & $[5, 6, 7]$
                         \\[.6ex] \hline
  \rule{0ex}{3ex}
   $6$    & $4$   & $6$   & $7$    
                         \\[.6ex] \hline\hline
  \rule{0ex}{3ex}
    $v_6^{\prime}g^{\prime}e_1^{\prime}\;\;(5)$
       & $v_6^{\prime}g^{\prime}e_2^{\prime}\;\;(3)$
          &  $v_6^{\prime}g^{\prime}e_3^{\prime}\;\;(4)$
              & $v_6^{\prime}v_2^{\prime}e_1^{\prime}\;\;(5)$
                   \\[.6ex] \hline
  \rule{0ex}{3ex}
    $[6, 7]$
     & $[6, 7, 12]$                 
        & $[6, 12]$
            & $[6, 7]$
                      \\[.6ex] \hline
  \rule{0ex}{3ex}
    $2$  & $5$   & $4$    & $2$ 
     \\[.6ex] \hline\hline
  \rule{0ex}{3ex}
    $v_6^{\prime}v_2^{\prime}e_2^{\prime}\;\;(3)$
         & $v_6^{\prime}v_2^{\prime}e_3^{\prime}\;\;(4)$
           & $v_6^{\prime}v_2^{\prime}c_1^{\prime}\;\;(4)$
             & $v_6^{\prime}v_2^{\prime}c_2^{\prime}\;\;(3)$
              \\[.6ex] \hline
  \rule{0ex}{3ex}
    $[6, 7, 13]$
         & $[6, 13]$
           & $[4, 6, 7]$
             & $[2, 4, 7]$  
                  \\[.6ex] \hline
  \rule{0ex}{3ex}
    $5$   & $4$   &  $4$  & $6$
                                        \\[.6ex] \hline
 \end{tabular}
 }  
 }  

 \bigskip

 \centerline{
   \parbox{9cm}{ {\sc Table 4-2-2.}
   Primitive cones $\tau^{\prime}$ with respect to $\varphi$ are
   listed in the $3k+1$ rows.
   The toric combinatorial complexity of $V(\tau^{\prime})$ is
   indicated by the number of the maximal cones of its fan
   (in parenthesis).
   The generating vertices for the polytope
   $\Delta^{\prime}_{\tau^{\prime}}$ are indicated in the
   $3k+2$ rows and the number of lattice points contained in it
   on the $3k$ rows.
   }
 } 
\end{minipage}

\newpage
 \begin{figure}[htbp]
  \setcaption{{\sc Figure 4-2-3.}
   \baselineskip 14pt
   Tha fan for the six different types of irreducible components
   $F^{\tau^{\prime}}_{\sigma}$ of fibers of $\widetilde{\varphi}$
   are listed in the first row.
   The restriction ${\cal L}_{F^{\tau^{\prime}}_{\sigma}}$
   of ${\cal L}_{\Delta^{\prime}}$ to these fibers are described
   by the polytope $\overline{\Delta^{\prime}_{\tau^{\prime}}}$
   indicated in the second row. Observe that the fan for
   $F^{\tau^{\prime}}_{\sigma}$ is a refinement of the normal fan
   of $\overline{\Delta^{\prime}_{\tau^{\prime}}}$ and that all the
   polytopes $\overline{\Delta^{\prime}_{\tau^{\prime}}}$ are some
   chopping of $\overline{\Delta^{\prime}_{\{0\}}}$.
  } 
  \centerline{\psfig{figure=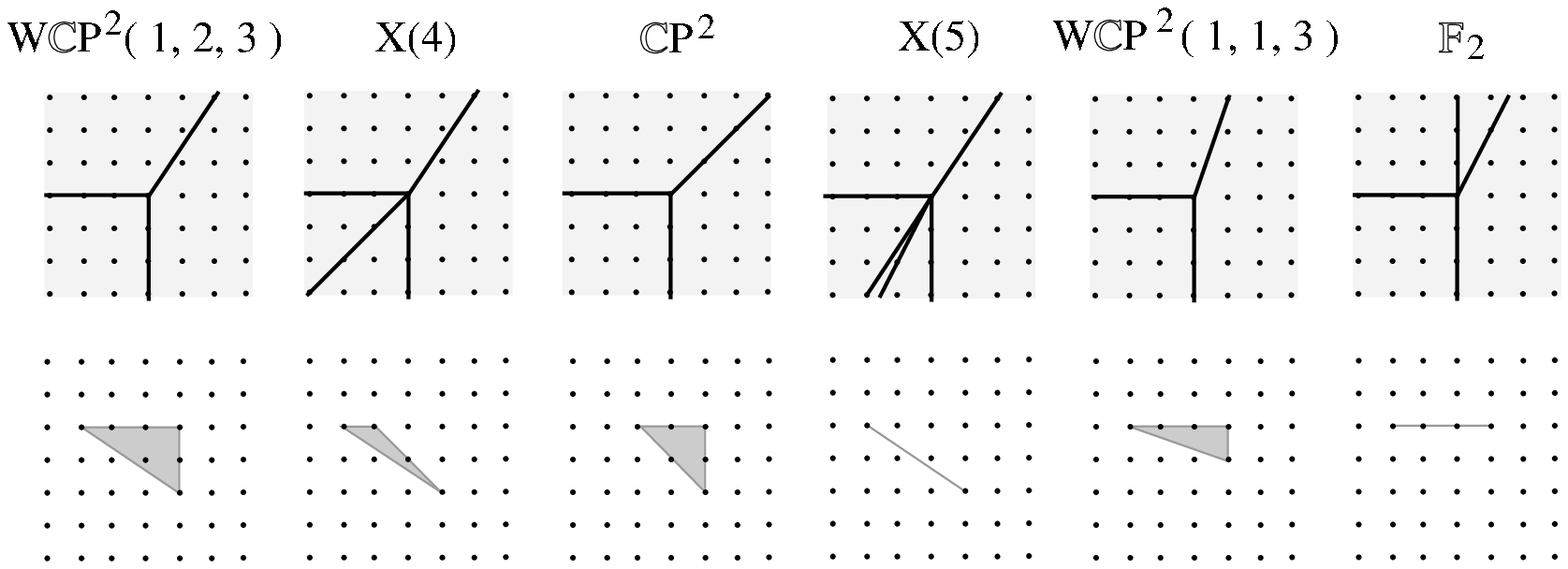,width=11cm,caption=}}
 \end{figure}

\begin{itemize}
\item
 {\it
 The discriminant associated to
 $\overline{\Delta^{\prime}_{\tau^{\prime}}}$ in
 ${\tau^{\prime}}^{\perp}/
      \mbox{\raisebox{-.4ex}{$\varphi^{\dagger}(\sigma^{\perp})$}}\,$}:
 They are related to the discriminant locus of the family of sections
 of the line bundle
 ${\cal L}_{\overline{\Delta^{\prime}_{\tau^{\prime}}}}$ over
 $F_{\sigma}^{\tau^{\prime}}$.
 \begin{quote}
  \hspace{-1.9em}(1)\hspace{1ex}
  For $\overline{\Delta^{\prime}_{\tau^{\prime}}}$ associated
  to $\WCP^2(1,2,3)$:
  After the integral affine transformation on $M^{\prime}$ such
  that $(1,1)\rightarrow (0,0)$, $(-2,1)\rightarrow (3,0)$,
  $(1,-1)\rightarrow (0,2)$, the local polynomial form of
  a section becomes
  $$
   f(x,y)\;
   =\;a_{00}+a_{10}x+a_{20}x^2+a_{30}x^3+a_{01}y+a_{11}xy+a_{02}y^2\,.
  $$
  Its discriminant is given by ([G-K-Z3] and [La]) 
  {\footnotesize
  \begin{eqnarray*}
   \lefteqn{
    \delta_{\overline{\Delta^{\prime}_{\tau^{\prime}}}}\;
    =\; -27\cdot 16\,a_{00}^2a_{02}^3a_{30}^2\,
        -\,64\,a_{00}a_{20}^3a_{02}^3\,
        -\, 64\,a_{10}^3a_{02}^3a_{30}\,
        -\, 27\,a_{01}^4a_{02}a_{30}^2 } \\
     && +\, a_{00}a_{11}^6\,
        +\, 16\, a_{10}^2a_{20}^2a_{02}^3\,
        +\,16\, a_{01}^2a_{20}^3a_{02}^2\,
        +\, a_{02}a_{10}^2a_{11}^4\,
        -\,a_{01}a_{10}a_{11}^5         \\
     && +\, a_{01}^2a_{11}^4a_{20}\,
        -\, a_{01}^3a_{11}^3a_{30}\,
        +\, 32\,\cdot 9\, a_{00} a_{02}^3a_{10}a_{20}a_{30} \\
     && +\, 48\,a_{00} a_{02}^2a_{11}^2a_{20}^2\,
        +\, 8\,\cdot 27\, a_{00}a_{01}^2a_{02}^2a_{30}^2\,,
        -\,72\,a_{01}^2a_{02}^2a_{10}a_{20}a_{30}      \\
     && -\,72\,a_{00}a_{02}^2a_{10}a_{11}^2a_{30}\,         
        -\,16\,a_{01}a_{02}^2a_{10}a_{11}a_{20}^2\,
        -\,8\,a_{02}^2a_{10}^2a_{11}^2a_{20}\,               \\
     && +\, 96\,a_{01}a_{02}^2a_{10}^2a_{11}a_{30}\,
        -\, 144\,a_{00}a_{01}a_{02}^2a_{11}a_{20}a_{30}\,  
        -\, 12\,a_{00}a_{02}a_{11}^4a_{20}    \\
     && +\, 8\,a_{01}a_{02}a_{10}a_{11}^2a_{20}\,   
        -\, 8\,a_{01}^2a_{02}a_{11}^2a_{20}^2\,
        -\, 30\,a_{01}^2a_{02}a_{10}a_{11}^2a_{30}  \\  
     && +\, 36\,a_{01}^3a_{02}a_{11}a_{20}a_{30}\,
        +\, 36\,a_{00}a_{01}a_{02}a_{11}^3a_{30}\,.
  \end{eqnarray*} } 

  \hspace{-1.9em}(2)\hspace{1ex}
  For $\overline{\Delta^{\prime}_{\tau^{\prime}}}$ associated
  to $X(4)$:
  After the integral affine transformation on $M^{\prime}$ such
  that $(-1,1)\rightarrow (0,0)$, $(-2,1)\rightarrow (0,1)$,
  $(-2,1)\rightarrow (2,0)$, the local polynomial form of a section
  becomes $f(x,y)\,=\,a_{00}+a_{10}x+a_{20}x^2+a_{01}y$.
  Its discriminant is given by ([G-K-Z3]) 
  $\delta_{\overline{\Delta^{\prime}_{\tau^{\prime}}}}\,=\,1$.
 
  \hspace{-1.9em}(3)\hspace{1ex}
  For $\overline{\Delta^{\prime}_{\tau^{\prime}}}$ associated
  to $\CP^2$:
  After the integral affine transformation on $M^{\prime}$ such
  that $(1,1)\rightarrow (0,0)$, $(-1,1)\rightarrow (2,0)$,
  $(1,-1)\rightarrow (0,2)$, the local polynomial form of
  a section becomes
  $$
   f(x,y)\,=\,a_{00}+a_{10}x+a_{20}x^2+a_{01}y+a_{11}xy+a_{02}y^2\,.
  $$
  Its discriminant is ([G-K-Z3])
  $$
   \delta_{\overline{\Delta^{\prime}_{\tau^{\prime}}}}\,
   =\,a_{20}a_{01}^2\,+\,a_{10}^2a_{02}\,+\,a_{11}^2a_{00}\,
       -\,a_{11}a_{10}a_{01}\,-\,4\,a_{20}a_{02}a_{00}\,.
  $$

  \hspace{-1.9em}(4)\hspace{1ex}
  For $\overline{\Delta^{\prime}_{\tau^{\prime}}}$ associated
  to $X(5)$:
  After an integral affine transformation on $M^{\prime}$ such
  that $(-2,1)\rightarrow (0,0)$, $(1,-1)\rightarrow (1,0)$,
  the local polynomial form of a section becomes
  $f(x)\,=\,a_0+a_1x$. Its discriminant is given by ([G-K-Z3])
  $\delta_{\overline{\Delta^{\prime}_{\tau^{\prime}}}}\,=\,a_1$.

  \hspace{-1.9em}(5)\hspace{1ex}
  For $\overline{\Delta^{\prime}_{\tau^{\prime}}}$ associated
  to $\WCP^2(1,1,3)$:
  After the integral affine transformation on $M^{\prime}$ such
  that $(1,1)\rightarrow (0,0)$, $(-2,1)\rightarrow (3,0)$,
  $(1,0)\rightarrow (0,1)$, the local polynomial form of
  a section becomes
  $$
   f(x,y)\,=\,a_{00}+a_{10}x+a_{20}x^2+a_{30}x^3+a_{01}y\,.
  $$
  Its discriminant is given by ([G-K-Z3])
  $\delta_{\overline{\Delta^{\prime}_{\tau^{\prime}}}}\,=\,1$.

  \hspace{-1.9em}(6)\hspace{1ex}
  For $\overline{\Delta^{\prime}_{\tau^{\prime}}}$ associated
  to ${\Bbb F}_2$:
  After an integral affine transformation on $M^{\prime}$ such
  that $(-2,1)\rightarrow (0,0)$, $(1,1)\rightarrow (3,0)$,
  the local polynomial form of a section becomes
  $f(x)\,=\,a_0+a_1x+a_2x^2+a_3x^3$.
  Its discriminant is given by ([G-K-Z3]) 
  $\delta_{\overline{\Delta^{\prime}_{\tau^{\prime}}}}\,
    =\, 27\,a_0^2a_3^2\,+\,4\,a_0a_2^3\,+\,4\,a_1^3a_3\,
        -\,a_1^2a_2^2\,-\,18\,a_0a_1a_2a_3$.
 \end{quote}
\end{itemize}

\begin{itemize}
\item
 {\it The $\Sigma^{\prime}$-regular sections in
  ${\cal L}_{\Delta^{\prime}}$ and their restrictions}$\,$:
 \begin{itemize}
  \item
  {\it A $\Sigma^{\prime}$-regular section $s$ of
   ${\cal L}_{\Delta^{\prime}}$
   and its zero-locus, the Calabi-Yau hypersurface $Y^{\prime}\,$}:
   The Hodge numbers of these Calabi-Yau hypersurfaces were computed
   in [B-C-dlO-G]$\,$:
   $h^{1,1}\,=\,8$, $h^{3,1}\,=\,2897$, $h^{2,1}\,=\,1$, and
   $h^{2,2}\,=\,11662$.
   The dimension of the polynomial complex moduli space
   ${\cal M}_{\mbox{\rm\scriptsize poly}}$ is computed by the formula
   ([Ba2] and [C-K])
   $$
    \dimm({\cal M}_{\rm\scriptsize poly})\;
     =\;|\Delta^{\prime}\cap M^{\prime}|-6-
         \sum_{\Theta^{\prime}}\,l(\Theta^{\prime})\;
     =\; 2897\,,
   $$
   where $\Theta^{\prime}$ are the codimension-$1$ faces of
   $\Delta^{\prime}$ and $l(\Theta^{\prime})$ is the number of
   lattice points in the relative interior of $\Theta^{\prime}$.
   Thus, the zero-locus of the $3365$-parameter family of
   sections of ${\cal L}_{\Delta^{\prime}}$ gives
   a $2897$-parameter family of Calabi-Yau $4$-varieties.

   \medskip
   {\it Remark.} The triangulation of ${\Delta^{\prime}}^{\ast}$ used
   so far missed $4$ points in ${\Delta^{\prime}}^{\ast}\cap N^{\prime}$.
   The latter come in when $X_{\Sigma^{\prime}}$ and hence $Y^{\prime}$
   are desingularized. The above counting of moduli is more related
   to what comes after desingularization.

  \item
  {\it The singular loci of $Y^{\prime}\,$}:
  Under the assumption that $Y^{\prime}$ is $\Sigma^{\prime}$-regular,
  the singular locus $\Sing(Y^{\prime})$ of $Y^{\prime}$ are induced
  from that of the ambient $X_{\Sigma^{\prime}}$. Thus,
  $\Sing(Y^{\prime})
    =Y^{\prime}_{[v_5^{\prime}b^{\prime}]}
     \sqcup Y^{\prime}_{[v_4^{\prime}b^{\prime}]}
     \sqcup Y^{\prime}_{[v_4^{\prime}e_1^{\prime}e_2^{\prime}]}$.
  The first two components, if non-empty, have co-dimension $2$
  in $Y^{\prime}$ while the last component, if nonempty, has
  co-dimension $3$.
  
  \item
  {\it The restriction of $s$ to $V(\tau^{\prime})$ for
   $\tau^{\prime}$ primitive in $\Sigma^{\prime}_{\sigma}$
   for some $\sigma\in\Sigma\,$}:
  Suppose $s$ is given by
  $$
   s\;=\;\sum_{m^{\prime}\in\Delta^{\prime}\cap M^{\prime}}\,
          a_{m^{\prime}}\chi^{m^{\prime}}\,,
    \hspace{1em}\mbox{(a total of $3365$ summands)}\,.
  $$
  Then
  $$
   s|_{V(\tau^{\prime})}\;
    =\;\sum_{m^{\prime}
              \in\Delta^{\prime}_{\tau^{\prime}}\cap M^{\prime}}\,
        a_{m^{\prime}} \chi^{m^{\prime}}\,,
     \hspace{1em}\mbox{\parbox[t]{4.6cm}{(only summands related to lattice
            points in $\Delta^{\prime}_{\tau^{\prime}}$ remain).}}
  $$
  The character $\chi^{m^{\prime}}$,
  $m^{\prime}\in\Delta^{\prime}_{\tau^{\prime}}$, is then rewritten
  as a character related to the lower dimensional torus associated
  to $V(\tau^{\prime})$ by projecting it to ${\tau^{\prime}}^{\perp}$.
  Let us give a concrete example.

  \bigskip
  {\bf Demonstration for $[v_1^{\prime}, e_2^{\prime}]$.}
  Let $\tau^{\prime}=[v_1^{\prime},e_2^{\prime}]$, then
  $$
   \Delta^{\prime}_{[v_1^{\prime},e_2^{\prime}]}\;
    =\;[\,m_6^{\prime},\, m_7^{\prime},\, m_{12}^{\prime},\,
                                                m_{13}^{\prime}\,]
  $$
  (cf.\ {\sc Table 4-2-2}) and
  {\footnotesize
  $$
   \Delta^{\prime}_{[v_1^{\prime},\,e_2^{\prime}]}\cap M^{\prime}\;
    =\; \left\{\,
         \begin{array}{llll}
          (0,0,0,-2,1)\,,  & (0,0,0,1,-1)\,, & (2,1,-1,-1,1)\,, \\
          (2,2,-1,-1,1)\,, & (4,2,-2,0,1)\,,  & (4,3,-2,0,1)\,, \\
          (4,4,-2,0,1)\,,  & (6,3,-3,1,1)\,,  & (6,4,-3,1,1)\,, \\
          (6,5,-3,1,1)\,,    & (6,6,-3,1,1)
         \end{array}\,\right\}\,.
  $$
  {\normalsize Fix}} 
  a basis, e.g.\ $((2,0,-1,1,0),(0,1,0,0,0),(6,0,-3,0,2))$,
  for $[v_1^{\prime},e_2^{\prime}]^{\perp}\cap M^{\prime}$
  and take, say $(0,0,0,-2,1)$, as the origin. Then
  {\footnotesize
  $$
   \Delta^{\prime}_{[v_1^{\prime},\,e_2^{\prime}]}\cap M^{\prime}\;
    =\; \left\{\,
         \begin{array}{llll}
          (0,0,0)\,,  & (3,0,-1)\,, & (1,1,0)\,,  & (1,2,0) \\
          (2,2,0)\,,  & (2,3,0)\,,  & (2,4,0)\,,  & \\
          (3,3,0)\,,  & (3,4,0)\,,  & (3,5,0)\,,  & (3,6,0)
         \end{array}\,\right\}\,.
  $$
  {\normalsize and}} 
  {\footnotesize
  \begin{eqnarray*}
   \lefteqn{
     s|_{V([v_1^{\prime},\,e_2^{\prime}])}\;
       =\; a_{(0,0,0)}\,+\,a_{(3,0,-1)}t_1^3t_3^{-1}\,
        +\,a_{(1,1,0)}t_1t_2\,+\,a_{(1,2,0)}t_1t_2^2 } \\
    &&  +\, a_{(2,2,0)}t_1^2t_2^2\, +\,a_{(2,3,0)}t_1^2t_2^3\,
        +\,a_{(2,4,0)}t_1^2t_2^4\,+\,a_{(3,3,0)}t_1^3t_2^3\,
        +\,a_{(3,4,0)}t_1^3t_2^4 \\   
    &&  +\,a_{(3,5,0)}t_1^3t_2^5\, +\,a_{(3,6,0)}t_1^3t_2^6\,.
  \end{eqnarray*}
  {\normalsize To}} 
  understand the meaning of this, fix a transverse dual basis, say
  $(0,0,0,1,0)$, $(0,1,0,0,0)$, $(0,0,-1,-1,-1))$,
  to the cone $[v_1^{\prime},e_2^{\prime}]$ in $N^{\prime}$.
  The $6$ maximal cones (cf.\ Table 4-2-2) in $\Sigma^{\prime}$ that
  contains $[v_1^{\prime},e_2^{\prime}]$ in the face are
  {\small
   $$
   v_1^{\prime}e_3^{\prime}e_1^{\prime}g^{\prime}e_2^{\prime},\;\;
   v_1^{\prime}e_3^{\prime}v_4^{\prime}g^{\prime}e_2^{\prime},\;\;
   v_1^{\prime}e_1^{\prime}v_4^{\prime}g^{\prime}e_2^{\prime},\;\;
   v_1^{\prime}e_3^{\prime}e_1^{\prime}v_2^{\prime}e_2^{\prime},\;\;
   v_1^{\prime}e_3^{\prime}v_2^{\prime}v_4^{\prime}e_2^{\prime},\;\;
   v_1^{\prime}e_1^{\prime}v_2^{\prime}v_4^{\prime}e_2^{\prime}\,.
  $$
  {\normalsize They}} 
  project to the following six maximal cones in
  $\Star([v_1^{\prime},\,e_2^{\prime}])$:
  $$
   \bar{e}_3^{\prime}\bar{e}_1^{\prime}\bar{g}^{\prime},\;\;
   \bar{e}_3^{\prime}\bar{v}_4^{\prime}\bar{g}^{\prime},\;\;
   \bar{e}_1^{\prime}\bar{v}_4^{\prime}\bar{g}^{\prime},\;\;
   \bar{e}_3^{\prime}\bar{e}_1^{\prime}\bar{v}_2^{\prime},\;\;
   \bar{e}_3^{\prime}\bar{v}_2^{\prime}\bar{v}_4^{\prime},\;\;
   \bar{e}_1^{\prime}\bar{v}_2^{\prime}\bar{v}_4^{\prime}\,.
  $$
  where
  {\small 
   $$
    \bar{e}_1^{\prime}=(1,0,3),\;  \bar{e}_3^{\prime}=(0,0,-1),\;
    \bar{v}_2^{\prime}=(2,-1,6),\; \bar{v}_4^{\prime}=(-1,0,0),\;
    \bar{g}^{\prime}=(-1,1,-3)\,.
   $$
  {\normalsize In}} 
  this way, ${\cal L}|_{V([v_1^{\prime},e_2^{\prime}])}$
  is described directly by the $3$-dimensional polytope$\,$:
  the convex hull of $(0,0,0)$, $(3,0,-1)$, $(3,3,0)$, and $(3,6,0)$
  and $s|_{V([v_1^{\prime}e_2^{\prime}])}$ as given is a combination
  of characters on the $3$-dimensional torus $({\Bbb C}^{\times})^3$
  associated to the cone: the origin of the rank $3$ lattice in
  describing $\Star([v_1^{\prime},\,e_2^{\prime}])$.  
  \newline$\mbox{\hspace{10cm}}\Box$
  \bigskip

  \item
  {\it The restriction of $s$ to $F_{\sigma}^{\tau^{\prime}}\,$}:
  \begin{itemize}
   \item
   {\it description of the irreducible components of a generic fiber
    of $\widetilde{\varphi}:Y^{\prime}\rightarrow X_{\Sigma}$
    over each of the $33$ toric orbits of $X_{\Sigma}\,$}:
   \begin{quote}
    \hspace{-1.9em}(1)\hspace{1ex}
    The essential feature depends only on the toric geometry of the
    fiber $F_{\sigma}$.

    \bigskip
    \hspace{-1.9em}(2)\hspace{1ex}
    {\it For $F_{\sigma}=X(4)\cup\CP^2\,$}:
    With labels in {\sc Figure 2-2-3},
    {\footnotesize
    $$
     \begin{array}{lcrll}
      {\cal L}_{c_1^{\prime}}={\cal O}(D_1+D_2+D_4) & \mbox{and}
        & (D_1+D_2+D_4)\cdot D_3=2 & \mbox{on $X(4)$}\,;  \\
      {\cal L}_{c_2^{\prime}}={\cal O}(D_5+D_6)     & \mbox{and}
        & (D_5+D_6)\cdot D_3=2  
        & \mbox{on $\mbox{\rm $\footnotesizeBbb C$P}^2$}\,.
    \end{array}
    $$
    {\normalsize These}} 
    data restrict the generic degenerate fiber
    of $Y^{\prime}$ that lies in $X(4)\cup\CP^2$.
    (To be continued in Item: Beyond.)

    \bigskip
    \hspace{-1.9em}(3)\hspace{1ex}
    {\it For $F_{\sigma}=X(5)\cup\WCP^2(1,1,3)\cup{\Bbb F}_2\,$}:
    With labels in {\sc Figure 2-2-3},
    {\footnotesize
    $$
    \hspace{-3em}
    \begin{array}{lcl}
     {\cal L}_{e_1^{\prime}}={\cal O}(D_1+D_2-D_3+D_5) & \mbox{and}
       & \left\{
          \begin{array}{rr}
           (D_1+D_2-D_3+D_5)\cdot D_3=4 \\
           (D_1+D_2-D_3+D_5)\cdot D_4=0
          \end{array}
         \right.      \\[.6ex]
     & & \mbox{on $X(5)$}\,;    \\[2ex]
     {\cal L}_{e_2^{\prime}}={\cal O}(-D_3+D_6+D_7)    & \mbox{and}
       & \left\{
          \begin{array}{rr}
           (-D_3+D_6+D_7)\cdot D_3=1 \\
           (-D_3+D_6+D_7)\cdot D_7=3
          \end{array}
         \right.   \\[.6ex]                  
     & & \mbox{on $\mbox{\rm W${\footnotesizeBbb C}$P}^2(1,1,3)$}\,;
                                                            \\[2ex]
     {\cal L}_{e_3^{\prime}}={\cal O}(-D_7+D_8+D_9)    & \mbox{and}
       & \left\{
          \begin{array}{rr}
           (-D_7+D_8+D_9)\cdot D_4=0 \\
           (-D_7+D_8+D_9)\cdot D_7=3
          \end{array}
         \right.   \\[.6ex]
     & & \mbox{on ${\footnotesizeBbb F}_2$}\,.
    \end{array}
    $$
    {\normalsize These}} 
    data restrict the generic degenerate fiber
    of $Y^{\prime}$ that lies in $X(5)\cup\WCP^2(1,1,3)\cup{\Bbb F}_2$.
    (To be continued in Item: Beyond.)
   \end{quote} 

   \medskip
   \item
   {\it restriction of a section of ${\cal L}_{\Delta^{\prime}}$ to
    $V(\tau^{\prime})$ and then rewrite it in the fibred form
    with respect to a map
    $\xi_{\sigma}^{\tau^{\prime}}:
      N/\mbox{\raisebox{-.4ex}{$N_{\sigma}$}} \rightarrow
      N^{\prime}/\mbox{\raisebox{-.4ex}{$N^{\prime}_{\tau^{\prime}}$}}\,$}:
    The details are already explained in the theory part. Let us give
    an explicit computation here for $[v_1^{\prime},e_2^{\prime}]$
    as an example.

   \bigskip
   {\bf Continuing Demonstration for $[v_1^{\prime},\,e_2^{\prime}]$.}
   Recall that $\varphi:[v_1^{\prime},e_2^{\prime}]\mapsto[d_4,r_1]$.
   Choose a compatible basis, say, $(d_4, r_1, (0,1,0))$ for $N$
   and its dual basis for $M$. Extend the transverse basis to
   $[v_1^{\prime},e_2^{\prime}]$ to the full basis, say
   $(0,0,0,1,0)$, $(0,1,0,0,0)$, $(0,0,-1,-1,-1)$, $(0,0,2,2,3)$,
   $(1,0,0,-2,-3)$ for $N^{\prime}$ and take its dual basis for
   $M^{\prime}$. Then (treating vectors as column vectors)
   $$
    \varphi\;=\; \left[\,
                  \begin{array}{rrrrr}
                   0  & 0  &  0  & 0  & -1 \\
                   0  & 0  & -1  & 2  &  0 \\
                   0  & 1  &  0  & 0  &  0 
                  \end{array}\,
                 \right]
   $$
   and the induced map on quotients
   $\overline{\varphi}_{d_4r_1}^{v_1^{\prime}e_2^{\prime}}:
     N^{\prime}/
      \mbox{\raisebox{-.4ex}{$N^{\prime}_{v_1^{\prime}e_2^{\prime}}$}}
     \rightarrow
     N/\mbox{\raisebox{-.4ex}{$N_{d_4r_1}$}}$
   is given by
   $\overline{\varphi}_{d_4r_1}^{v_1^{\prime}e_2^{\prime}}=(0, 1, 0)$
   and the dual
   ${\overline{\varphi}_{d_4r_1}^{v_1^{\prime}e_2^{\prime}}}^{\dagger}:
     [d_4,r_1]^{\perp}\rightarrow{v_1^{\prime}e_2^{\prime}}^{\perp}$
   is given by
   ${\overline{\varphi}_{d_4r_1}^{v_1^{\prime}e_2^{\prime}}}^{\dagger}
     =(0, 1, 0)^{t}$.
   From these, one obtains that
   $$
    \overline{\Delta^{\prime}_{[v_1^{\prime},e_2^{\prime}]}
     \cap M^{\prime}}
     =\{(0,0),(1,0),(2,0),(3,0),(3,-1)\}\,,
   $$
   as indicated in {\sc Figure 4-2-3}, up to an integral affine
   transformation. This describes the line bundle
   ${\cal L}|_{F^{v_1^{\prime}e_2^{\prime}}_{d_4r_1}}$. The set of
   cones in $\Star([v_1^{\prime},\,e_2^{\prime}])$ that lie in 
   ${\overline{\varphi}^{v_1^{\prime}e_2^{\prime}}_{d_4r_1}}^{-1}(0)$
   are given by
   $\bar{e}_3^{\prime}\bar{e}_1^{\prime}$, 
   $\bar{e}_3^{\prime}\bar{v}_4^{\prime}$, 
   $\bar{e}_1^{\prime}\bar{v}_4^{\prime}$,
   and their faces. Since
   $\bar{e}_1^{\prime}\rightarrow (1,3)$,
   $\bar{e}_3^{\prime}\rightarrow (0,-1)$,
   $\bar{v_4}^{\prime}\rightarrow (-1,0)$
   in the $2$-dimensional description,
   this fan describes $\WCP^2(1,1,3)$.
   Now choose
   $\xi_{d_4r_1}^{v_1^{\prime}e_2^{\prime}}:
     N/\mbox{\raisebox{-.4ex}{$N_{d_4r_1}$}} \rightarrow
     N^{\prime}/
     \mbox{\raisebox{-.4ex}{$N^{\prime}_{v_1^{\prime}e_2^{\prime}}$}}$
   to be, say, $(0,1,0)^t$.
   This induces a rewriting of $s|_{V([v_1^{\prime},\,e_2^{\prime}])}$
   as
   {\footnotesize
   \begin{eqnarray*}
    \lefteqn{
      s|_{V([v_1^{\prime},\,e_2^{\prime}])}\;
        =\; a_{(0,0,0)}\,+\,a_{(3,0,-1)}t_1^3t_3^{-1}\,
         +\,a_{(1,1,0)}t_1t_2\,+\,a_{(1,2,0)}t_1t_2^2  } \\
     && \hspace{2em}  
         +\,a_{(2,2,0)}t_1^2t_2^2\,,+\,a_{(2,3,0)}t_1^2t_2^3\,
         +\,a_{(2,4,0)}t_1^2t_2^4\,+\,a_{(3,3,0)}t_1^3t_2^3\, \\
     && \hspace{2em}
         +\,a_{(3,4,0)}t_1^3t_2^4\,+\,a_{(3,5,0)}t_1^3t_2^5\,
         +\,a_{(3,6,0)}t_1^3t_2^6\,. \\
     && =\; a_{(0,0,0)}\,
            +\,\left(\,a_{(1,1,0)}t_2\,+\,a_{(1,2,0)}t_2^2\,\right)\,t_1\,\\
     &&    \hspace{2em}   
            +\,\left(\,a_{(2,2,0)}t_2^2\,+\, a_{(2,3,0)}t_2^3\,
                               +\, a_{(2,4,0)}t_2^4\,\right)\,t_1^2 \\
     &&     \hspace{2em}
            +\,\left(\,a_{(3,3,0)}t_2^3\,+\,a_{(3,4,0)}t_2^4\,
            +\,a_{(3,5,0)}t_2^5\,+\,a_{(3,6,0)}t_2^6\,\right)\,t_1^3\,\\
     &&     \hspace{2em}
            +\,a_{(3,0,-1)}t_1^3t_3^{-1} \\
     && =\; b_{(0,0)}(t_2)\, +\,b_{(1,0)}(t_2)\,t^1\,
            +\,b_{(2,0)}(t_2)\,t_1^2\, +\,b_{(3,0)}(t_2)\,t_1^3\, \\
     &&     \hspace{2em}
            +\,b_{(3,-1)}(t_2)\,t_1^3t_3^{-1}\,.
   \end{eqnarray*}
   {\normalsize This}} 
   is an expression of a section in
   ${\cal L}|_{F^{v_1^{\prime}e_2^{\prime}}_{d_4r_1}}$
   with coefficient functions $b_{(\,\cdot\,,\,\cdot\,)}(t_2)$
   depending on the base parametrized by $t_2\in{\Bbb C}^{\times}$.
   \newline$\mbox{\hspace{9cm}}\Box$
   \bigskip
  \end{itemize}
 \end{itemize}
\end{itemize}

\begin{itemize}
 \item
 {\it The induced fibration $Y^{\prime}\rightarrow X_{\Sigma}$
 and the stratification of the discriminant locus of the fibration}$\,$:
 Once the section $s$ is $X_{\Sigma}$-orbit-wise written in the form
 of sections of ${\cal L}_{F^{\tau^{\prime}}_{\sigma}}$ with
 coefficients functions of $O_{\sigma}$, using the discriminant
 associated to $\overline{\Delta^{\prime}_{\tau^{\prime}}}$, one
 morally can study the stratification of the discriminant locus
 of the fibration. Unfortunately, for the current example, the
 number of parameters involved are so huge that it is above any
 tool known to us to analyze it.
\end{itemize}

\begin{itemize}
 \item
 {\it Beyond}$\,$:
 Here we select only {\it the issue of desingularization} of this
 example for some discussions.
 \begin{quote}
  \hspace{-1.9em}(1)\hspace{1ex}
  Recall that
  {\footnotesize
  $$
   \Sing(X_{\Sigma^{\prime}})\;
   =\; V([v_5^{\prime}b^{\prime}])\,
      \sqcup\, V([v_4^{\prime}b^{\prime}])\,
      \sqcup\, V([v_4^{\prime}e_1^{\prime}e_2^{\prime}])\;
   \simeq\; X_{\Sigma}\sqcup X_{\Sigma}\,
           \sqcup\, \mbox{\rm ${\footnotesizeBbb C}$P}^1
                     \times\mbox{\rm $\footnotesizeBbb C$P}^1\,.
  $$
  {\normalsize To}} 
  desingularize $X_{\Sigma^{\prime}}$, one adds
  {\small
   $$
    \begin{array} {l}
     b_3^{\prime}\;=\;(v_5^{\prime}+b^{\prime})/2\;
                   =\;(0,0,0,1,1)\,,     \\[.6ex]
     b_1^{\prime}\;=\;(v_4^{\prime}+2\,b^{\prime})/3\;
                   =\;(0,0,0,1,2)\,,\quad
     b_2^{\prime}\;=\;(2\,v_4^{\prime}+b^{\prime})/3\;
                   =\;(0,0,0,0,1)\,,     \\[.6ex]
     e_4^{\prime}\;=\;(v_4^{\prime}+e_1^{\prime}+e_2^{\prime})/3\;
                   =\;(0,0,1,1,2)
    \end{array}
   $$
   {\normalsize to}} 
  $\Sigma^{\prime}(1)$. Since these lattice points lie in
  ${\Delta^{\prime}}^{\ast}$, they determine a unique refined
  triangulation of ${\Delta^{\prime}}^{\ast}$, which in turn determines
  a refinement $\widetilde{\Sigma^{\prime}}$ of $\Sigma^{\prime}$.
  The following rules tell us how a simplicial cone in
  $\Sigma^{\prime}$ splits to simplicial cones in
  $\widetilde{\Sigma^{\prime}}$:
  $$
   \begin{array}{rcl}
    [\,\cdots,\, v_5^{\prime},\, b^{\prime}\,]
     & \longrightarrow
     & [\,\cdots,\,v_5^{\prime},\, b_3^{\prime}\,]\,
       +\, [\,\cdots,\, b_3^{\prime},\, b^{\prime}\,]      \\[.6ex]
    [\,\cdots,\, v_4^{\prime},\, b^{\prime}\,]
     & \longrightarrow
     & [\,\cdots,\, v_4^{\prime},\, b_2^{\prime}\,]\,
       +\, [\,\cdots,\, b_2^{\prime},\, b_1^{\prime}\,]\,
       +\, [\,\cdots,\, b_1^{\prime},\, b^{\prime}\,]      \\[.6ex]
    [\,\cdots,\, v_4^{\prime},\, e_1^{\prime},\, e_2^{\prime}\,]
     & \longrightarrow
     & [\,\cdots,\, e_1^{\prime},\, e_2^{\prime},\, e_4^{\prime}\,]\,
       +\, [\,\cdots,\, v_4^{\prime},\, e_2^{\prime},\,
                                        e_4^{\prime}\,] \\
    && \hspace{8em} 
        +\, [\,\cdots,\, v_4^{\prime},\, e_1^{\prime},\, e_4^{\prime}\,]
                                                          \\[.6ex]
   \end{array}
  $$
  All other types of cones remain the same.
  Denote $\widetilde{X_{\Sigma^{\prime}}}$ the toric variety associated
  to the refined fan. One can check that all the new maximal cones have
  multiplicity $1$, thus $\widetilde{X_{\Sigma^{\prime}}}$ is a smooth
  $5$-fold. Since $Y^{\prime}$ is $\Sigma^{\prime}$-regular, its proper
  transform $\widetilde{Y^{\prime}}$ in $\widetilde{X_{\Sigma^{\prime}}}$
  is also smooth.

  \medskip
  \hspace{-1.9em}(2)\hspace{1ex}
  Since
  $\varphi:b_3^{\prime}, b_1^{\prime}, b_2^{\prime}\,\rightarrow\, 0,\;
    e_4^{\prime}\,\rightarrow\, r_1$,
  the set of primitive cones over $[\cdots,r_1]$ with respect to
  $\varphi$ now becomes
  $$ 
   \{\, [\cdots^{\prime},\, e_1^{\prime}],\,
        [\cdots^{\prime},\, e_2^{\prime}],\,
        [\cdots^{\prime},\, e_3^{\prime}],\,
        [\cdots^{\prime},\, e_4^{\prime}]\,\}\,.
  $$
  The primitive cones over all other types of cones in $\Sigma$
  remain the same.

  \medskip
  \hspace{-1.9em}(3)\hspace{1ex}
  Following the same argument in earlier parts of the example, one
  has the following proper transform/resolution of components of
  fibers of $\widetilde{\varphi}\,$:
  {\footnotesize
  $$
   \hspace{-5em}
   \begin{array}{ccl}
    \footnotesizeWCP^2(1,2,3) & \longleftarrow
      & \widetilde{\footnotesizeWCP^2(1,2,3)}\;
        =\;X_{\{(1,1),\,(2,3),\,(1,2),\,(0,1),\,(-1,0),\,(0,-1)\}}
          \\[.6ex]
     && \simeq\;
        \mbox{\rm {\footnotesizeBbb C}P}^2\,
        \sharp\,3\,\overline{\mbox{\rm {\footnotesizeBbb C}P}^2}
      \\[2ex]
    X(4)   & \longleftarrow
      & \widetilde{X(4)}\;
        =\;X_{\{(1,1),\,(2,3),\,(1,2),\,(0,1),\,(-1,0),\,
                                           (-1,-1),\,(0,-1)\}}
            \\[.6ex]
     && \simeq\;
        \mbox{\rm {\footnotesizeBbb C}P}^2\,
        \sharp\,4\,\overline{\mbox{\rm {\footnotesizeBbb C}P}^2}
         \\[2ex]
   \mbox{\rm {\footnotesizeBbb C}P}^2  & \longleftarrow
      & \mbox{\rm {\footnotesizeBbb C}P}^2    \\[2ex]
    X(5)  & \longleftarrow
      & \widetilde{X(5)}\;
        =\;X_{\{(1,1),\,(2,3),\,(1,2),\,(0,1),\,(-1,0),\,(-1,-1),\,
                                 (-2,-3),\,(-1,-2),\,(0,-1)\}}
          \\[.6ex]
     && \simeq\;
        \mbox{\rm {\footnotesizeBbb C}P}^2\,
        \sharp\,6\,\overline{\mbox{\rm {\footnotesizeBbb C}P}^2}
          \\[2ex]
    \footnotesizeWCP^2(1,1,3)   & \longleftarrow
      & \widetilde{\footnotesizeWCP^2(1,1,3)}\;
        =\;X_{\{(1,3),\,(0,1),\,(-1,0),\,(0,-1)\}}
          \\[.6ex]
     && \simeq\;
        {\footnotesizeBbb F}_3
         \\[2ex]
    {\footnotesizeBbb F}_2   & \longleftarrow
      & {\footnotesizeBbb F}_2        \\[2ex]
    \bullet & \longleftarrow & \mbox{\rm {\footnotesizeBbb C}P}^2
   \end{array}
  $$
  {\normalsize The}} 
  intersection relations of these components are induced from the
  original components and are indicated in {\sc Figure 4-2-4}.

  \medskip
  \hspace{-1.9em}(4)\hspace{1ex}
  The restriction of the pullback line bundle on
  $\widetilde{X_{\Sigma^{\prime}}}$, still denoted by ${\cal L}$,
  to these fibers are still described by the same polytopes, except
  for the new fiber component $\CP^2$ from the relative stars
  $\Star_{[\cdots,\,r_1]}{[\cdots^{\prime},\,e_4^{\prime}]}$,
  the restriction of ${\cal L}$ to which is trivial.
  Note that none of the restriction is the anticanonical bundle
  except for the generic fiber
  {\footnotesize $\widetilde{\footnotesizeWCP^2(1,2,3)}$}.
  \begin{quote}
   \hspace{-1.9em}(4.1)\hspace{1ex}
   {\it For $F_{\sigma}=\widetilde{X(4)}\cup\CP^2\,$}:
   With labels in {\sc Figure 4-2-4},
   {\footnotesize
   $$
   \hspace{-6em}
   \begin{array}{l}
    {\cal L}_{c_1^{\prime}}={\cal O}(D_1+D_2+D_4+E_1+E_2+E_3) \\[.6ex]
    \hspace{6ex} \mbox{and} \hspace{2ex}
         (D_1+D_2+D_4+E_1+E_2+E_3)\cdot D_3=2
          \hspace{2ex} \mbox{on $\widetilde{X(4)}$}\,;  \\[2ex]
    {\cal L}_{c_2^{\prime}}={\cal O}(D_5+D_6)   \\[.6ex]
    \hspace{6ex} \mbox{and} \hspace{2ex}
         (D_5+D_6)\cdot D_3=2
          \hspace{2ex} \mbox{on $\mbox{\rm $\footnotesizeBbb C$P}^2$}\,.
   \end{array}
   $$
   } 

   \hspace{-1.9em}(4.2)\hspace{1ex}
   {\it For $F_{\sigma}
        =\widetilde{X(5)}\cup\widetilde{\WCP^2(1,1,3)}
         \cup{\Bbb F}_2\cup\CP^2\,$}:
   With labels in {\sc Figure 4-2-4},
   {\footnotesize
   $$
   \hspace{-6em}
   \begin{array}{l}
    {\cal L}_{e_1^{\prime}}={\cal O}(D_1+D_2-D_3+D_5+E_1+E_2+E_3) \\[.6ex]
      \hspace{6em}\mbox{and}\hspace{2ex}
        \left\{
          \begin{array}{rr}
           (D_1+D_2-D_3+D_5+E_1+E_2+E_3)\cdot D_3=1 \\
           (D_1+D_2-D_3+D_5+E_1+E_2+E_3)\cdot D_4=0
          \end{array}
         \right.
        \hspace{2ex} \mbox{on $\widetilde{X(5)}$}\,;    \\[2ex]
    {\cal L}_{e_2^{\prime}}={\cal O}(-D_3+D_6+D_7)   \\[.6ex]
      \hspace{6em}\mbox{and}\hspace{2ex}
         \left\{
          \begin{array}{rr}
           (-D_3+D_6+D_7)\cdot D_3=1 \\
           (-D_3+D_6+D_7)\cdot D_7=3
          \end{array}
         \right.
        \hspace{2ex} \mbox{on
         $\widetilde{\mbox{\rm W${\footnotesizeBbb C}$P}^2(1,1,3)}$}\,;
                                                            \\[2ex]
    {\cal L}_{e_3^{\prime}}={\cal O}(-D_7+D_8+D_9)  \\[.6ex] 
      \hspace{6em} \mbox{and}\hspace{2ex}
         \left\{
          \begin{array}{rr}
           (-D_7+D_8+D_9)\cdot D_4=0 \\
           (-D_7+D_8+D_9)\cdot D_7=3
          \end{array}
         \right.           
        \hspace{2ex} \mbox{on ${\footnotesizeBbb F}_2$}\,;
                                                           \\[2ex]
    {\cal L}_{e_4^{\prime}}={\cal O}
        \hspace{2ex} \mbox{on $\mbox{{\footnotesizeBbb C}P}^2$}\,.
   \end{array}
   $$
   } 
  \end{quote}

  \medskip
  \hspace{-1.9em}(5)\hspace{1ex}
  Recall the adjunction formula $g=1+\frac{1}{2}(K_S+C)\cdot C$ that
  computes the genus $g$ of a smooth embedded curve $C$ in the divisor
  class $[C]$ of a smooth surface $S$ via the intersection number of
  the canonical divisor $K_S$ of $S$ and the self-intersection number
  of $C$.
  One can check that all the smooth curves in
  $F^{\widetilde{\tau^{\prime}}}_{\sigma}$ cut out by sections of
  ${\cal L}_{F^{\widetilde{\tau^{\prime}}}_{\sigma}}$, where
  $\widetilde{\tau^{\prime}}$ is any primitive cone in
  $\widetilde{\Sigma^{\prime}}$, are rational, as they should be.
  Since each type of fibers appear over a codimension-$1$ stratum
  of the discriminant locus of the fibration
  $\widetilde{Y^{\prime}}\rightarrow X_{\Sigma}$, from the intersection
  numbers computed in Sub-item (4) above and the Kodaira's table
  of singular elliptic fibers of a $1$-parameter family elliptic
  fibration, one concludes that$\,$:
  {\it
  For the elliptic fibration 
  $\widetilde{Y^{\prime}}\rightarrow X_{\Sigma}$ of Calabi-Yau
  $4$-fold induced from the toric morphism
  $X_{\widetilde{\Sigma^{\prime}}}\rightarrow X_{\Sigma}$,
  the generic degenerate fiber associated to the
  $\widetilde{X(4)}\cup\CP^2$ toric fiber has two irreducible
  components and is of type $\mbox{\it III}$ (or $\widetilde{A_1}$)
  in Kodaira's table
  while the generic degenerate fiber associated to the
  $\widetilde{X(5)}\cup\widetilde{\WCP^2(1,1,3)}\cup{\Bbb F}_2\cup\CP^2$
  toric fiber has five irreducible components
  (one in $\widetilde{X(5)}$, one in $\widetilde{\WCP^2(1,1,3)}$ and
   three in ${\Bbb F}_2$) and is of type $I_0^{\ast}$
  (or $\widetilde{D_4}$). }
 \end{quote}
\end{itemize}
\vspace{-1.5em}
\begin{figure}[htbp]
 \setcaption{{\sc Figure 4-2-4.}
  \baselineskip 14pt
  The toric divisors associated to $1$-cones in the toric fiber
   of $\widetilde{X_{\Sigma^{\prime}}}\rightarrow X_{\Sigma}$
  and the Kodaira labels for the generic fiber and the generic
   degenerate fibers of the elliptic Calabi-Yau hypersurface
   $Y^{\prime}$ after toric resolution of the singularities.
 } 
 \centerline{\psfig{figure=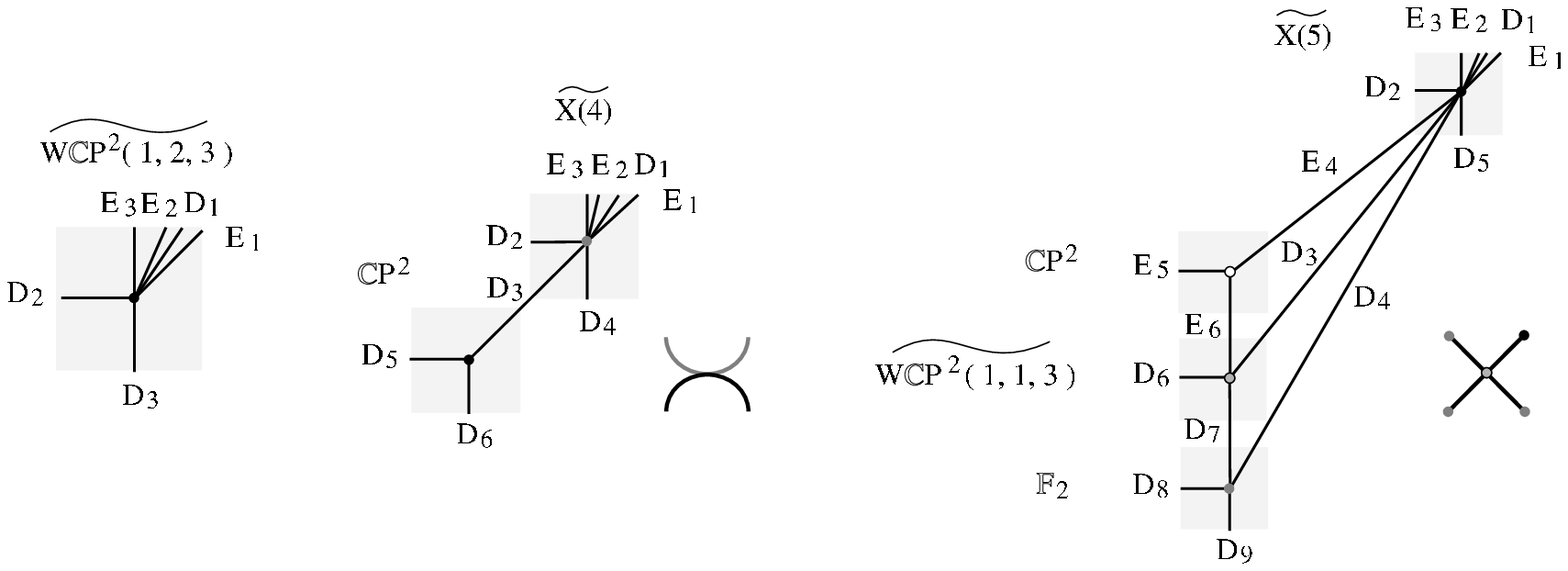,width=11cm,caption=}}
\end{figure}

\bigskip

This concludes our long discussions of the example.

\noindent\hspace{12cm} $\Box$

\bigskip

\section{Remarks on further study.}

At the end, we list five problems from our list for further pursuit.

\bigskip

\noindent 
{\bf Problem 5.1 [link to Kreuzer-Skarke].}
 Elaboration of the work of Kreuzer and Skarke ([K-S1] and [K-S2]),
 which classifies higher dimensional reflexive polytopes that admit
 fibration, with the computational scheme presented in this paper.
 This would give further information of fibrations in their
 list.

\bigskip

\noindent
{\bf Problem 5.2 [fibred Calabi-Yau variety].}
 The above huge number of concrete examples may provide us a hint
  for deeper understanding of general fibred Calabi-Yau varieties
  - toric and non-toric alike -.
 These examples may also shed some light on the generalization of
  Kodaira's table of elliptic fibrations to higher dimensional
  families.
 
\bigskip

\noindent 
{\bf Problem 5.3 [complete intersection].}
 Extension of the current work to higher rank equivariant vector
  bundles over a toric variety to take care of fibred Calabi-Yau
  varieties from complete intersections of hypersurfaces in a toric
  variety and a toric morphism.
 Realization of this will complete the first step toward
  understanding fibred Calabi-Yau variety in the toric category.

\bigskip

\noindent 
{\bf Problem 5.4 [effect on quantum cohomology].}
 In view of the goal that understanding a fibration is meant to
  understanding the original space in terms of lower dimensional ones
  and the data how they pile together, one would like to work out 
  a generalized quantum K\"{u}nneth formula for the quantum cohomology
  of fibred toric variety (cf.\ [Ba1] and [K-M]) and also for
  the hypersurfaces therein. 
 This problem currently looks extremely difficult.

\bigskip

\noindent 
{\bf Problem 5.5 [Het/F-theory duality].}
 Calabi-Yau hierarchy from toric morphisms:
  a purely toric realization of the Het/F-theory duality.
  This will be the analogue of the work of Batyrev on mirror
  symmetry of Calabi-Yau manifolds of the same dimensions.
 This problem is a final dream of all the participants of the
  heterotic-string/F-theiry duality in the toric category. 

\bigskip

With these problems in mind for further investigations, let us
conclude this paper and turn to another new pursuit.


\end{document}